\documentclass[11pt]{amsart}
\usepackage{amsmath,amssymb,amscd,fullpage,pb-diagram,hyperref}
\thanks{The first named author was supported  by NSA
Grant MSPF-08Y-172 and NSF Grant DMS-1001792 during the period this work was completed.
The second named author was supported by ISF Grant 1138/10
during the period this work was completed.}
\theoremstyle{plain}

\newcommand{\iq}[1]{\int\limits_{#1\quo}}
\newcommand{\wt}{\widetilde}
\newcommand{\bpm}{\begin{pmatrix}}
\newcommand{\epm}{\end{pmatrix}}
\newcommand{\bsm}{\begin{smallmatrix}}
\newcommand{\esm}{\end{smallmatrix}}
\newcommand{\bspm}{\left(\begin{smallmatrix}}
\newcommand{\espm}{\end{smallmatrix}\right)}
\newcommand{\bm}{\begin{matrix}}
\newcommand{\bbm}{\begin{bmatrix}}
\newcommand{\ebm}{\end{bmatrix}}
\newcommand{\on}{\operatorname}

\newcommand{\bs}{\backslash}
\newcommand{\iFA}{\int\limits_{F\bs \A}}
\newcommand{\jquo}{G^J(F) \bs \wt G^J(\A)}
\newcommand{\mquo}{SL_2(F) \bs \wt{SL}_2(\A)}

\newcommand{\cH}{\mathcal{H}}

\newtheorem*{theorem*}{Theorem}
\newtheorem{lemma}{Lemma}[section]

\newtheorem{proposition}[lemma]{Proposition}

\newtheorem{theorem}[lemma]{Theorem}
\newtheorem{definition}[lemma]{Definition}

\newtheorem{corollary}[lemma]{Corollary}

\newcommand{\ccc}{\mathbb{C}}
\newcommand{\iso}{\cong}
\usepackage[usenames]{color}

\newtheorem{thm}[equation]{Theorem}
\newtheorem{cor}[equation]{Corollary}
\newtheorem{cortopf}[equation]{Corollary-to-the-Proof}
\newtheorem{prop}[equation]{Proposition}
\newtheorem{lem}[equation]{Lemma}
\newtheorem{rmk}[equation]{Remark}
\newtheorem{rmks}[equation]{Remarks}

\newtheorem{defn}[equation]{Definition}
\numberwithin{equation}{subsection}

\newcommand{\C}{\mathbb{C}}
\newcommand{\G}{\mathbb{G}}
\newcommand{\A}{\mathbb{A}}

\newcommand{\Z}{\mathbb{Z}}
\newcommand{\R}{\mathbb{R}}
\newcommand{\N}{\mathbb{N}}
\newcommand{\schwartz}{\mathcal{S}}

\newcommand{\Ad}{\operatorname{Ad}}
\newcommand{\Invt}{\operatorname{Invt}}
\newcommand{\Vol}{\operatorname{Vol}}

\newcommand{\Gal}{\operatorname{Gal}}
\newcommand{\diag}{\operatorname{diag}}
\newcommand{\pr}{\operatorname{pr}}
\newcommand{\Ind}{\operatorname{Ind}}
\newcommand{\Aut}{\operatorname{Aut}}
\newcommand{\Hom}{\operatorname{Hom}}
\newcommand{\Inn}{\operatorname{Inn}}
\newcommand{\inc}{\operatorname{inc}}
\newcommand{\sym}{\operatorname{sym}}

\newcommand{\osc}{{\omega_{\baseaddchar}}}
\newcommand{\tpp}{{\theta_{\phi}^{\baseaddchar}}}

\newcommand{\SL}{\widetilde{SL}}
\newcommand{\maxunip}{U_{\max}}
\newcommand{\descgroup}[1]{N_{#1}}
\newcommand{\desclevi}[1]{L_{#1}}

\newcommand{\baseaddchar}{\psi_0}
\newcommand{\uniper}{\mathcal{U}}
\newcommand{\Res}{\operatorname{Res}}
\newcommand{\res}{\operatorname{res}}
\newcommand{\quo}{(F \backslash \A)}
\newcommand{\residuerep}{\mathcal{E}_{-1}(\tau,\omega )}

\newcommand{\primeideal}{\mathfrak{p}}

\newcommand{\middlestageparabolic}{P_1}
\newcommand{\descentparabolic}[1]{Q_{#1}}
\newcommand{\middlestagechar}{\hat\mu\delta_{\middlestageparabolic}^{\frac12}}

\renewcommand{\Re}{\operatorname{Re}}

\newcommand{\mainglnrep}{\tau}
\newcommand{\somerep}{\pi}

\newcommand{\ourgroup}{G}
\newcommand{\poleloc}{\underline{\mathbf{\frac{1}{2}}}}
\newcommand{\normal}{$ is a normal subgroup of $}
\newcommand{\FJ}{F\!J^{\baseaddchar}}
\newcommand{\Sieg}{U_{\on{Si}}}
\newcommand{\FC}[2]{^{(#1, \psi_{#1,#2})}}
\begin{document}

\author{Joseph Hundley} \email{jhundley@math.siu.edu}
\address{Math. Department, Mailcode 4408,
Southern Illinois University Carbondale, 1245 Lincoln
Drive Carbondale, IL 62901}

\author{Eitan Sayag}
\email{eitan.sayag@gmail.com}
\address{Institute of Mathematics, The Hebrew University of Jerusalem,
Jerusalem, 91904, Israel}


\title[Descent Construction for GSpin groups]{Descent Construction for GSpin groups}

\date{\today}

\maketitle
\tableofcontents

\begin{abstract}

In this paper we provide an extension of the theory of {\sl descent}
of Ginzburg-Rallis-Soudry to the context of essentially self-dual
representations, that is representations which are
isomorphic to the twist of their own contragredient
 by some Hecke character.
Our theory supplements the recent work of Asgari-Shahidi on
the functorial lift from (split and quasisplit forms of) $GSpin_{2n}$  to $GL_{2n}.$
\end{abstract}

\section{Introduction}

A fundamental problem in modern number theory is the determination of the spectral decomposition of the right regular representation of $G(\A)$ on the function space
$L^{2}(G(F) \backslash G(\A))$, where $G$ is a reductive $F$-group and $\A$ is the adele ring of the global field $F.$
 The constituents of the decomposition are special representations of $G(\A)$ that are called automorphic representations and encode analytic, geometric and number theoretic information.
The {\it Langlands program}, a vast generalization of class-field theory, offers a description
of these constituents in terms of certain homomorphisms from the conjectural Langlands group $L_{F}$ to a complex group $^LG$, the $L$-group associated with $G.$
Schematically, one expects a meaningful correspondence between the two sets  $$ \{\pi \text{ automorphic representation of } G(\A) \} \longleftrightarrow \{\text{admissible }\phi: L_{F} \to \;^L\!\,G\}$$
and sometimes it is possible to make this expectation exact.

\subsection{Functoriality}

One approach towards making the Langlands philosophy an exact predicition is via a conjectural relation called {\it functorial transfer}.
Namely,
whenever we have a homomorphism $r:\,^L\!H \to \,^L\!\,G$ one should expect,
in view of the correspondence above,  a relation between
automorphic representations on $H(\A)$ and those on $G(\A).$ To give a more precise description of the exepected correspondence we recall that 
an automorphic representation $\pi$ admits a presentation as $\pi=\otimes \pi_{v}$ where each $\pi_{v}$ is an irreducible admissible representation
of $G(F_{v})$. Furthermore for almost all places $v$ of $F$ , the subgroup  $K_{v}=G(\frak{o}_{v}) \subset G(F_{v})$  is a maximal compact subgroup and 
 the representations $\pi_{v}$ admit a fixed vector with respect to $K_{v}$. Representations $\rho$ of $G(F_{v})$ with this property are called {\it unramified} and to each one may attach a semisimple conjugacy class $S(\rho) \subset {}\,^L\!\,G$ called the Satake parameters of $\rho.$
For an automorphic representation $\pi=\otimes \pi_{v}$ we denote $S_{v}(\pi)=S(\pi_{v}).$ The ordered collection $(S_{v}(\pi))$ of Satake parameters serves as a fingerprint of the automorphic representation $\pi$ and is used to compare automorphic representations. Given $r:\,^L\!H \to \,^L\!\,G$ as above and $\pi=\otimes \pi_{v}$  an
 automorphic representations on $G(\A)$ the {\it functoriality conjecture} predicts the existence of an  automorphic representations $\sigma$ on $H(\A)$ 
which is a weak lift of $\pi$, meaning
$$r(S_{v}(\pi))=S_{v}(\sigma)$$ holds for almost all places $v.$

%

The problem of establishing the existence of functorial lifts relating automorphic representations on different groups attracted a lot of research and several methods have been used to establish such a relation.
Roughly speaking, one can classify the attempts to establish the functoriality conjecture into three major approachs: the method of the trace formula, the method of  the converse theorem supplemented with a theory of $L$-functions, and various methods of explicit constructions of automorphic representations.
In recent years, all these methods have been developed and refined and many new instances of functoriality were established.

In the present paper we focus on the method of explicit constructions by the method of descent. We extend the descent method of Ginzburg, Rallis, and
Soudry to GSpin groups and are able to obtain information about the image of the functorial lift constructed in \cite{Asg-Sha1} and \cite{Asg-Sha3}. Below we provide more information.

\subsection{Self dual representations of $GL_{n}(\A)$ and the descent method}

Given a classical group $H$ (e.g. $Sp(2n),SO(2n+1),SO(2n)$) the dual group $^LH$ is a classical group as well  (e.g. $SO(2n+1,\C), Sp(2n,\C), SO(2n,\C)$) and naturally embedded in a general linear group. Now this standard embedding $^LH \to 
 GL_{N}(\C) ={}\; ^LGL_{N}$
 should conjecturally yield a functorial lift from each automorphic representation $\sigma$ of $H(\A)$ to an automorphic representation $\tau$ of $GL(N,\A).$ 
Analyzing the Satake parameters of these lifts one expect $\tau$ to be self dual, that is $\tilde{\tau}\cong \tau.$

In
\cite{Cog-K-PS-Sha} the method the of converse theorem is used to show
the existence of a weak functorial lift from generic cuspidal
automorphic representations of classical groups to automorphic
representations of the general linear group.

The theory of descent for self-dual cuspidal representations of the
general linear group $GL_{n}(\A)$ was developed in a sequence of
remarkable works \cite{GRS1}-\cite{GRS5}.
The definitive account of this work is now
available in \cite{GRS-Book}.
In these works the authors
showed, among other things, that
every self-dual irreducible cuspidal automorphic
representation $\pi$ of $GL_n(\A)$ is a weak functorial lift
from a cuspidal representation of some classical group.  
Moreover, this was accomplished by explicitly 
constructing 
 a space $\sigma(\pi)$ of cuspidal
automorphic functions  which weakly lifts to
a  $\pi.$
Thus the construction of $\sigma(\pi)$ from $\pi$
reverses the {\it lifting} constructed in \cite{Cog-K-PS-Sha}, 
and for this reason is known as the ``descent.''
The precise classical group on which $\sigma(\pi)$ 
is defined is governed by poles of the symmetric 
and exterior square $L$ functions
of $\pi.$  For example,  if $L(\pi,\wedge^2,s)$ has a pole at $s=1$ then $\sigma(\pi)$ is defined 
on $SO_{2n+1}(\A).$

Among other things, 
the authors of \cite{GRS4} were able to use their descent construction 
to describe the
image of the functorial lift of \cite{CKP-SS}.

Thus, the conjunction of the descent method and associated integral
representations with the method of the
converse theorem provides a very detailed description of the image of
functoriality corresponding to the standard embedding of $^LG \to
GL_{N}(\ccc)$ with $G$ a classical group. For an excellent survey we
refer the reader to \cite{So-Paris}.

The recent proof of the fundamental Lemma (and its twisted versions) by \cite{Ngo}, and the work of Arthur on the trace formula and its application to functoriality is expected to lead
to a complete proof of the functoriality conjecture for 
the case when $r$ is the standard embedding 
of any classical group.  
We refer the reader to \cite{Arthur-Book} for the expected results.
The method of the trace formula has the advantage that 
it is not restricted to generic representations of classical groups.   On the other hand,  it is indirect and one of the advantage of the descent method is that it provides an explicit realization of the automorphic representation constructed as a space of functions.

\subsection{Beyond classical}
Observe that each classical group is contained in a corresponding 
similitude group.  One is naturally led to consider the case 
when $r$ is the standard embedding of a similitude 
group.
 However, the dual group of a 
similitude classical group is not another similitude classical group, but rather a new type of group called a GSpin group.  Thus in order to consider the case when 
$r$ is the standard embedding of a similitude classical group, one must deal with automorphic forms defined
on GSpin groups.
Recently,
Asgari and Shahidi proved in \cite{Asg-Sha1} 
and \cite{Asg-Sha3} the existence of a weak
functorial lifting from each quasisplit GSpin group to  the appropriate general linear group.
The representations 
 obtained in this fashion are essentially self dual. 
Moreover, in the special case of $GSp(4)$ they were able to show in
\cite{Asg-Sha2} that this weak functorial lift is in fact strong in
an appropriate sense.

\subsection{Descent construction for essentially self dual representations}

In this paper we extend the descent method of Ginzburg, Rallis, and
Soudry to GSpin groups. As a bonus, for $n\ge 2$  we can provide a ``lower bound''
on the image of each of the functorial lifts constructed by  Asgari and Shahidi.   For the case of $GSpin_5 = GSp_4,$ these results were obtained
by another method in \cite{Gan-Tak}.
Since this work was completed, Asgari and Shahidi have shown that this ``lower bound'' is, in fact, the full image \cite{Asg-Sha3}.



We now describe a few properties of GSpin groups.  For more details see section \ref{Notation: Gspin}. 
There is a unique quasisplit form of $GSpin_{2n+1},$ 
while
quasisplit forms of $GSpin_{2n}$ are in natural one-to-one correspondence
with quadratic characters $\chi:\A^\times/F^\times \to \pm 1.$ We will denote the corresponding
group by  $GSpin_{2n}^\chi.$ We also denote by $F[\chi]/F$ the quadratic extension corresponding to
$\chi.$

We then obtain  
maps
$$r: {}\;^L(GSpin_{2n}^\chi)=GSO_{2n}(\C)\rtimes \Gal(F[\chi]/F) \rightarrow GL_{2n}(\C) ={}\; ^LGL_{2n}.$$
$$r: {}\;^L(GSpin_{2n+1})=GSO_{2n}(\C) \rightarrow GL_{2n}(\C) ={}\; ^LGL_{2n}.$$
sending the nontrivial element of $\Gal(F[\chi]/F)$ to
$$\begin{pmatrix}I_{n-1}&&&\\ &&1&\\ &1&&\\ &&&I_{n-1}\end{pmatrix}.$$
In each of the theorems below the notion of weak lift is defined relative to the appropriate map $r$ as above.

To discern the form of $GSpin_{2n}$  to which a given representation $\tau$
 will descend, we observe that $\tau \iso \tilde \tau \otimes \omega$
 implies $\omega_\tau^2=\omega^{2n}.$  Here $\omega_\tau$ denotes the
 central character of $\tau.$  Hence $\omega_\tau/\omega^n$ is some
 quadratic character $\chi.$



Let us now state the main results of the present paper. For simplicity of exposition we consider the odd case and the even case separately. The results are analoguous to the ones of \cite{GRS4} describing endoscopic lifts from orthogonal and symplectic groups. In our formulation we  use the notion of  $\bar{\omega}$- symplectic and  $\bar{\omega}$-orthogonal  representations which is introduced in \ref{def:ess-symp}

We begin with main result of part 1 of this work. 
\begin{theorem*}[MAIN THEOREM: ODD CASE]\label{oddt:main}
For $r \in \N,$ take $\tau_1, \dots, \tau_r$ to be irreducible unitary
automorphic cuspidal representations of $GL_{2n_1}(\A), \dots, GL_{2n_r}(\A),$
respectively, and let $\tau=\tau_1\boxplus \dots \boxplus \tau_r$
be the isobaric sum (see section \ref{s:descrOfTau}).
 Let $\omega$ denote a Hecke character.
Suppose that
\begin{itemize}
\item $\tau_i$ is $\bar{\omega}$- symplectic
 for each $i,$ and
 \item
 $\tau_i \iso \tau_j \Rightarrow i = j.$
 \end{itemize}
Then there exists an irreducible generic cuspidal
automorphic representation $\sigma$ of
$GSpin_{2n+1}(\A)$ such that
\begin{itemize}
\item
$\sigma$ weakly lifts to $\tau,$  and
\item
the central character  of $\sigma$ is $\omega.$
\end{itemize}
\end{theorem*}

Next, we formulate the main result of part 2. 

\begin{theorem*}[MAIN THEOREM: EVEN CASE]\label{t:main-even}
For $r \in \N,$ take $\tau_1, \dots, \tau_r$ to be irreducible unitary
automorphic cuspidal representations of $GL_{2n_1}(\A), \dots, GL_{2n_r}(\A),$
respectively, and let $\tau=\tau_1\boxplus \dots \boxplus \tau_r$
be the isobaric sum (see section \ref{s:descrOfTau}).
Let $n=n_1+\dots +n_r,$ and assume that $n\ge 2.$
 Let $\omega$ denote a Hecke character, which is not the square of another Hecke character.
Suppose that
\begin{itemize}
\item $\tau_i$ is $\omega^{-1}$-orthogonal
 for each $i,$ and
 \item
 $\tau_i \iso \tau_j \Rightarrow i = j.$
 \end{itemize}
 For each $i,$ let $\chi_i=\omega_{\tau_i}/\omega^{n_i}$ (which is quadratic),
 and let $\chi=\prod_{i=1}^r \chi_i.$
Then there exists an irreducible generic cuspidal
automorphic representation $\sigma$ of
$GSpin_{2n}^\chi(\A)$ such that
\begin{itemize}
\item
$\sigma$ weakly lifts to $\tau,$  and
\item
the central character $\omega_{\sigma}$ of $\sigma$ is $\omega.$
\end{itemize}
\end{theorem*}

\subsection{Applications}

In a series of works, Asgari and Shahidi studied the functorial lifts from Spinor groups to the general linear group. We first formulate their result
for the odd case.
The next theorem is Theorem 1.1 (p. 140) of \cite{Asg-Sha1}.
\begin{thm}\label{AS}
Let $\pi$ be an irreducible globally generic cuspidal automorphic representation of $GSpin_{2n+1}(\A).$ Then there exists
 an automorphic representation $\tau=AS(\pi)$ of $GL_{2n}(\A)$  which is a weak lift of $\pi$ with respect to the map $r$ of \eqref{functoriality}. Furthermore, denoting by $\omega_{\pi}, \omega_{\tau}$ the central characters
 of $\pi,\tau$ we have $\omega_{\tau}=\omega_{\pi}^{n}.$
Finally, $AS(\pi)$ and $\widetilde{AS(\pi)}\otimes \omega_\pi$ are nearly equivalent. 
\end{thm}



The present work allows a description of the image of the lifting constructed by Asgari-Shahidi. Thus we have the following

\begin{corollary}
The  image
of the functorial lift $AS$ (see Theorem \ref{AS}) contains the set of all representations $\tau_1 \boxplus \dots \boxplus \tau_r $
such that
\begin{itemize}
\item $\tau_i \iso \tau_j \Rightarrow i = j,$
\item there is a Hecke character $\omega$ such that
$\tau_i$ is $\bar{\omega}$- symplectic
 for each $i.$
 \end{itemize}
\end{corollary}

The above result of the present work was used  in \cite{Asg-Rag}.
\begin{rmk}
Since this work was completed, it has
 been shown in \cite{Asg-Sha3} that the image of the functorial lift $AS$
 is actually equal to the set described above.
 \end{rmk}

For the even case, the results are similar, though the 
sequence of events is somewhat different.
The next theorem is
Theorem 5.16 of \cite{Asg-Sha3}.
\begin{thm}\label{ASeven}
Let $\pi$ be an irreducible globally generic cuspidal automorphic
representation of $GSpin_{2n}^{\chi}(\A).$ Then there exists
 an automorphic representation $\tau=AS_{\chi}(\pi)$ of $GL_{2n}(\A)$  which is a weak lift of $\pi$ with respect to the map $r$ of
 \eqref{functorialityeven}. Furthermore, denoting by $\omega_{\pi}, \omega_{\tau}$ the central characters
 of $\pi,\tau$ we have $\omega_{\tau}=\omega_{\pi}^{n}\chi.$
Finally, $AS_{\chi}(\pi)$ and $\widetilde{AS_{\chi}(\pi)}\otimes \omega_\pi$ are nearly equivalent. 
\end{thm}
The contribution of this work is a lower bound 
for the image of this functorial lift.
\begin{corollary}
The  image of the functorial lift $AS_{\chi}$ described in Theorem
\ref{ASeven} contains the set of all representations $\tau_1
\boxplus \dots \boxplus \tau_r $ such that
\begin{itemize}
\item $\tau_i \iso \tau_j \Rightarrow i = j,$
\item there is a Hecke character $\omega$ such that
$\tau_i$ is $\omega^{-1}$- orthogonal
 for each $i.$
 \end{itemize}
\end{corollary}

As with the odd case, it is
 shown in \cite{Asg-Sha3} that the image of the functorial lift $AS_\chi$
 is actually equal to the set described above.

\subsection{The descent construction and the structure of the argument}

We give here a brief outline of the ingredients one use to prove the main theorems. We refer the reader to \cite{HS} where we have provided a 
detailed description for the scheme of the argument for the even case. 

The construction of the descent representations
relies on the notion of Fourier coefficient, as defined in
 \cite{MR1981592}, \cite{MR2214128}  (cf. also the ``Gelfand-Graev'' coefficients of
\cite{So-Paris}).
    For purposes
 of presenting certain of the global arguments, it seems convenient to embed
 these Fourier coefficients into a slightly larger family of functionals,
 which we shall refer to as ``unipotent periods.''

  Suppose that $U$ is a unipotent subgroup of $G$ and
$\psi$ is a character of $U(F) \backslash U(\A).$  We define the
corresponding {\it unipotent period} to be the map from smooth,
left $U(F)$-invariant functions on $G(\A)$ to smooth,
left $(U(\A),\psi)$-equivariant functions,  given by
$$\varphi\mapsto \varphi^{(U,\psi)}$$
where
$$\varphi^{(U,\psi)}(g) := \int_{U(F)\backslash U(\A)} \varphi( ug ) \psi(u) \; du.$$



Let $\mathcal S$ be a set of
unipotent periods.  We will say that another unipotent period
$(U,\psi)$ {\it is spanned by} $\mathcal S$ if the implication
\begin{equation}\label{spanning of uniper}
\left( \varphi^{(N,\vartheta)}\equiv 0 \;\;\; \forall (N,\vartheta) \in \mathcal S \right)
\quad \Longrightarrow \quad \varphi^{(U,\psi)} \equiv 0
\end{equation}
is valid for any automorphic function $\varphi.$


 For simplicity of the exposition
we assume that we are trying to construct a descent for an
$\bar \omega$-orthogonal 
 cuspidal
representation, $\tau$ of the group $GL_{2n}(\A).$
 Recall that the central character of $\tau$ is 
 equal to $\omega^n \chi_\tau$ for some quadratic character
 $\chi_\tau.$ 
 
 We can conveniently describe the
method in the following steps:

\begin{enumerate}
\item Construction of a family of descent representations 
of $GSpin^{\chi}_{4n-2\ell}(\A)$ for $\ell \geq n,$
and $\chi$ any quadratic character.

\item Vanishing of the descent representations
for all $\ell > n$ and all $\chi \ne
 \chi_{{\tau}}.$


\item Cuspidality and genericity (hence nonvanishing) of the descent representation of
$GSpin_{2n}^{\chi_{ \tau}}(\A)$.

\item Matching of spectral parameters at unramified places.
\end{enumerate}

%

{For the construction of the descent representations, we begin, in section \ref{s:EisensteinSeries}, with an Eisenstein series on the group $GSpin_{4n+1}$ that is induced from
 a Levi $M$ isomorphic to $GL_{2n} \times GL_{1}$. The representations $\tau$, properly twisted, is used
 as the inducing data and a pole of
 $L(s, \tau, \sym^2 \times \bar\omega)$
  allows us to construct a residue
 representation $\mathcal{E}_{-1}(\tau,\omega)$ of $GSpin_{4n+1}(\A).$}
  Next,  for each $\ell$ and $\chi$ we give 
  in section \ref{section:  vanishing of deeper descents (even)}
  an embedding of
  $GSpin_{4n-2\ell}^{\chi}$ into $GSpin_{4n+1},$ and
construct, using formation of a Fourier coefficient, a space of
functions  on this subgroup of $GSpin_{4n+1}.$

 Let $DC_{\omega}^\chi(\tau)$ denote the space of functions on the quasisplit form $GSpin_{2n}^\chi$ of $GSpin_{2n}$ corresponding to the character $\chi.$
 Then, using a local argument involving twisted Jacquet modules
which is given in section \ref{localresults} ,  we
prove that $DC^\chi_\omega(\tau)$ is
zero, except {when $\chi=\omega_\tau/\omega^n.$}
A similar argument is used to prove the vanishing of the 
``deeper descents,'' defined on $GSpin_{4n-2\ell}^\chi$
for all $\ell > n$ and all $\chi.$

If $\chi = \chi_{\tau},$  
then  
 we show that $DC^\chi_\omega(\tau)$ is
nonzero, even generic, and all of the functions in
it are cuspidal.  
This is based on relations of unipotent periods of the 
type given in \eqref{spanning of uniper}, which 
are proved in section \ref{s: appV}.
It follows that $DC^\chi_\omega(\tau)=\oplus \sigma_{i}$
decomposes as a direct sum of
irreducible automorphic cuspidal representations $\sigma_{i}$ of $GSpin_{2n}^\chi$.
We then show,
using again a local argument in section \ref{localresults}, that each of these irreducible constituents $\sigma_i$ lifts
weakly to $\mainglnrep$ by the functorial lifting 
$AS_{\chi_\tau}.$

To treat the general case, the representation $\tau$ may be an isobaric
sum of  several cuspidal representations
$\tau_1, \dots, \tau_r$ and one considers a more complicated Eisenstein series.
{The main differences in the construction are that the
 residue is a multi-residue, and the notation is more complicated.  But the idea of the construction
 is similar.} See sections \ref{s:even} - \ref{s: appV}.
 
 The $\bar \omega$-symplectic case is similar, and
 is covered in sections \ref{oddt:main} to \ref{s:appII}.

\subsection{Structure of the paper}

In constructing this paper, we have tried to emphasize self-containment and readability over efficiency.  Therefore we write out in detail many arguments on Eisenstein series which may be regarded as standard and well known to the experts.  Also, we have,  written up the odd and even cases separately,
and in full detail, even though there is a great overlap between the two.  At each step, the statements and proofs in the odd and even cases are variations of one another.  However, we felt it was desirable to give complete details for both cases, and that would be easier to follow one argument at a time than two variations developed in parallel.  The only significant deviation from this principle is the treatment of Eisenstein series, for which we have
included a complete and detailed proof only in the even case.  The proof in the odd case is omitted because it is virtually identical, {\it and} because both proofs consist mainly of specializing general results on Eisenstein
series to a particular case.

If one is interested in reading only the odd case,
and is willing to take as given some standard
facts about Eisenstein series, this case is contained
in parts 1 and 2.
For a self-contained
treatment of the even case, one should read parts 1 and 3.

We now quickly summarize the logic of the paper. In section
\ref{sec: formulation} we formulate the main result in the odd case.
In sections
\ref{sec: notations}
and \ref{sec: four}
we set up some general notions,
 and introduce
$GSpin$
groups.  We collect here information about root datum,
Weyl group, and unipotent subgroups.
Section \ref{s:uniper} gives some general material 
concerning
 and unipotent periods which may be useful in other contexts. 
With these general matters completed, we begin the 
detailed treatment of the odd case.
 Section \ref{notation and statement for odd case}
 contains the precise statement of the main theorem in 
 the odd case, and a little bit of additional notation.
 In section
\ref{sec: unramified} we consider the Satake parameters of
essentially self-dual unramified representations on $GL_{2n}.$ We
show that these are always parameters associated with
representations of smooth representations of $GSpin$ groups. This
is the local underpinning of the descent construction. In section
\ref{sec: Eis} we introduce a specific family of Eisenstein
series, and formulate their main properties, in Theorem \ref{thm:Eis}.
In section \ref{s: main, odd} we prove the main
result for the odd case.  Sections  \ref{oddlocalresults}
and \ref{s:appII} are appendices to the odd case, containing
proofs of various technical results.
In the remaining sections we repeat all
of these steps for the ``even'' case (descent from $GL_{2n}$ to some
quasisplit form of $GSpin_{2n}$).

\subsection{Acknowledgements}
The authors wish to thank the following people for helpful
conversations: Dubravka Ban, William Banks, Daniel Bump, Wee Teck
Gan, Herv\'e Jacquet, Erez Lapid, Omer Offen, Yiannis
Sakellaridis,  Gordan Savin, Freydoon Shahidi, and Christian Zorn.  In addition, we wish to thank
Mahdi Asgari, Jim Cogdell, Anantharam Raghuram and Freydoon
Shahidi for their interest, which stimulated the work.

Without
 David Ginzburg and David Soudry's many careful and patient explanations of the
``classical'' case-- $\omega =1$--  this work would not have been completed.
It is important to point out that not all of the arguments shown to
us have appeared in print.   We mention in particular the computation of
Jacquet modules in Appendix V,  and the nonvanishing argument in
Appendix VI.
Nevertheless, in each case the specialization of
our arguments to the case $\omega=1$ may be correctly attributed to
Ginzburg, Rallis, Soudry (with any errors or stylistic blemishes introduced being our
own responsibility).

This
work was undertaken while both authors were in Bonn at the
Hausdorff Research Institute for Mathematics, in connection with a
series of lectures of Professor Soudry's.  They wish to thank the
Hausdorff Institute and Michael Rapoport for the opportunity. 
Finally, the second author wishes to thank Prof. Erez Lapid for many
enlightening discussions on the subject matter of these notes.

\part{General matters}
\section{Some notions related to Langlands functoriality}\label{s:main statement odd}\label{sec: formulation}

To formulate the main result of the paper we need the
terminology of weak lift, the notion of essentially self dual
representations and the notion of isobaric sums. We review them
here briefly but the reader can go directly to \ref{oddt:main}.

\subsection{Weak Lift}
\label{ss:weak lift}
Let $G$ and $H$ be reductive groups defined over a number field
$F.$ Let $\pi \iso \otimes_v' \pi_v$ be an automorphic
representation of the group $G(\A),$ where $\A$ is the ring of
adeles of $F.$

At almost all places $v$ the smooth representation $\pi_v$ of
$G(F_{v})$ is unramified. The semi-simple conjugacy
class in the $L$-group $^LG$ associated by the Satake isomorphism
to $\pi_v$ will be denoted $t_{\pi_v}.$

We say that an automorphic representation $\sigma$ of $G(\A)$ is a
{\it weak lift} of the
automorphic representation $\tau\iso \otimes'_v\tau_v$ of $H(\A)$ with respect to a map $r:^LH \to ^LG$  if for almost all
places, $r(t_{\sigma_{v}})\subset t_{\tau_{v}}$.

We now specialize to the case $H=GSpin_{2n+1}$ (for definition
and properties see \cite{AsgariCanJMath} and section \ref{Notation: Gspin})
$G=GL_{2n}$ and the inclusion
\begin{equation}\label{functoriality}
r: \,^L\!H={}^L\!(GSpin_{2n+1})=GSp_{2n}(\C) \hookrightarrow GL_{2n}(\C) ={}\; ^L\!GL_{2n}={}^L\!G.
\end{equation}

\subsection{Essentially self-dual: $\eta$-orthogonal and $\eta$-symplectic}

To formulate our main
result we introduce the notion of $\eta$-orthogonal and $\eta$-symplectic
 representations.
Let $\mainglnrep$ be an
irreducible automorphic cuspidal representation of $GL_{2n},$
and let $\eta$ be a Hecke character.
Recall that the theory of self dual representations is well
developed. Suppose that $\mainglnrep$ is {\it essentially
self-dual}, i.e. that

{
$$\tilde\mainglnrep=\mainglnrep \otimes \eta$$ Where $\tilde\mainglnrep$
is the contragredient of $\mainglnrep$ and $\eta$ is some Hecke character.}

\begin{defn}
Let $\mainglnrep$ be an irreducible automorphic cuspidal
representation of $GL_{2n}$
and $\eta$ a Hecke character.
Consider the (exterior) tensor product representation $\tau
\boxtimes \eta $  of the group
$GL_{2n}(\A) \times GL_1(\A).$
The finite Galois form of the
$L$ group of $GL_{2n}(\A) \times GL_1(\A)$ is the group $GL_{2n}(\C) \times GL_1(\C).$

We write
 $L( s, \mainglnrep, \wedge^2 \times \eta)$
 ( resp.  $L(s,
\mainglnrep, sym^2 \times \eta)$)
 for the $L$ function attached to $\tau
\boxtimes \eta $ and the representation of
$GL_{2n}(\C) \times GL_1(\C)$
where $GL_{2n}(\C)$ acts by the exterior (resp. symmetric) square
representation  and $GL_1(\C)$ acts by scalar
multiplication.
\end{defn}


\begin{lemma}
Let $\mainglnrep$ be an irreducible automorphic cuspidal
representation of $GL_{2n}$
and $\eta$ a Hecke character.
 Consider the $L$ functions $L( s,
\mainglnrep, \wedge^2 \times \eta)$ and $L(s, \mainglnrep, sym^2
\times \eta).$
\begin{enumerate}
\item If $\tilde\mainglnrep$ is isomorphic to $\mainglnrep \otimes \eta$ for some
Hecke character $\eta$ then exactly one of these $L$-functions has
a simple pole at $s=1$ while the other is holomorphic and
nonvanishing.
\item
if $\tilde \tau$ is not isomorphic to $\tau \otimes \eta$ then
both functions are holomorphic at $s=1.$
\end{enumerate}
\end{lemma}
\begin{proof}
It follows from \cite{MW2}, the corollaire on p.667,
that
 $L(s, \mainglnrep \times \mainglnrep \otimes \eta)$ has a simple pole
at $s=1$ if $\mainglnrep \otimes \eta = \widetilde \tau$
and is holomorphic at $s=1$ otherwise.

Now, $L( s, \mainglnrep \times \mainglnrep \otimes
\eta)$ is the Langlands $L$ function of on $M_{2n\times 2n}(\C)$ in which $GL_{2n}(\C)$ acts by
$g \cdot X = g X \;^tg$
 and
$GL_1(\C)$ acts by scalar multiplication.
But this representation is reducible, decomposing into the subspaces
of  skew-symmetric and symmetric matrices.
Thus  we have
$$L( s, \mainglnrep \times \mainglnrep
\otimes \eta)=L( s, \mainglnrep, \wedge^2 \times \eta)L(s,
\mainglnrep, sym^2 \times \eta).$$

It is proved in \cite{Shahidi-PacJMath}, theorem 4.1, p.
321, that
neither
$L(s, \tau, \wedge^2 \times \eta)$
 nor $L(s, \tau, \wedge^2 \times \eta)$
vanishes on the line $\Re(s)=1.$  The result follows

\end{proof}

We are ready for a definition

\begin{definition}\label{def:ess-symp}
\begin{enumerate}
\item We will say that $\tau$ is {\it $\eta$-symplectic} in case $L( s,
\mainglnrep, \wedge^2 \times \eta)$ has a pole at $s=1.$

\item We will say that $\tau$ is {\it $\eta$-orthogonal} in case $L(s,
\mainglnrep, sym^2 \times \eta)$ has a pole at $s=1.$
\end{enumerate}
\end{definition}

Thus, if $\tau$ is essentially self dual, there exists a Hecke
character $\eta$ such that $\tau$ is either $\eta$-symplectic or
$\eta$-orthogonal.
\begin{rmks}\begin{enumerate}
\item Above we have defined ``$\eta$-symplectic'' and 
``$\eta$-orthogonal'' using completed $L$ functions.  We 
could also have used partial $L$ functions.  The notions
obtained are equivalent.   
(See item {\bf(ii):}
near the end of the proof of propostion \ref{relRk1toLFcn-GSpin}, and its proof.)
\item{There is another natural notion of ``orthogonal/symplectic representation.''  Specifically,
one could say that an automorphic representation is orthogonal/symplectic if the space
it acts on supports an invariant symmetric/skew-symmetric form.  The two notions
appear to be related, but do not coincide.  See \cite{PrasadRamakrishnan}.
}
\item{ There is a third approach to defining a local factor for
$L(s, \tau, \wedge^2 \times \eta),$ which is to apply the ``gcd'' construction
described in \cite{GelbartShahidi} section I.1.6, p. 17, to the integrals in
\cite {JS}.  As far as we know this is not written down anywhere.}
\item{
An integral representation for
$L(s, \tau, sym^2)$ was given in \cite{BG}.  The problem of
extending this to $L(s, \tau, sym^2\times \eta)$ has been solved by Banks \cite{Banks-GL3,Banks-Thesis} in
the case $n=3$ and Takeda \cite{Takeda-twistedSym2}
for general $n.$
}
\end{enumerate}
\end{rmks}

\subsection{Isobaric sums}

We recall the following result.

\begin{theorem}
\begin{enumerate}

\item
 every irreducible automorphic representation of $GL_n(\A)$
is isomorphic to a subquotient of $Ind_{P(\A)}^{GL_n(\A)} \tau_1
|\det_1|^{s_1} \otimes \dots \otimes  \tau_r |\det_r|^{s_r}$ for
some real numbers $s_1, \dots s_r$ and irreducible unitary
automorphic cuspidal representations $\tau_1, \dots, \tau_r$ of
$GL_{n_1}(\A), \dots , GL_{n_r}(\A)$ respectively, such that $n_1+
\dots + n_r = n.$ Here $P$ is the standard parabolic of $GL_n$
corresponding to the ordered partition $(n_1, \dots, n_r)$ of $n.$

\item
In the case when $s_i=0$ for all $i,$ this induced representation
is irreducible.

\item
The representations obtained when $s_i=0$ for all $i$ by numbering a given set of
cuspidal representations in different ways are isomorphic.

\item
 If two induced representations as in (1), with $s_i=0$ for all $i,$ are isomorphic,
then they are obtained from two numberings of the same set of
cuspidal representations

\end{enumerate}

\end{theorem}
\begin{proof}
\begin{enumerate}
\item This follows from proposition 2 of \cite{Langlands-Notion}.

\item This follows from the irreducibility of all the local induced
representations, which is Theorem 3.2 of \cite{Jacquet-Generic}.

\item
This follows from the fact that the standard intertwining operator
between them does not vanish, which follows from \cite{MW1},
II.1.8 (meromorphically continued in IV.1.9(e)), and IV.1.10(b).
In IV.3.12 these elements are combined to prove that the
intertwining operator does not have a pole.  The proof that it
does not have a zero is an easy adaptation.
\item
This follows from \cite{JacquetShalika-EPandClassnII},
Theorem 4.4, p.809.
\end{enumerate}
\end{proof}

\begin{definition}
Let $\tau_1, \dots, \tau_r$ be irreducible unitary cuspidal representations
of $GL_{n_1}(\A), \dots, GL_{n_r}(\A).$
The
{\it isobaric sum} of  $\tau_1\boxplus
\dots \boxplus \tau_r$
of these representations is the
irreducible unitary representation $\tau$ of $GL_n(\A)$
(where $n=n_1+\dots+n_r$)
given in (1) above, with $s_i=0$ for all $i.$
\label{s:descrOfTau}
\end{definition}
\begin{rmk}
A more general notion of ``isobaric representation'' was
introduced in \cite{Langlands-Maerchen}, but we don't need
it.
\end{rmk}

\section{Notation}\label{s:notation}\label{sec: notations}
\subsection{General}\label{ss:GeneralNotation}
Throughout most of the paper, $F$ will denote a number field.

We denote by $J$ the matrix, of any size, with ones on the diagonal running from
upper right to lower left, and by $J'$ the matrix $\left( \begin{smallmatrix} &J \\-J&
\end{smallmatrix} \right)$ of any even size.
We also employ the notation $^tg$ for the transpose of $g$ and
$_tg$ for the ``other transpose'' $J\; ^t g J.$
We employ the shorthand $G\quo = G(F) \backslash G(\A),$ where $G$ is any
$F$-group.

We denote the Weyl group of the reductive group $G$ by $W_G$ or by $W,$
when the meaning  is clear from {the} context.
{Given $\pi$ a local representation or an automorphic
representation,  we denote by $\tilde \pi$ the contragredient of $\pi$, and
by $\omega_\pi$ its central character.}



\subsection{Notions of  ``genericity''}\label{s:NotionsOfGenericity}
Let $G$ be a quasisplit reductive group over the number field $F,$ and  $\maxunip$ a
maximal unipotent subgroup.
First let $\psi_v$ be a generic character
(cf. \cite{Kim-Survey}, p. 147, and also \cite{Shahidi-ParkCity}, p.304)
of $\maxunip(F_v)$ for some place $v$ of $F,$ and
 $(\pi_v,V)$ a representation of $G(F_v).$  We say that $\pi_v$ is $\psi_v$-generic if it supports
 a nontrivial $\psi_v$-Whittaker functional (i.e., a  $\maxunip(F_v)$-equivariant linear map
 $V \to \C_{\psi_v},$
 which is continuous in an appropriate topology, see \cite{Shahidi-ParkCity}, p. 304.
 Here $\C_{\psi_v}$ denotes the one-dimensional $\maxunip(F_v)$-module
 with action via the character $\psi_v.$
 )
 Now let
 $\pi \iso \otimes_v{}'\pi_v$ be a
 automorphic representation of $G(\A),$
 and let
 $\psi=\prod_v\psi_v$ be a  character of $\maxunip\quo,$
 by which we mean a character $\maxunip(\A)$
 which is trivial on $\maxunip(F).$

Ignoring topological considerations, it
 is easy to see that
 the space $\Hom_{\maxunip(\A)}(V_\pi, \C_\psi)$ is nontrivial iff each of the spaces
 $\Hom_{\maxunip(F_v)}(V_{\pi_v}, \C_{\psi_v})$ is.
   However, it turns out that the more important
 issue is not whether there exists {\it some}  nontrivial $\psi$-Whittaker functional,
 but whether the specific $\psi$-Whittaker functional
given by
$$\varphi \mapsto \int_{\maxunip\quo}\varphi(u)\bar\psi(u)\;du$$
is nonvanishing.
We refer to this Whittaker functional  as the $\psi$-Whittaker {\it integral.}
(See \cite{GelbartSoudry} for an example where the Whittaker integral vanishes,
but a nonzero Whittaker functional exists.)

We would like to take this opportunity to draw attention to the subtle fact  that there are two slightly different notions of global genericity
for automorphic representations in common usage.
The first states that a representation is globally $\psi$-generic
if it supports a nonzero $\psi$-Whittaker integral.
The second-- which was the notion originally introduced in
\cite{PS-MultOne}-- requires that
a cuspidal representation  be {\it orthogonal to the kernel} of the
$\psi$-Whittaker integral in $L^2_{\text{cusp}}(G\quo),$ in order to be called ``generic.''
Clearly, the latter condition implies the former (except for the zero representation).

A nice feature of the stronger formulation is that the condition  defines a subspace
of $L^2_{\text{cusp}}(G\quo),$ which one may term
the $\psi$-generic spectrum.  Furthermore, this subspace
 satisfies multiplicity one, even if $L^2_{\text{cusp}}(G\quo)$ does not.
 (Cf. \cite{PS-MultOne})
A nice feature of the weaker formulation is that it does not rely on the $L^2$-pairing, and hence
no technicalities arise in applying the notion to non-cuspidal forms and representations.

Throughout most of this paper, we shall say that a representation ``is $\psi$-generic''  if it supports
a nonzero $\psi$-Whittaker integral, and ``is generic'' if it
satisfies this condition for {\it some} $\psi.$  We shall say that a cuspidal representation
is ``in the $\psi$-generic spectrum''
if it is orthogonal to the kernel of the $\psi$-Whittaker integral.

Let $P_0 = N_G(\maxunip).$  If $P_0(F_v)$ permutes the characters of $U(F_v)$ transitively,
then we may refer to a representation as ``generic'' or ``non-generic'' without reference to a
specific $\psi_v,$ and without ambiguity.  The same applies to both notions of global genericity, in
the case when $P_0(F)$ permutes the characters of $\maxunip\quo$ transitively.  This condition
is satisfied by $GL_n$ and  $GSpin_{2n+1},$ but {\it not} by $GSpin_{2n}.$

\section{The Spin groups $GSpin_{m}$ and their quasisplit forms}
\label{sec: four}\label{Notation: Gspin}

We shall now define $GSpin$ groups {by introducing their root datum}. They will be defined as the
groups whose duals are the similitude classical groups
$GSp_{2n}(\C), GSO_{2n}(\C).$ Thus {we begin by describing} the based root
data of these classical groups, but employ notation appropriate to the role these groups will eventually
have as the duals of $GSpin$ groups.

\subsection{Similitude groups ${GSp_{2n},\;GSO_{2n}}$}

We first define the {\sl similitude orthogonal and symplectic
groups} to be
$$GO_m = \{ g \in GL_m: g J \;^t g = \lambda(g) J
\text{ for some }\lambda(g) \in \G_m \},$$
$$GSp_{2n} = \{ g \in GL_{2n}: g J' \;^t g = \lambda(g) J'
\text{ for some }\lambda(g) \in \G_m \}.$$
For each of these groups the map $g \mapsto \lambda(g)$ is a rational character
called the {\sl similitude factor.}

{
\begin{rmk}
If $m$ is odd then $GO_m$ is in fact isomorphic
to $SO_m \times GL_1.$  This case will play no further role.
\end{rmk}
}
The group $GO_{2n}$
is disconnected; indeed let $GSO_{2n}$ be the subgroup generated by $SO_{2n}$ and
$\left\{\left( \begin{smallmatrix} \lambda I_n & \\ & I_n \end{smallmatrix} \right) :
\lambda \in \G_m\right\}.$ Then $GSO_{2n}$ is connected and of index two in $GO_{2n}.$


The next lemma is straightforward.
\begin{lemma}

Let $$T=\{t=\diag(t_1, \dots, t_n, \lambda t_n^{-1} , \dots,
\lambda t_1^{-1})\}.$$

For $i = 1$ to $n,$ let $e_i^*(t)=t_{i}$ and $e_0^*(t)=\lambda.$

Let {$X^\vee=Hom(T,\G_m)$} denote the character lattice of $T$. Then
\begin{itemize}
\item $T$ is a maximal torus in both groups $GSp_{2n}$ and $GSO_{2n}$

\item $\{e_i^{*}: i = 0 $ to $n \}$ is a basis for $X^\vee$.
\end{itemize}
\end{lemma}
Let {$X=Hom(\G_m,T)$} be the cocharacter lattice of $T,$ and let $\{e_i: i = 0 $
to $n \}$ be the basis of $X$ dual to the basis $\{e_i^{*}: i = 0 $ to $n \}$ of $X^\vee.$

Each similitude classical group has a Borel subgroup equal to the
set of upper triangular matrices which are in it. In each case we
employ this choice of Borel, and
 let $\Delta^\vee$ denote the set of simple roots and $\Delta$
the set of simple coroots.

\begin{lemma}
\begin{enumerate}
\item For $GSp_{2n}$
$$\Delta^\vee = \{e_i^* - e_{i+1}^*, \; i = 1 \text{ to }n-1\} \cup \{2 e_n^*-e_0^*\}.$$
$$\Delta = \{ e_i-e_{i+1},\; i = 1 \text{ to }n-1\} \cup \{ e_n\}.$$
\item For $GSO_{2n}$
$$\Delta^\vee = \{e_i^* - e_{i+1}^*, \; i = 1 \text{ to }n-1\} \cup \{ e_{n-1}^*+e_n^*-e_0^*\}.$$
$$\Delta = \{ e_i-e_{i+1},\; i = 1 \text{ to }n-1\} \cup \{ e_{n-1}+e_n\}.$$
\end{enumerate}
\end{lemma}
\label{rootDataSection}

\subsection{Definition of split ${GSpin}$}

We first deal with the {\it split} forms.
By p.274 of \cite{Springer} 
$F$-split connected reductive algebraic groups are classified by
based root data. Thus the following makes sense.

\begin{definition}
\begin{enumerate}
\item The group $GSpin_{2n+1}$ is the $F$-split connected
reductive algebraic group having based root datum dual to that of
$GSp_{2n}.$

\item The group  $GSpin_{2n}$ is the $F$-split connected
reductive algebraic group having based root datum dual
to that of $GSO_{2n}.$
\end{enumerate}
\end{definition}

We note that an alternative description of the same group is used
in \cite{AsgariCanJMath}.

\begin{proposition}
The group $GSpin_m$ is identical to the quotient of $Spin_m \times
GL_1$ as given in \cite{AsgariCanJMath}.
\end{proposition}
This is Proposition 2.4 on p. 678 of \cite{AsgariCanJMath}.

\subsection{Definition of quasi-split ${GSpin}$}

Now we turn to {\it quasi-split} forms. By the classification
results in Chapter 16 of \cite{Springer} (especially 16.3.2,
16.3.3 16.4.2), one finds that $GSpin_{2n+1}$ is in fact the
unique {\it quasi-split } $F$-group having  based root datum dual
to that of $GSp_{2n}.$

Furthermore, there is a 1-1 correspondence between quasi-split $F$
groups $G$ such that $^LG^0 =GSO_{2n}(\C)$ and homomorphisms from
$\Gal( \bar F/ F)$ to the group $S$ of automorphisms of the
lattice $X(T)$ which permute the set $\Delta$ of simple roots
according to an automorphism of the Dynkin diagram.

\begin{lemma}
\begin{enumerate}
\item $S=\{1,\nu\}$ is of order two.

\item  $\nu(e_{n-1}-e_n)=e_{n-1}+e_n,
\nu(e_{n-1}+e_n)=e_{n-1}-e_n$ and $\nu ( \alpha) = \alpha$
if $\alpha$ is any of the remaining simple roots.

\item
$\nu(e_i)=\begin{cases}
e_i&i\ne 0,n\\
-e_n&i=n\\
e_0+e_n&i=0
\end{cases}
$

\item
$\nu(e_i^*)=\begin{cases}
e_i^*&i\ne n\\
e_0^*-e_n^*&i=n\\
\end{cases}
$

\item
lattices of $F$-rational characters is spanned by $\{e_i:0<i<n\}
\cup \{2e_0+e_n\}.$

\item the lattice of $F$-rational cocharacters is spanned by $\{e_i:0\le i <n\}.$

\end{enumerate}
\end{lemma}
\begin{proof}
One easily checks that the
formulae given in (4) define a nontrivial element $\nu$ of $S,$
and that the effect on the dual lattice and the roots in it
is as described in (2), (3).
   If $n \ne 4,$ then the Dynkin diagram has only two automorphisms,
   and (1) follows.   When $n=4,$ one must check that $S$ can
not contain any element of order $3.$  This follows from the
fact that such an element would induce an automorphism of $SO_8(\C)$
of order three, and no such automorphism exists.    This completes the proof of (1) in this case.  Items (5) and (6) are easily checked, given that an element
of $X$ or $X^\vee$ is $F$-rational if and only if it is $\Gal(\bar F/F)$-stable.
\end{proof}

By class field theory
homomorphisms from $\Gal( \bar F/ F)$ to a group
with two elements are in one-to-one correspondence with quadratic characters
$\chi: \A_F^\times / F^\times \to \{\pm 1\}.$
We denote the $F$-group corresponding to the character
$\chi$  by $GSpin_{2n}^\chi.$
The $F$-group corresponding to the trivial character is the unique split $F$-group
having the specified root datum, and is also denoted simply by
$GSpin_{2n}.$

To save space, the group $GSpin_m$ will usually be denoted $G_m.$

Observe that in either the odd or even case $e_0^*$ is a generator for the lattice of cocharacters
of the center of $\ourgroup_m.$

Because we define $\ourgroup_m$ in the manner we do, it comes equipped with a choice of Borel subgroup
and maximal torus,
as do various reductive subgroups we shall consider below.  In each case, we
denote the Borel subgroup of the reductive group $G$ by $B(G),$ and the maximal
torus by $T(G).$

A straightforward adaptation of the proof of Theorem 16.3.2 of \cite{Springer}
shows that there exist surjections $\ourgroup_m \rightarrow SO_m$
defined over $F.$  We fix one such and denote it $\pr.$  We require that
$\pr$ is such that $\pr( B( \ourgroup_m))$ consists of upper triangular matrices.

\begin{rmk}
For those familiar with the construction of
spin groups  as a subgroups of the multiplicative
groups of  Clifford algebras, we remark that our
$GSpin$ groups can be constructed in the same way by including the nonzero scalars in the Clifford
algebra as well.
However, this description will not play a role for us.
\end{rmk}

\subsection{Unipotent subgroups of $GSpin_m$ and their characters}
The kernel of
$\pr$ consists of semisimple elements.
In particular, the number of unipotent elements
 of a fiber is zero or one, and it's one if and only if the element of $SO_m$ is unipotent.
 In other words, $\pr$ yields a bijection of unipotent elements
 (indeed, an isomorphism of unipotent subvarieties), and we may specify unipotent
 elements or subgroups by their images under $\pr.$
 This also defines coordinates for any unipotent element or subgroup, which we use when
 defining characters.  Thus, we write $u_{ij}$ for the $i,j$ entry of $\pr(u).$

Above we fixed a specific isomorphism of a subgroup of $G_{2m}$ with $GL_m.$  If
$u$ is a unipotent element of of this subgroup this identification with an $m\times m$
matrix gives a second definition of $u_{ij}$ This is not a problem, however, as
the two definitions agree.

 Most of the unipotent groups we consider are subgroups of the maximal unipotent
 of $\ourgroup_m$ consisting of elements $u$ with $\pr(u)$ upper triangular.
 We denote this group $\maxunip.$
 A complete set of coordinates is $\{u_{ij}: 1 \leq i < j \leq m - i \}.$
We denote the opposite maximal unipotent by $\overline{\maxunip}.$  It consists of
all unipotent elements of $\ourgroup_m$ such that $\pr(u)$ is {\sl lower} triangular.

We fix once and for all a character $\psi_0$ of $\A/F.$  We use this character together
with the coordinates just above to specify characters of our unipotent subgroups.

When specifying subgroups of $\maxunip$ and their characters, the restriction to
$\{(i,j): 1 \leq i < j \leq m - i \}$ is implicit.

It will also sometimes be necessary to describe unipotent subgroups such that only a
few of the entries in the corresponding elements of $SO_m$ are nonzero.  For this
purpose we introduce the notation
$e'_{ij}= e_{ij}-e_{m+1-j,m+1-i}.$  One may check
that for all $i \neq j$ and $a \in F,$ the matrix $I+ae'_{ij}$ is an element of $SO_m(F),$ unless $m$ is odd and $i+j=m,$ in which case $I+ae'_{ij}+\frac{a^2}2 e_{i,m+1-i}
\in SO_m(F).$

\section{``Unipotent periods''}\label{s:uniper}
We now introduce the framework within which, we believe, certain of the computations involved in the descent construction can be most easily understood.

Let $G$ be a reductive algebraic group defined over a number field $F$.  If $U$ is a unipotent
subgroup of $G$
and $\psi_U$ is a character of $U\quo,$
(by which we mean a character of $U(\A)$ trivial on $U(F)$)
we define the {\sl unipotent period} $(U,\psi_U)$ associated
to this pair to be given by the formula
$$\varphi^{(U,\psi_U)}( g) := \int_{U\quo} \varphi(ug) \psi_U(u) du.$$
Clearly, $\varphi$ must be restricted to a space of
left $U(F)$-invariant locally integrable functions.

Let $\uniper$ denote the set of unipotent periods.  For $V$ a space of functions
defined on $G(\A),$ put
$$\uniper^\perp( V)=\{(U,\psi)\in \uniper: \varphi^{(U,\psi)} \equiv 0 \; \forall \varphi \in V\}.$$
When $V$ is the space of a representation $\somerep$ we will employ also the notation
$\uniper^\perp( \somerep).$  We also employ the notation $(U,\psi) \perp V$ for
$(U,\psi ) \in \uniper^\perp(V)$
and similarly $(U,\psi) \perp \somerep.$

We also require a vocabulary to express relationships among unipotent periods.
 We shall say that
$$(U,\psi_U)  \in \langle ( U_1, \psi_{U_1}) , \dots, (U_n, \psi_{U_n}), \dots \rangle$$
if $V \perp (U_i, \psi_{U_i}) \forall i \Rightarrow V \perp (U, \psi_U).$
Clearly, if
$(U_1,\psi_{U_1}) \in \langle (U_2, \psi_2), (U_3,\psi_3) \rangle,$ and
$(U_2, \psi_2) \in \langle (U_4,\psi_4), (U_5, \psi_5) \rangle$
then $(U_1, \psi_1) \in \langle (U_3, \psi_3),(U_4, \psi_4) ,(U_5,\psi_5)\rangle.$

We also use notation $(U_1, \psi_1) | (U_2, \psi_2),$  or the language
``$(U_1, \psi_1)$ divides $(U_2, \psi_2)$,'' `` $(U_2, \psi_2)$ is
divisible by $(U_1, \psi_1)$ ''
for $(U_2, \psi_2) \in \langle (U_1, \psi_1) \rangle.$
Finally, $(U_1, \psi_1 ) \sim(U_2, \psi_2)$
means $(U_1, \psi_1) | (U_2, \psi_2)$ and  $(U_2, \psi_2) | (U_1, \psi_1).$
This is an equivalence relation and we shall refer to unipotent periods which
are related in this way as ``equivalent.''

It is {\sl sometimes } possible to compose unipotent periods.  Specifically, if
$f^{(U_1,\psi_1)}$ is left-invariant by $U_2(F),$ then one may consider
$(f^{(U_1,\psi_1)})^{(U_2,\psi_2)}.$  We denote the composite
 by $(U_2, \psi_2) \circ (U_1, \psi_1).$

If $(U,U)$ is the commutator subgroup of $U,$
then $U/(U,U)$ is a vector group in the sense of
\cite{Springer}, p. 51.  The normalizer $N_G(U)$
of $U$ in $G$ acts on $U/(U,U).$  This is an
$F$-rational representation of $N_G(U).$  Denote
the dual $F$-rational representation of $N_G(U)$
by $(U/(U,U))^*.$  Then the choice
of $\baseaddchar$ identifies $F$ with
the space of characters of $F\backslash \A$ and thus
identifies $(U/(U,U))^*(F)$  with the space of
characters of $U(\A)$ which are trivial on $U(F).$
This gives this space the additional structure of
an $F$-rational representation of $N_G(U).$

Now, suppose that $U$ is the unipotent radical of a parabolic $P$ of $G$ with Levi $M.$
\begin{lem}
Then $(U/(U,U))^*$ is isomorphic, as an $F$-rational
representation of $M,$ to
$\overline U/(\overline U,\overline U),$
where $\overline{U}$ denotes the unipotent radical of the parabolic $\overline{P}$
of $G$ opposite to $P.$
\end{lem}
\begin{proof}
The exponential map induces an $M$-equivariant
isomorphism $\mathfrak{u}/[\mathfrak{u},\mathfrak{u}]
\to U/(U,U),$ defined over $F.$
The Killing form $\b B$ induces $M$-equivariant isomorphisms
$\mathfrak{u}^* \to \overline{\mathfrak{u}}$
and
$([\mathfrak{u},\mathfrak{u}])^*
\to [\overline{\mathfrak{u}},\overline{\mathfrak{u}}].$
If
$\mathfrak{v} := \{ X \in \overline{ \mathfrak{u}}
: \b B(X,Y) = 0 \; \forall \; Y \in [\mathfrak{u},\mathfrak{u}]\},$
then $\mathfrak v$
is canonically identified with $(\mathfrak{u}/[\mathfrak{u},\mathfrak{u}])^*,$
and maps isomorphically onto  $\overline U/(\overline U,\overline U).$
\end{proof}

This idea can be extended to more general unipotent groups.  See section \ref{the section with the extension}.

For $\vartheta$ a character of $U(\A)$
which is trivial on $U(F),$
 let $M^\vartheta$ denote the stabilizer of
$\vartheta$ (regarded as a point in $\overline{U}/(\overline{U},\overline{U})(F)$) in $M.$
So $M^\vartheta$ is an algebraic subgroup of $M$ defined over $F.$

\begin{defn}
\label{def:fc}
{In the set up of the previous paragraph, let}
$FC^\vartheta:  C^\infty( G\quo) \to C^\infty( M^\vartheta \quo )$ be given by
$$FC^\vartheta(\varphi)(m) = \varphi^{(U,\vartheta)} ( m )
= \int_{U\quo} \varphi( um ) \vartheta( u ) du.$$
\end{defn}
{
Clearly,
\begin{lemma}
The map $FC^\vartheta$ is an $M^\vartheta( \A)$-equivariant map.
\end{lemma}
}

\subsection{A Lemma Regarding Unipotent Periods}
There is a natural action of $G(F)$ on the space of unipotent periods
$\uniper$ given by $\gamma\cdot( U, \psi)
= (\gamma U \gamma^{-1}, \gamma \cdot \psi)$ where $\gamma\cdot \psi( u)
= \psi(\gamma^{-1} u \gamma).$
We shall refer to this action as ``conjugation.''  Obviously, unipotent periods which
are conjugate are equivalent.\label{ss:ConjOfUniper}

\begin{lem}  \label{uniperlemma}
Suppose $U_1 \supset U_2 \supset (U_1,U_1)$ are unipotent subgroups
of a reductive algebraic group $G$.  Suppose $H$ is a subgroup of $G$ and let $f$ be a smooth
left $H(F)$-invariant function on $G(\A).$  Suppose $\psi_2$ is a character of $U_2$ such
that $\left.\psi_2 \right|_{(U_1,U_1)}\equiv 0.$  Then the set $\res^{-1}(\psi_2)$ of characters
of $U_1$ such that the restriction to $U_2$ is $\psi_2$ is nontrivial.  (Here ``res'' is
for ``restriction'' not ``residue''.)
The
elements of $\res^{-1}(\psi_2)$ are permuted by the action of $N_H(U_1)(F).$
The  following are equivalent.
\begin{enumerate}
\item{ $f^{(U_2, \psi_2)}\equiv 0$}
\item{ $f^{(U_1, \psi_1)} \equiv 0 \; \forall \psi_1 \in \res^{-1}( \psi_2) $}
\item{ For each $N_H(U_1)(F)$-orbit $\mathcal{O}$ in $\res^{-1}( \psi_2)$
$\exists \psi_1 \in \mathcal{O}$ with $f^{(U_1, \psi_1)}\equiv 0 $}
\end{enumerate}
\end{lem}
\begin{proof}
It is obvious that 1 implies 2 and 3, and  that
2 and  3 are equivalent.
 Consider
$$f^{(U_2,\psi_2)}(u_1g) =
\int_{U_2\quo} f( u_2 u_1g) \psi_2(u_2) du_2,$$
regarded as a function of $u_1.$  It is left $u_2$ invariant and hence gives rise to
a function of the compact abelian group $U_2(\A)U_1(F) \backslash U_1(\A).$
Denote this function by $\phi(u_1).$  Then
$$\phi(0) = \sum_\chi \int_{U_2(\A)U_1(F) \backslash U_1(\A)}  \phi(u_1 ) \chi(u_1) du_1,$$
where ``$0$'' denotes the identity in $U_2(\A)U_1(F) \backslash U_1(\A),$
and the sum is over characters of
$U_2(\A)U_1(F) \backslash U_1(\A).$  This, in turn, is equal to
$$\sum_\chi \int_{D}  \int_{U_2\quo} f( u_2 u_1g) \psi_2(u_2) du_2
 \chi(u_1) du_1,$$
 for $D$ a fundamental domain for the above quotient in $U_1(\A).$
 The group  $U_1/(U_1,U_1)(F)$ is an $F$-vector space (cf. section \ref{s:uniper})
 which can be decomposed
 into $U_2/(U_1,U_1)(F)$ and a complement.  The $F$-dual of this vector space
 is identified, via the choice of $\psi_0,$ with the
  space of characters of $U_1(\A)$ which are trivial on $U_1(F).$  It follows that the
  sum above is equal to
 $$=\sum_{\psi_1 \in \res^{-1}(\psi_2) }
 \int_{U_1\quo} f(u_1 g) \psi_1(u_1) du_1.$$
 The matter of replacing the sum over $\chi$ by one over $\psi_1 \in \res^{-1}(\psi_2)$
 is clear from regarding $U_1/(U_1,U_1)(F)$ as a vector space which can be decomposed
 into $U_2/(U_1,U_1)$ and a complement.
Now 2 $\Rightarrow$ 1 is immediate.
\end{proof}
\begin{cor}  \label{unipercor}
 If $N_G(H)$ permutes the elements of $\res^{-1}( \psi_2)$ transitively,
then $(U_2, \psi_2) \sim (U_1, \psi_1)$ for every $\psi_1 \in \res^{-1}(\psi_2).$\end{cor}

\begin{defn}  Many of the applications of the above corollary are of a special type,
and it will be convenient to introduce a term for them.
The special situation is the following:  one has three unipotent periods $(U_i,\psi_i)$
for $i=1,2,3,$ such that $U_2 =U_1\cap U_3$ and $\psi_1 |_{U_2}=\psi_3|_{U_2} = \psi_2.$
Furthermore, $U_1$ normalizes $U_3$ and permutes transitively the set of characters $\psi_3'$
such that $\psi_3' |_{U_2}= \psi_2,$ and the same is true with the roles of $1$ and $3$ reversed.
In this situation, the identity
$$(U_1,\psi_1) \sim (U_2,\psi_2) \sim( U_3,\psi_3),$$
(which follows from Corollary \ref{unipercor}) will be called a {\bf swap,}
and we say that $(U_1, \psi_1)$ ``may be swapped for'' $(U_3, \psi_3),$
and vice versa.
\label{d:swap}
\end{defn}

\begin{rmk}
This notion of a ``swap'' is essentially an alternate
take on ``exchange of roots,'' as described
in \cite{GL}, section 2.1
(and references therein; see also \cite{GRS-Book}).
\end{rmk}

\subsection{Relation of unipotent periods  via theta functions}
\subsubsection{Initial remarks}
\label{the section with the extension}
We described above how a character of $U\quo,$  may be thought of as an element
of an $F$-vector space equipped with an algebraic action of $N_G(U).$
It was further shown how to identify such characters
with the $F$ points  of a distinguished subspace $\mathfrak{v}$  of $\mathfrak{g}$ in the special case when $U$ is the unipotent
radical of a parabolic subgroup.
This is advantageous, because the space $\mathfrak{g}_F$ is equipped with an action of $G(F),$ which is
 compatible with the action of $G(F)$ on unipotent periods by
conjugation (as in \ref{ss:ConjOfUniper}).

The same basic idea may be extended to arbitrary
unipotent subgroups using a suitable involution,
as we now explain.
First, every unipotent subgroup of $G$ is conjugate
to a subgroup of $\maxunip.$  So, assume $U \subset \maxunip.$
There is an automorphism $\iota$ of $G,$ unique up
to conjugation by an element of $T,$ which maps
$T$ to $T$ in such a fashion that $\iota(t)^\alpha
= t^{-\alpha}$ for all $t\in T$ and all roots $\alpha,$
and which maps the one-dimensional unipotent
subgroup $U_\alpha$ to $U_{-\alpha}.$
For $GSpin$ groups we can normalize $\iota$ by requiring that
the automorphism induced on $SO_m$ be transpose
inverse.
The automorphism $\iota$ induces an automorphism
of $\mathfrak g$ which we denote by the same letter.
Define $\b B'(X,Y) = \b B(X, \iota(Y))$
where $\b B$ is the Killing form.
Then the restriction of $\b B'$ to $\mathfrak{u}_{\max}$
is nondegenerate.
If $U$ is any unipotent subgroup of $G,$
define $\overline{U}= \iota(U).$ Note that if $U$ is
the unipotent radical of a parabolic subgroup,
then $\overline{U}$ is the unipotent radical of
the opposite parabolic, so this definition of $\overline{U}$
extends the previous usage.  Then $\b B$ identifies
the perp space of $[\mathfrak u, \mathfrak u]$
in
$\overline{\mathfrak{u}}$ with $(U/(U,U))^*.$
Thus, each character of $U(F)\backslash U(\A)$ is attached
to an element of $\overline{\mathfrak{u}}_F \subset\mathfrak g_F.$

Let $(U_1, \psi_1)$ and $(U_2, \psi_2)$ be
two unipotent periods and suppose that
the equivalence $(U_1, \psi_1)\sim (U_2, \psi_2)$
may be proved using conjugation and swapping.
Then it's fairly easy to see that $\psi_1$ and $\psi_2$ correspond to points in the same orbit of $G(F)$ acting on $\mathfrak g_F.$
The manner in which $U_1$ and $U_2$ are  related is not as easy to describe, but
one may note for example that they will have the same dimension.

The purpose of this section \label{section with theta functions} is to explain how the theory of the Weil representation may be used to obtain subtler relations among unipotent periods.
 This method was shown
to us by  David Ginzburg.

\subsubsection{Preliminaries on the Jacobi group
and Weil representation}
Let $\cH$ denote the Heisenberg group in three
variables, which we define as the set of triples
$(x,y,z)\in \G_a^3$ equipped with the product
$(x_1, y_1, z_1)(x_2, y_2, z_2) = (x_1+x_2, y_1+y_2,
z_1+z_2+x_1y_2-y_1x_2).$
The group is two-step nilpotent:  its center, $Z$  is equal
to its commutator subgroup, and consists
of all elements of the form $(0,0,z).$
The quotient $\cH/Z$ is a vector group.
The commutator defines a skew-symmetric bilinear
form $\cH/Z\times \cH/Z \to Z.$

The group $SL_2$ may be identified
with the group of all automorphisms of $\cH/Z$ which fix
this symplectic form, or, equivalently, with
the group of all automorphisms of $\cH$ which
restrict to the identity on $Z.$
The semidirect product
$\cH \rtimes SL_2$
of $\cH$ and $SL_2$ is sometimes called the Jacobi group.
Following \cite{Berndt-Schmidt}, we denote this
group by $G^J$ and identify it with
$\left\{ \bspm 1&*&*\\ &g_0 &*\\&&1 \espm : g_0 \in SL_2 \right\}
\subset Sp_4.$

The group $\cH(\A)$ has a unique
isomorphism class of
unitary representations such that the central
character is $(0,0,z) \mapsto \baseaddchar(z).$
Any representation in this class
extends uniquely to a projective
representation of $G^J(\A),$ or to a genuine
representation of
$\wt G^J(\A):=\cH(\A) \rtimes \SL_2(\A).$
Here $\SL_2(\A)$ denotes the
 metaplectic double cover of $SL_2(\A).$

There is a
representation $\osc$
known as the Weil representation
which is  in this class, given by action on the
space $\schwartz(\A)$ of Schwartz functions on $\A$ such that
\begin{equation}\label{oscillator action}
\osc(0,y,0)\phi(\xi ) = \baseaddchar(y\xi ) \phi(\xi),\qquad
\osc \begin{pmatrix}1&x \\ &1
\end{pmatrix}
\phi( \xi) = \baseaddchar(x\xi^2) \phi(x)
 \text{ and}\end{equation}
$$\osc(x,0,0)\phi(x)=\phi(\xi+x), \qquad (x,y, \xi \in \A, \;\; \phi \in \schwartz(\A)).$$
(We identify $\cH(\A)$ and $\SL_2( \A)$ with their
images in $\cH(\A) \rtimes \SL_2( \A),$
and write $\osc(x,y,z)$ or $\osc(g)$ as appropriate.
The group
$\left \{ \left( \begin{smallmatrix} 1&x\\ &1 \end{smallmatrix} \right) \mid x \in \A\right\}$ lifts into $\SL_2( \A)$
and we identify it with its image, writing $\osc \left( \begin{smallmatrix} 1&x\\ &1 \end{smallmatrix} \right)$ as well.
)

It is known that $SL_2( F)$ lifts into $\SL_2( \A),$
and hence $G^J(F)$ lifts into $\wt G^J(\A).$
The  representation $\osc$ has an automorphic realization given by theta functions
$$\tpp(g) = \sum_{\xi \in F} \osc(g) \phi(\xi), \qquad ( g \in  \wt G^J(\A).$$
For each $\phi \in \schwartz(\A),$ the corresponding
theta function is
a genuine function $G^J(F) \bs \wt G^J(\A),$ i.e.,
it does not factor through the projection
to $G^J( \A).$

It is further known that
$\left \{ \left(\begin{smallmatrix}1& x \\ 0&1 \end{smallmatrix} \right): x \in \A\right\}$
lifts into $\SL_2( \A).$
A smooth function $f: SL_2( F) \backslash \SL_2( \A)
\to \C$
or $\jquo \to \C$
has a Fourier expansion
$$
f =\sum_{a \in F}
f^{\psi_a}, \;\;\text{ where } \qquad
f^{\psi_a}(\tilde h):=
\int_{\quo}
f\left( \left(\begin{matrix}1& x \\ 0&1 \end{matrix} \right)
\tilde h
\right) \psi(ax)\, dx.
$$
Since $U(\A)$ and $SL_2( F)$ generate
$\SL_2( \A),$ it follows that
$$f^{\psi_a}=0\; \forall a \in F^\times
\implies
f \text{ is constant on }\SL_2(\A)$$
and thus
\begin{equation}
\label{genuine => not constant}
f \text{ is genuine }
\implies f^{\psi_a}\ne 0
\text{ for some }a \in F^\times.
\end{equation}
\subsubsection{Fourier-Jacobi coefficients}
Define the Fourier-Jacobi coefficient
mapping
$$
F\!J^{\baseaddchar}: C^\infty(\jquo)\times \schwartz(\A) \to C^\infty(SL_2(F) \bs \SL_2(\A)).
$$
$$
F\!J^{\baseaddchar}(f, \phi)(\tilde g_0)
:=\!\!\! \iq{\cH}f(u\tilde g_0 ) \tpp(u\tilde g_0) \, du,
\quad
\left(\tilde g_0 \in \SL_2(\A), \; f\in C^\infty(\jquo), \; \phi \in \schwartz(\A)\right)
$$
Recall that a function $\jquo \to \C$ (resp. $\mquo \to \C$)
is said to be genuine if it does not factor through
the projection to $G^J\quo$ (resp. $SL_2\quo$).
Since $\tpp$ is genuine, it follows that $F\!J^{\baseaddchar}(f, \phi)$ is genuine if and only if $f$ is not,
with the only exception being that $\FJ(f, \phi)$ may
equal zero (which is of course not genuine) when $f$ is
not genuine.

Now, let $$ \Sieg= \left\{
\bpm I&X \\ &I
\epm : X = \bpm y&z \\ r&y \epm
\right\}, \qquad
\psi_{\Sieg,a}:= \baseaddchar( z+ar).
$$
and for $f \in C^\infty( \jquo),$ define
$$f\FC{\Sieg}{a}(\tilde g_0)
:= \iq{\Sieg} f(u\tilde g_0) \psi_{\Sieg}(u) \, du$$

\begin{lem}\label{theta lemma Sp4 version}
Given $f \in C^\infty( \jquo),$ we have
$$
\left(\FJ(\phi, f)\right)^{\psi_a}(\tilde g_0)
= \int_\A
f^{(\Sieg, \psi_{\Sieg, a})}((x,0,0)\tilde g_0)
\left[\osc(\tilde g_0)\phi\right](x)
\, dx
$$
\end{lem}
\begin{proof}
$$\begin{aligned}
\left(\FJ(\phi, f)\right)^{\psi_a}(\tilde g_0)
&= \int_{F\bs \A}
\FJ(\phi, f)\left(\bpm 1&r_1\\&1 \epm
\tilde g_0\right)
\baseaddchar(ar_1)
\, dr_1\\
&= \int_{F\bs \A}
\iq{\cH}
f\left(u\bpm 1&r_1\\&1 \epm
\tilde g_0\right)
\tpp\left(u\bpm 1&r_1\\&1 \epm
\tilde g_0\right)
\, du\,
\baseaddchar(ar_1)
\, dr_1.
\end{aligned}
$$

It is convenient to identify unipotent elements with
their images in $Sp_4.$  In this notation we have
\begin{equation}\label{theta lemma Sp4 eq 1}
\iFA\iFA\iFA\iFA
f\left(\bpm 1&x&y+xr_1&z\\&1&r_1&y\\&&1&-x\\&&&1 \epm
\tilde g_0\right)
\tpp\left(\bpm 1&x&y+xr_1&z\\&1&r_1&y\\&&1&-x\\&&&1 \epm
\tilde g_0\right)
\baseaddchar(ar_1)
\, dx
\, dy\, dz
\, dr_1\end{equation}
Now, for any $g \in \wt G^J(\A),$
$$
\tpp(g) = \sum_{\xi \in F}
\osc(g) \phi(\xi)
= \osc\left(\bpm 1&\xi &&\\&1&&\\&&1-\xi&\\&&&1\epm g\right)\phi(0).
$$
Plug this in, use the invariance of $f$ by $G^J(F),$
and collapse summation in $\xi$ with integration
in $x.$  It follows that \eqref{theta lemma Sp4 eq 1}
  is equal to
  \begin{equation}\label{theta lemma Sp4 eq 2}\notag
\int\limits_\A\iFA\iFA\iFA
f\left(\bpm 1&x&y+xr_1&z\\&1&r_1&y\\&&1&-x\\&&&1 \epm
\tilde g_0\right)
\tpp\left(\bpm 1&x&y+xr_1&z\\&1&r_1&y\\&&1&-x\\&&&1 \epm
\tilde g_0\right)
\baseaddchar(ar_1)
\, dy\, dz
\, dr_1
\, dx
\end{equation}
One has
$$
\bpm 1&x&y+xr_1&z\\&1&r_1&y\\&&1&-x\\&&&1 \epm
=\bpm 1&&y'&z'\\&1&r_1&y'\\&&1&\\&&&1 \epm
\bpm 1&x&&\\&1&&\\&&1&-x\\&&&1 \epm,
y'=y+xr_1, z'= z+xy'
$$
and from \eqref{oscillator action}
$$\osc \bpm 1&&y'&z'\\&1&r_1&y'\\&&1&\\&&&1 \epm
\phi_1(0) = \psi(z') \phi_1(0)
\qquad (\forall y',z',r_1\in \A, \phi_1\in \schwartz(\A)).
$$
The result follows.
\end{proof}
\begin{cor}
For $f \in C^\infty(\jquo),$ we have
$$
f\FC\Sieg a \equiv 0
\iff
\FJ(\phi, f)^{\psi_a}\equiv 0 \qquad \forall \phi\in \schwartz(\A).
$$\label{theta lemma Sp4 first cor}
\end{cor}
\begin{proof}
This follows from lemma \ref{theta lemma Sp4 version},
because a smooth function whose integral
against every Schwartz function is zero is the zero function (and vice versa).
\end{proof}
\begin{cor}\label{theta lemma Sp4 second cor}
For $f \in C^\infty(\jquo),$ not genuine, we have
$$
f\FC\Sieg 0 \not\equiv 0
\implies
f\FC\Sieg a \not\equiv 0, \text{ some }a \in F^\times.
$$
\end{cor}
\begin{proof}
Indeed, for each $\phi \in \schwartz(\A),$ the
function  $\FJ(\phi, f)$ is either $0$
or genuine.  If $f\FC\Sieg 0 \not\equiv 0$
then it follows from
corollary \ref{theta lemma Sp4 first cor} that
$\FJ(\phi, f)$ is nonzero for some $\phi.$
It then follows from \eqref{genuine => not constant}
that $\FJ(\phi, f)^{\psi_a}$ is nonzero
for some $a \in F^\times,$ and corollary \ref{theta lemma Sp4 first cor}
completes the proof.
\end{proof}

\part{Odd case}

\section{Notation and statement}
\label{notation and statement for odd case}

We are now ready to formulate the main result of the paper
in the odd case.

\begin{theorem*}[MAIN THEOREM: ODD CASE]\label{oddt:main}
For $r \in \N,$ take $\tau_1, \dots, \tau_r$ to be irreducible unitary
automorphic cuspidal representations of $GL_{2n_1}(\A), \dots, GL_{2n_r}(\A),$
respectively, and let $\tau=\tau_1\boxplus \dots \boxplus \tau_r.$
 Let $\omega$ denote a Hecke character.
Suppose that
\begin{itemize}
\item $\tau_i$ is $\bar{\omega}$- symplectic
 for each $i,$ and
 \item
 $\tau_i \iso \tau_j \Rightarrow i = j.$
 \end{itemize}
Then there exists an irreducible generic cuspidal
automorphic representation $\sigma$ of
$GSpin_{2n+1}(\A)$ such that
\begin{itemize}
\item
$\sigma$ weakly lifts to $\tau,$  and
\item
the central character  of $\sigma$ is $\omega.$
\end{itemize}
\end{theorem*}
In fact, a refinement of this theorem with 
an explicit description of $\sigma$ is given in theorem 
\ref{oddt:maintheorem}, and proved in section 
\ref{s: main, odd}.

\begin{rmk}
The case $n=1$ is trivial
 because $GSpin_3 = GSp_2 = GL_2,$ so the
inclusion $r$ is simply the identity map.
Clearly, $r$ must be one and
$\sigma = \tau_1.$
   Henceforth, we assume $n\ge 2.$  The careful reader will find places where
this assumption is crucial to the validity of the argument.
\end{rmk}


\subsection{Siegel Parabolic}

We will construct an Eisenstein series on $\ourgroup_{2m}$ induced from
a standard parabolic $P=MU$ such that  $M$ is isomorphic to $GL_m\times GL_1.$
There are two such parabolics.  We
choose the one in which we delete the root $e_{m-1}+e_m$ and the
coroot $e_{m-1}^* + e_m^* - e_0^*$ from the based root datum. {\sl
We shall refer to this parabolic as the ``Siegel.''}\label{odddefOfSiegel}

\begin{rmks}
\begin{itemize}
\item
We can identify the based root datum of the Levi $M$ with that of
$GL_m\times GL_1$ in such a fashion that $e_0$ corresponds to
$GL_1$ and does not appear at all in $GL_m.$
We can then identify $M$ itself with $GL_m\times GL_1$ via
 a particular choice of isomorphism which is compatible with this
 and with the
usual usage of $e_i, e_i^*$ for characters, cocharacters of the standard
torus of $GL_m.$
\item
Having made this identification, a Levi $M'$ which is contained in $M$
will be identified with $GL_1 \times GL_{m_1} \times \dots GL_{m_k},$
(for some $m_1, \dots, m_k \in \N$ that add up to $m$) in the natural
way:  $GL_1$ is identified with the $GL_1$ factor of $M,$ and then
$GL_{m_1} \times \dots GL_{m_k}$ is identified with the subgroup of
$M$ corresponding to block diagonal elements with the specified
block sizes, in the specified order.      \label{oddr:identifications}
\item
The lattice of rational characters of $M$ is spanned by
the maps $(g,\alpha)\mapsto \alpha$ and $(g,\alpha)\mapsto \det g.$
Restriction defines an embedding $X(M)\to X(T),$ which sends these
maps to $e_0$ and $(e_1+ \dots +e_m),$ respectively.  By
abuse of notation, we shall refer to the rational character of $M$
corresponding to $e_0$
as $e_0$ as well.
\item
{$\delta_{P}(g,\alpha)=\det g^{(m-1)}$, with $\delta_{P}$ the modulus function of $P.$}
\end{itemize}
\end{rmks}

\label{section with definition of dagger}
The group $\ourgroup_{2n}$ has an involution $\dagger$ which reverses the last two
simple roots.  The effect is such that
$$\pr( ^\dagger g) = \left( \begin{matrix} I_{n-1}&&&\\ &&1& \\ &1&& \\ &&&I_{n-1} \end{matrix}\right)
\pr(g)  \left( \begin{matrix} I_{n-1}&&&\\ &&1& \\ &1&& \\ &&&I_{n-1} \end{matrix}\right).$$

As is well known, there is a group $Pin_{4n} \supset Spin_{4n}$ such that $\pr$
extends to a two-fold covering $Pin_{4n} \rightarrow O_{4n}.$  The involution $\dagger$
can be realized as conjugation by a preimage of the above permutation matrix.

\subsection{Weyl group of $GSpin_{2m}$; it's action on standard Levis and their representations}
\label{odds:WeylGroup}
\begin{lem} \label{lem: WeylGroupIso}
The Weyl group of $G_m$ is canonically identified with that of $SO_m.$
\end{lem}
\begin{proof}
For this lemma only, let $T$ denote the torus of $SO_m$ and $\tilde T$ that of
$G_m.$  Then the following diagram commutes:
$$\begin{diagram}
\node{Z_{G_m}(\tilde T) } \arrow{e,>}\arrow{s,>} \node{N_{G_m}(\tilde T)} \arrow{s,>}\\
\node{Z_{SO_m}( T) } \arrow{e,>} \node{N_{SO_m}( T).}
\end{diagram}$$
Both horizontal arrows are inclusions and both vertical arrows are $\pr.$
\end{proof}
One easily checks that every element of the Weyl group of $SO_{2n}$
is represented by a permutation matrix.  We denote the permutation
associated to $w$ also by $w.$
The set of permutations $w$ obtained is precisely the set of
permutations $w \in \mathfrak S_{2n}$ satisfying,
  \begin{enumerate}
\item $w( 2n+1-i) = 2n+1-w(i)$ and
\item $\det w = 1$ when $w$ is written as a $2n \times 2n$ permutation matrix.
\end{enumerate}
It is well known that the Weyl group of $SO_{2n}$ (or $G_{2n}$)
is isomorphic to
$\mathfrak S_n \rtimes \{\pm 1 \}^{n-1}.$
To fix a concrete isomorphism, we identify $p \in \mathfrak S_n$ with an $n\times n$ matrix in
the usual way, and then with
$$\begin{pmatrix}p&\\ &_tp^{-1}\end{pmatrix} \in SO_{2n}.$$
We identify $\underline\epsilon=(\epsilon_1,\dots,\epsilon_{n-1}) \in \{\pm 1\}^{n-1}$ with the
permutation $p$ of $\{1,\dots, 2n\}$ such that
$$p(i)=\begin{cases}
i&\text{ if }\epsilon_i=1\\
2n+1-i&\text{ if }\epsilon_i=-1,
\end{cases}$$
 where $\epsilon_n$ is defined to be $\prod_{i=1}^{n-1}\epsilon_i.$
 We then identify $(p,\underline\epsilon)\in \mathfrak S_n \times \{\pm 1\}^{n-1}$
 (direct product of sets) with $p\cdot \underline\epsilon\in W_{SO_{2n}}.$

 With this identification made,
 \begin{equation}\label{odde:so2nweyl}
 (p,\underline\epsilon)\cdot \begin{pmatrix}
 t_1&&&&&\\
&\ddots&&&&\\
&&t_n&&&\\
&&&t_n^{-1}&&\\
&&&&\ddots&\\
&&&&&t_1^{-1}\end{pmatrix}\cdot (p,\underline\epsilon)^{-1}
= \begin{pmatrix}
 t_{p^{-1}(1)}^{\epsilon_{p^{-1}(1)}}&&&&&\\
&\ddots&&&&\\
&&t_{p^{-1}(n)}^{\epsilon_{p^{-1}(n)}}&&&\\
&&&t_{p^{-1}(n)}^{-\epsilon_{p^{-1}(n)}}&&\\
&&&&\ddots&\\
&&&&&t_{p^{-1}(1)}^{-\epsilon_{p^{-1}(1)}}\end{pmatrix}.
\end{equation}
\begin{lem}\label{oddl:WeylAction}
Let $(p,\underline\epsilon)\in\mathfrak S_n\rtimes \{\pm1\}^{n-1}$ be idenified
with an element of $W_{SO_{2m}}=W_{G_{2m}}$ as above.  Then
the action on the character and cocharacter lattices of $G_{2m}$ is given as follows:
\begin{eqnarray*}
(p,\underline\epsilon)\cdot e_i&=&\begin{cases}
e_{p(i)}&i>0, \epsilon_{p(i)}=1,\\
-e_{p(i)}&i>0, \epsilon_{p(i)}=-1,\\
e_0+\sum_{\epsilon_{p(i)}=-1} e_{p(i)}& i=0.
\end{cases}\\
(p,\underline\epsilon)\cdot e_i^*&=&
\begin{cases}e_{p(i)}^*&i>0, \epsilon_{p(i)}=1,\\
e_0^*-e_{p(i)}^*&i>0, \epsilon_{p(i)}=-1,\\
e_0^*&i=0.
\end{cases}
\end{eqnarray*}
\end{lem}
\begin{rmk}
Much of this can be deduced from \eqref{odde:so2nweyl}, keeping in mind that $w\in W_G$
acts on cocharacters by $(w\cdot\varphi)(t)=w\varphi(t)w^{-1}$ and on characters by
$(w\cdot \chi)(t) = \chi(w^{-1}tw).$  However, it is more convenient to give a different proof.
\end{rmk}
\begin{proof}
Let $\alpha_i=e_i-e_{i+1}, i=1$ to $n-1$ and $\alpha_n=e_{n-1}+e_n.$
Let $s_i$ denote the elementary reflection in $W_{G_{2n}}$ corresponding to
$\alpha_i.$
Then it
is easily verified that $s_1,\dots,s_{n-1}$
generate a group isomorphic to $\mathfrak S_n$ which acts on $\{e_1,\dots, e_n\}\in X(T)$
and $\{e_1^*,\dots, e_n^*\}\in X^\vee(T)$ by permuting the indices and acts trivially
on $e_0$ and $e_0^*.$
Also
\begin{eqnarray*}
s_{n}\cdot e_i&=&\begin{cases}
e_i&i\ne n-1,n,0\\
e_0+e_n+e_{n-1}&i=0\\
-e_n&i=n-1\\
-e_{n-1}&i=n
\end{cases}\\
s_{n}\cdot e_i^*&=&\begin{cases}
e_i^*&i\ne n-1,n\\
e_0^*-e^*_n&i=n-1\\
e_0^*-e^*_{n-1}\phantom{+e_n}&i=n.
\end{cases}
\end{eqnarray*}
If $\underline \epsilon \in \{\pm1\}^{n-1}$ is such that $\#\{i:\epsilon_i=-1\}=1$ or $2,$
then $\underline \epsilon$ is conjugate to $s_n$ by an element of the subgroup
isomorphic to $\mathfrak S_n$ generated by $s_1,\dots,s_{n-1}.$  An arbitrary element
of $\{\pm 1\}^{n-1}$ is a product of elements of this form, so one is able to
deduce the assertion for general $(p,\underline\epsilon).$
\end{proof}

Observe that the $\mathfrak S_n$ factor in the semidirect product is precisely the Weyl
group of the Siegel Levi.

In the study of {Jacquet modules of induced representations as well as in the study of}  intertwining operators and Eisenstein series (e.g.,  section \ref{odds:EisensteinSeries}
below),  one encounters a certain subset of the Weyl group associated to a standard
Levi, $M.$  Specifically,
$$W(M):= \left\{ w \in W_{G_{2n}}\left| \begin{array}{l}
w \text{ is of minimal length in }w\cdot W_M\\
wMw^{-1}\text{ is a standard Levi of }G_{2n}
\end{array}\right.\right\}.$$
For our purposes, it is enough to consider the case when $M$ is a subgroup of the
Siegel Levi.  In this case it is isomorphic to $GL_{m_1} \times \dots \times GL_{m_r}
\times GL_1$
for some integers $m_1, \dots, m_r$ which add up to $n,$ and we shall only need
to consider the case when $m_i$ is even for each $i.$
(This, of course, forces $n$ to be even as well.)
\begin{lem}\label{oddl:WeylActionii}
For each $w \in W(M)$ with $M$ as above, there exist a permutation
$p \in \mathfrak S_r$ and and element $\underline \epsilon \in \{\pm1\}^r$ such
that, if $m \in M =(g,\alpha)$ with $\alpha \in GL_1$ and
$$g=\begin{pmatrix} g_1 &&\\&\ddots&\\&&g_r\end{pmatrix} \in GL_n,$$
then
$$wmw^{-1} =(g',\alpha \cdot \prod_{\epsilon_i=-1}\det g_i)
\quad g'=\begin{pmatrix} g'_1&&\\&\ddots&\\&&g'_r\end{pmatrix},$$
where $$g'_i\approx\begin{cases}g_{p^{-1}(i)}&\text{ if }\epsilon_{p^{-1}(i)} = 1,\\
_tg_{p^{-1}(i)}^{-1}& \text{ if }\epsilon_{p^{-1}(i)} = -1.
\end{cases}$$
Here $\approx$ has been used to denote equality up to an inner automorphism.
The map $(p,\underline \epsilon)\mapsto w$ is a bijection between
$W(M)$ and $\mathfrak S_r\times \{\pm 1\}^r.$
(Direct product of sets:  $W(M)$ is not, in general, a group.)
\end{lem}
\begin{proof}
We first prove that $wMw^{-1}$ is again contained in the Siegel Levi.

The Levi $M$ determines an equivalence relation $\sim$ on the set of indices,
$\{1, \dots, n\}$
defined by the condition that $i\sim i+1$ iff $e_i-e_{i+1}$ is
an root of $M.$  View $w$ of $W(M)$  as a permutation
of $\{1, \dots, 2n\}.$ Because $w$ is of minimal length,
$i\sim j, \;i<j\Rightarrow w(i)<w(j).$  Because
$wMw^{-1}$ is a standard Levi, we may deduce that if
$i\sim i+1$   then $w(i+1)=w(i)+1,$ except possibly when $w(i)=n-1,$ in which case
$w(i+1)$ could, {\it a priori} be $n+1.$  However, it is easy to check that in the
special case when all $m_i$ are even, the condition $\det w =1$
forces $w(i+1)=w(i)+1$ even if $w(i)=n-1.$
It follows that  $wMw^{-1}$ is  contained in the Siegel Levi.

When viewed as elements of $\mathfrak S_n \rtimes \{\pm 1\}^{n-1},$ the elements
of $W(M)$ are those pairs $(p,\underline\epsilon)$ such that
$i\sim j\Rightarrow \epsilon_i=\epsilon_j,$
and
$i\sim i+1\Rightarrow p(i+1)=p(i)+\epsilon_i.$ This gives the identification with $\mathfrak S_r\times \{\pm 1\}^r.$

It is clear that the precise value of $g_i'$ is determined only up to conjugacy by an
element of the torus (because we do not specify a particular representative
for our Weyl group element).
By Theorem 16.3.2 of \cite{Springer}, it may be discerned, to this level of
precision, by looking at the effect of $w$ on the based root datum of $M.$
The result now follows from Lemma \ref{oddl:WeylAction}.
\end{proof}
\begin{cor}\label{oddc:weylAction}
Let $w\in W(M)$ be associated to $(p,\underline\epsilon) \in \mathfrak S_r\times \{\pm 1\}^r$
as above.  Let $\tau_1, \dots, \tau_r$ be irreducible
cuspidal representations of $GL_{m_1}(\A),\dots, GL_{m_r}(\A),$
respectively,
and let $\omega$ be a Hecke character.
Then our identification of $M$ with $GL_{m_1} \times \dots\times GL_{m_r}\times GL_1$
determines an identification of $\bigotimes_{i=1}^r \tau_i \boxtimes \omega$ with
a representation of $M(\A).$  Let $M'=wMw^{-1}.$  Then
$M'$ is also identified, via \ref{oddr:identifications} with
$GL_{m_{p^{-1}(1)}} \times \dots\times GL_{m_{p^{-1}(r)}}\times GL_1,$
and we have
$$\bigotimes_{i=1}^r \tau_i \boxtimes \omega\circ Ad(w^{-1})
=\bigotimes_{i=1}^r \tau_i' \boxtimes \omega,$$
where
$$\tau_i'=\begin{cases}\tau_{p^{-1}(i)}&\text{ if }\epsilon_{p^{-1}(i)} = 1,\\
\tilde \tau_{p^{-1}(i)} \otimes \omega & \text{ if }\epsilon_{p^{-1}(i)} = -1.\end{cases}$$
\end{cor}
\begin{proof}
 The
contragredient $\tilde \tau_i$ of $\tau_i$ may be realized as an action on
the same space of functions as $\tau_i$ via $g \cdot \varphi(g_1)
= \varphi ( g_1 \;_tg^{-1}).$
This follows from strong multiplicity one and the analogous statement
for local representations, for which see
 \cite{MR0404534} page 96, or
\cite{MR0425030} page 57.  Combining this fact with the Lemma, we obtain
the Corollary.
\end{proof}

\section{Unramified Correspondence}
\label{odds:combolemma}\label{sec: unramified}
\begin{lem}
\label{oddl:combolemma}
Suppose that $\tau\iso \otimes_v' \tau_v$
is an $\bar\omega$-symplectic irreducible cuspidal
automorphic representation of
$GL_{2n}(\A).$  Let $v$ be a place such that $\tau_v$
is unramified.  Let $t_{\tau,v}$ denote the semisimple conjugacy class in
$GL_{2n}(\C)$ associated to $\tau_v.$
Let $r: GSp_{2n}(\C) \to GL_{2n}(\C)$ be the natural inclusion.
Then $t_{\tau,v}$ contains elements
of the image of $r.$
\end{lem}
\begin{proof}
For convenience in the application, we take $GL_{2n}$ to be identified
with a subgroup of the Levi of the Siegel parabolic as in section \ref{odddefOfSiegel}.
Since
$\tau_v$ is both unramified and generic, it is isomorphic to
$\Ind^{GL_{2n}(F_v)}_{B(GL_{2n})(F_v)} \mu$ for some unramified
character $\mu$ of the maximal torus $T(GL_{2n})(F_v)$ such that this induced representation
is irreducible.
(See \cite{cartiercorvallis}, section 4,  \cite{Zelevinsky}  Theorem 8.1, p. 195.)
Let $\mu_i = \mu \circ e_i^*.$

Since $ \tau \iso\tilde\tau \otimes \omega,$ it follows that $\tau_v \iso\tilde\tau_v \otimes \omega_v$
and from this we deduce that $\{\mu_i: 1 \leq i \leq 2n\}$ and $\{ \mu_i^{-1} \omega  : 1 \leq i
\leq 2n\}$ are the same set.

By Theorem 1, p. 213 of \cite{JS}, we have $\prod_{i=1}^{2n} \mu_i = \omega^n.$

Now, what we need to prove is the following:  if $S$ is a set of $2n$ unramified
characters of
$F_v,$ such that
\begin{enumerate}
\item{$\prod_{i=1}^{2n} \mu_i = \omega^n$}
\item{For each $i$ there exists $j$ such that $\mu_i = \mu_j^{-1} \omega $}
\end{enumerate}
then there is a permutation $\sigma: \{ 1, \dots, 2n \} \to \{1, \dots, 2n\}$ such that
$\mu_{\sigma(i) } = \omega  \mu_{2n - \sigma(i)}^{-1}$ for $i=1$ to $n.$
This we prove by induction on $n.$  When $n=1,$ we know that
$\mu_1 = \mu_i^{-1} \omega $ for $i=1$ or $2.$  If $i=2$ we are done,
while if $i=1$ we use $\omega  = \mu_1 \mu_2$ to obtain $\mu_1=\mu_2,$
and the desired assertion.  Now, if $n>1$ it is sufficient to show that
there exist $i\neq j$ such that $\mu_i =\mu_j^{-1} \omega .$  If there exists
$i$ such that $\mu_i \neq \mu_i^{-1} \omega $ then this is clear.  On the
other hand, there are
exactly
 two unramified characters $\mu$ such that
$\mu = \mu^{-1} \omega.$  The result follows
\end{proof}

{The above argument easily yields:}
\begin{cor}
Suppose $\tau= \tau_1 \boxplus \dots \boxplus \tau_r$
with $\tau_i$ an $\bar\omega$-symplectic irreducible cuspidal
automorphic representation of
$GL_{2n_i}(\A),$ for each $i.$  Then the same conclusion holds.
\end{cor}

\section{Eisenstein series I: Construction and main statements}\label{sec: Eis}
\label{odds:EisensteinSeries}
The main purpose of this section is to construct, for each integer $n\ge 2$ and
Hecke character $\omega,$ a map from the set of all
isobaric representations $\tau$ satisfying the hypotheses of theorem
\ref{oddt:main}
into the residual spectrum of $G_{4n}.$  We use the same notation $\mathcal{E}_{-1}
(\tau, \omega)$ for all $n.$
The construction is given by a multi-residue of an Eisenstein series in
several complex variables, induced from the cuspidal representations
$\tau_1, \dots, \tau_r$ used to form $\tau.$  (Note that by
\cite{JacquetShalika-EPandClassnII}, Theorem 4.4, p.809, this data is
recoverable from $\tau.$)

Let $\omega$ be a Hecke character.
Let $\tau_1, \dots, \tau_r$  be a irreducible cuspidal automorphic representations of
$GL_{2n_1}, \dots, GL_{2n_r},$ respectively.

For each $i,$ let $V_{\tau_i}$ denote the space of cuspforms on which
$\tau_i$ acts.  Then pointwise multipication
$$\varphi_1 \otimes \dots \otimes \varphi_r
\mapsto \prod_{i=1}^r \varphi_i$$
extends to an isomorphism between the abstract tensor
product $\bigotimes_{i=1}^r V_{\tau_i}$ and
the space of all functions
$$\Phi(g_1, \dots, g_r) = \sum_{i=1}^N c_i \prod_{j=1}^r \varphi_{i,j} ( g_j)
\quad c_i \in \C, \; \varphi_{i,j} \in V_{\tau_j}\; \forall i,j.$$
(This is an elementary exercise.)
We consider the representation $\tau_1 \otimes \dots \otimes  \tau_r$
of $GL_{2n_1}\times  \dots\times  GL_{2n_r},$ realized on this latter space,
which we denote $V_{\otimes\tau_i}.$

Let $n=n_1+ \dots +n_r.$

We will construct an Eisenstein series on $\ourgroup_{4n}$ induced from
the subgroup $P=MU$ of the Siegel parabolic such that
$M \iso  GL_{2n_1} \times \dots \times GL_{2n_r}\times GL_1.$
 Let $s_1, \dots s_r$ be a complex variables.  Using the
 identification of $M$ with $GL_{2n_1} \times \dots \times GL_{2n_r}\times GL_1$ fixed in section
 \ref{odddefOfSiegel} above,
we define  an action of $M(\A)$
on the space of $\tau_1 \otimes \dots \otimes  \tau_r$ by
\begin{equation}\label{odde:ProductAction}
( g_1, \dots, g_r, \alpha)\cdot
\prod_{j=1}^r \varphi_j(h_j) = \left( \prod_{j=1}^r \varphi( h_jg_j ) |\det g_j|^{s_j}
\right)
\omega  (\alpha).
\end{equation}
We denote this representation of $M(\A),$ by
$(\bigotimes_{i=1}^r  \tau_i\otimes |\det{}_i|^{s_i} ) \boxtimes \omega .$
(We use $\boxtimes$ to distinguish the ``outer'' tensor product with $\omega$ from the ``inner'' tensor product with  $\det{}_i|^{s_i}.$  Recall that if $V_1, V_2$ are two representations of the same group $G,$ then the ``outer'' tensor product$V_1\boxtimes V_2$ is
the representation of $G\times G$ on the tensor product of the two spaces, while the ``inner'' tensor product $V_1 \otimes V_2$ is the representation of $G$ on the same  space, acting diagonally.)

To shorten the notation, we write $\underline g=(g_1, \dots, g_r).$
Then \eqref{odde:ProductAction} may be shortened to
$$(\underline g, \alpha) \cdot \Phi(\underline h) =\Phi(\underline h \cdot \underline g)
  \left( \prod_{j=1}^r  |\det g_j|^{s_j}\right)
\omega  (\alpha).$$
We shall also employ the shorthand
$\underline s = (s_1, \dots, s_r),$  and $\underline\tau=(\tau_1, \dots, \tau_r).$

For each $\underline s$ we have the induced representation \label{odds:InducedReps}
$\Ind_{P(\A)}^{\ourgroup_{4n}(\A)} (\bigotimes_{i=1}^r  \tau_i\otimes |\det{}_i|^{s_i} )  \boxtimes \omega ,$
(normalized induction)
of $\ourgroup_{4n}(\A).$  The standard realization of this representation
is action by right translation on the space
$V^{(1)}(\underline s,\bigotimes_{i=1}^r \tau_i \boxtimes\omega )$ given by
$$\left\{ \tilde F: \ourgroup_{4n}(\A) \to V_\tau, \text{ smooth } \left|
\tilde F( (\underline g, \alpha) h )(\underline g') = \tilde F(h)(\underline g'\underline g)
\omega (\alpha) |\delta_P|^{\frac12}\prod_{i=1}^r |\det g_i|^{s_i}
\right.\right\}.$$

Where
\begin{equation}\label{deltaP-odd}
|\delta_P|^{\frac12}=\prod_{i=1}^r|\det g_i|^{n-\frac12+\sum_{j=i+1}^r n_i - \sum_{j=1}^{i-1}n_i}
\end{equation}
makes the induction normalized.

A second useful realization is action by right translation on
$$V^{(2)}(\underline s,\bigotimes_{i=1}^r \tau_i \boxtimes\omega)=
\left\{ f:\ourgroup_{4n}(\A) \to \C, \left|
f( h)= \tilde F(h)(e) ,\tilde F \in V^{(1)}(\underline s,\underline \tau,\omega ) \right.\right\}.$$
{Where $e \in GL_{2n}(\A)$ is the identity.}

These {vector spaces}  fit together into a
{vector}  bundle over $\C^r.$  So a section of this bundle is a function $f$ defined
on $\C^r$ such that $f(\underline s) \in
V^{(i)}(\underline s,\bigotimes_{i=1}^r \tau_i \boxtimes\omega)$
($i=1$ or $2$)
for each $\underline s.$
Fix a maximal compact subgroup $K$ of $G_{4n}(\A)$
satisfying the conditions
required in \cite{MW1} (see pages 1 and 4).  Intersecting
$K$ with $M(\A)$ for a standard Levi $M \subset G_{4n},$
we fix maximal compact subgroups of these groups as well.
We shall only require the use of flat, $K$-finite sections, which
are defined as follows.  Take $f_0 \in V^{(i)}(\underline 0,\bigotimes_{i=1}^r \tau_i \boxtimes\omega)$
$K$-finite,
and define $f(\underline s)(h)$ by
$$f(\underline s) (u (\underline g, \alpha) k ) = f_0 (u (\underline g, \alpha) k )\prod_{i=1}^r|\det g_i |^{s_i}
$$for  $u \in U(\A), \underline g \in GL_{2n_1}(\A)
\times \dots \times GL_{2n_r}(\A)
, \alpha \in \A^\times,
k \in K.$
This is well defined.  (I.e., although $g_i$ is not uniquely determined in the decomposition,
$|\det g_i|$ is.  Cf. the definition of $m_P$ on p.7 of \cite{MW1}.)

We begin with a flat $K$ finite section of the bundle of representations realized on the
spaces $V^{(2)}(\underline s,\bigotimes_{i=1}^r \tau_i \boxtimes\omega).$

\begin{rmk}
Clearly, the function $f$ is determined by $f(\underline s^*)$ for any choice
of base point $\underline s^*.$
In particular, any function
of $f$ may be regarded as a function of $f_{\underline s^*} \in
V^{(2)}(\underline s^*,\bigotimes_{i=1}^r \tau_i \boxtimes\omega),$
for any particular value of $\underline s^*.$  We have exploited this fact with $\underline s^*=0$ to streamline the definitions.
{\em A posteriori} it will become clear that the point
$\underline s^*=\poleloc:=(\frac12,\dots,\frac12)$
is of particular importance, and
we shall then switch to $\underline s^* = \poleloc.$
\end{rmk}

For such $f$ the sum
$$E(f)(g)(\underline s) : = \sum_{\gamma \in P(F) \backslash G(F) }
f(\underline s) (\gamma g )$$
converges for all $\underline s$ such that
$\Re(s_r), \Re(s_i-s_{i+1}), i=1$ to $r-1$ are all sufficiently large.
(\cite{MW1}, \S II.1.5, pp.85-86).   It has meromorphic continuation to
$\C^r$ (\cite{MW1} \S IV.1.8(a), IV.1.9(c),p.140).
These are our  Eisenstein series.  We collect some of their well-known
properties in the following theorem.

\begin{thm} \label{thm:Eis}
We have the following:
\begin{enumerate}
\item The function
\begin{equation}\label{odde:bigmultiresidue}
\prod_{i\ne j}(s_i+s_j-1)
\prod_{i=1}^r (s_i-\frac 12)E(f)(g)(\underline s)
\end{equation}
is  holomorphic at $s =\poleloc.$ (More precisely, while $E(f)(g)$ may have singularities,
there is a holomorphic function defined on an open neighborhood of $\underline s=\poleloc$
which agrees with \eqref{odde:multiresidue} on the complement of the  hyperplanes
$s_i=\frac 12,$ and $s_i+s_j=1.$)
\label{oddpoleIsSimple}
\item\label{oddpoleConditions}
The function \eqref{odde:bigmultiresidue} remains holomorphic (in the same sense) when
$s_i+s_j-1$ is omitted, provided $\tau_i \not \cong \omega \otimes \tilde \tau_j.$
It remains holomorphic when $s_i-\frac12$ is omitted, provided $\tau_i$ is not $\bar\omega$-
symplectic.  Furthermore, each of these sufficient conditions is also necessary, in that
the holomorphicity conclusion will fail, for some $f$ and $g,$ if any of the factors is
omitted without the corresponding condition on $\underline\tau$ being satisfied.
From this we deduce that if
\begin{equation}\label{odde:DistinctAndOmegaSymplectic}
\text{ the representations }
\tau_1, \dots, \tau_r\text{ are
all distinct and }\bar\omega\text{-symplectic,}
\end{equation}
then the function
\begin{equation}\label{odde:multiresidue}
\prod_{i=1}^r (s_i-\frac 12)E(f)(g)(\underline s)
\end{equation}
is holomorphic at $s =\poleloc$  for all $f,g$ and nonvanishing at $s=\poleloc$ for
some $f,g.$
\item{
Let us now assume condition \eqref{odde:DistinctAndOmegaSymplectic} holds, and
regard $f$ as a function of  \linebreak $f_{\poleloc} \in
V^{(2)}(\poleloc,\bigotimes_{i=1}^r \tau_i \boxtimes\omega).$
Let $E_{-1}(f_{\poleloc})(g)$ denote the value of the function
\eqref{odde:multiresidue} at $\underline s=\poleloc$ (defined by analytic continuation).
 Then $E_{-1}(f)$ is an $L^2$
function for all  $f_{\poleloc} \in
V^{(2)}(\poleloc,\bigotimes_{i=1}^r \tau_i \boxtimes\omega).$
\label{oddResidueIsL2}
}
\item{
\label{oddResidueMapIsAnIntOp}
The function $E_{-1}$ is an intertwining operator from
$\Ind_{P(\A)}^{\ourgroup_{4n}(\A)} (\bigotimes_{i=1}^r  \tau_i\otimes |\det{}_i|^{\frac12} )
 \boxtimes \omega$ into the space of $L^2$ automorphic forms.}
\item{
\label{oddresidueSupportsPeriod}
If $\residuerep$ is the image of $E_{-1},$ and $\psi_{LW}$ is the character of
$\maxunip$ given by $\psi_{LW}(u) = \baseaddchar( \sum_{i=1}^{2n-1} u_{i,i+1} ),$
then $(\maxunip, \psi_{LW}) \notin \uniper^\perp( \residuerep).$
}
\item
\label{oddindepOfOrder}
The space of functions $\residuerep$ does not depend on the
order chosen on the cuspidal representations $\tau_1, \dots, \tau_r.$  Thus
it is well-defined as a function of the isobaric representation $\tau.$
\end{enumerate}
\label{oddt:EisensteinSeriesProperties}
\end{thm}
\begin{rmk}
By induction in stages, the induced representation $\Ind_{P(\A)}^{\ourgroup_{4n}(\A)}
(\bigotimes_{i=1}^r  \tau_i\otimes |\det{}_i|^{\frac12} )  \boxtimes \omega,$
which comes up in part \eqref{oddResidueMapIsAnIntOp} of the theorem can also
be written as $\Ind_{P_{\text{Sieg}}(\A)}^{G_{4n}(\A)} \tau\otimes|\det|^{\frac12}\boxtimes\omega,$
where $\tau=\tau_1\boxplus \dots\boxplus \tau_r$ as before, and $P_{\text{Sieg}}$ is the Siegel
parabolic.  (Cf. section \ref{s:descrOfTau}.)  Here, we also exploit the identification of the Levi
$M_{\text{Sieg}}$ of $P_{\text{Sieg}}$ with $GL_{2n}\times GL_1$ fixed in
\ref{oddr:identifications}.\label{oddr:InducedFromSiegel}
\end{rmk}
\begin{proof}
The proof is virtually identical to the proof of theorem \ref{t:EisensteinSeriesProperties}.
In two places the proof of theorem \ref{t:EisensteinSeriesProperties}
is slightly more complicated, and therefore we include
complete details for that case, and for this case only describe
the differences.

One must change ``$4n+1$'' to ``$4n,$''
obviously, and one must change ``$\omega^{-1}$-orthogonal'' to ``$\omega^{-1}$-symplectic.''  The twisted exterior square $L$ function plays the role of the twisted symmetric square $L$ function.  The expression for $|\delta_P|^{\frac12}$ is \eqref{deltaP-odd} instead of \eqref{deltaP-even}.  The rational character
$\varepsilon_i, (1 \le i \le r)$ as in \eqref{restns of pos roots} is no
longer the restriction of a positive root, and therefore every restricted root
is indivisible. This simplifies various statments.
Finally, the analogue of remark \ref{blindtoparity} is simpler, since
a representation of $GL_m$ can be $\omega^{-1}$-symplectic only
if $m$ is even (with no condition on $\omega$).
\end{proof}

\section{Descent Construction}\label{s: main, odd}

\subsection{Vanishing of {\it deeper} descents and the descent representation}
In this section, we shall make use of
remark \ref{oddr:InducedFromSiegel}, and regard  $\residuerep$ as affording an automorphic
realization of the representation induced from the representation
$\tau \otimes |\det|^{\frac12}\boxtimes \omega$ of the Siegel Levi.  Thus we may dispense with the smaller
Levi denoted by $P$ in the previous section, and in this section we denote the Siegel
parabolic more briefly by $P=MU.$

Next we describe certain unipotent periods of $\ourgroup_{2m}$
which play a key role in the argument.
For $1 \leq \ell < m,$ let $\descgroup{\ell}$ be the subgroup of $\maxunip$ defined
by $u_{ij}=0$ for $i > \ell.$  (Recall that according to the convention above, this
refers only to those $i,j$ with $i<j\leq m-i.$)
This is the unipotent radical of a standard parabolic $Q_\ell$ having Levi $\desclevi{\ell}$
isomorphic
to $GL_1^\ell \times \ourgroup_{2m -2\ell}.$

Let $\vartheta$ be a character of $\descgroup{\ell}$ then we may define
$$DC^\ell( \tau, \omega , \vartheta)=FC^\vartheta \residuerep.$$

\begin{thm}\label{oddt:DeeperDescentsVanish}
Let $\omega$ be a Hecke character.
Let $\tau=\tau_1\boxplus\dots \boxplus\tau_r$ be an isobaric sum of $\bar\omega$-symplectic irreducible cuspidal automorphic representations $\tau_1, \dots, \tau_r,$
of $GL_{2n_1}(\A),\dots GL_{2n_r}(\A),$ respectively.
If $\ell \geq n,$ and $\vartheta$ is in general position, then
$$DC^\ell( \tau, \omega , \vartheta)=\{0\}.$$
\end{thm}
\begin{proof}
By Theorem \ref{oddt:EisensteinSeriesProperties}, \eqref{oddResidueIsL2}
 the
representation $\residuerep$ decomposes discretely.
Let $\pi\iso \otimes_v'\pi_v$ be one
of the irreducible components, and $p_\pi: \residuerep \to \pi$
the natural projection.

Fix a place $v_0$ such which $\tau_{v_0}$ and $\pi_{v_0}$ are  unramified.
For any
$\xi^{v_0} \in \otimes_{v \neq v_0}'Ind_{P(F_v)}^{G_{4n}(F_v)} \tau_v \otimes
|\det |^{\frac12}_v  \boxtimes \omega _v $ we define a map
$$i_{\xi^{v_0}}:Ind_{P(F_{v_0})}^{G_{4n}(F_{v_0})} \tau_{v_0} \otimes
|\det |^{\frac12}_{v_0} \boxtimes
\omega _{v_0}
\to Ind_{P(\A)}^{G_{4n}(\A)} \tau   \otimes |\det |^{\frac12} \boxtimes \omega $$
by
$i_{\xi^{v_0}}( \xi_v) = \iota( \xi_{v_0}\otimes \xi^{v_0} ),$
where $\iota$ is an isomorphism of the restricted product
$\otimes_v' Ind_{P(F_v)}^{G_{4n}(F_v)} \tau_v  \otimes |\det |^{\frac12}_v \boxtimes \omega _v $
with the global induced representation
$Ind_{P(\A)}^{G_{4n}(\A)} \tau   \otimes |\det |^{\frac12} \boxtimes \omega .$
Clearly
$$\residuerep =
 E_{-1} \circ \iota( \otimes_v' Ind_{P(F_v)}^{G_{4n}(F_v)}
\tau_v  \otimes |\det |^{\frac12}_v \boxtimes \omega _v ).$$
For any decomposable vector $\xi = \xi_{v_0}\otimes \xi^{v_0},$
$$p_\pi \circ E_{-1} \circ \iota(\xi)
=p_\pi \circ   E_{-1} \circ i_{\xi^{v_0}}(\xi_{v_0}).$$

Thus, $\pi_{v_0}$ is a quotient of
$Ind_{P(F_{v_0})}^{G_{4n}(F_{v_0})} \tau_{v_0} \otimes |\det |^{\frac12}_{v_0} \boxtimes
\omega _{v_0},$ and hence (since we took $v_0$ such that $\pi_{v_0}$ is unramified)
it is isomorphic to the unramified constituent
$^{un}Ind_{P(F_{v_0})}^{G_{4n}(F_{v_0})} \tau_{v_0} \otimes |\det |^{\frac12}_{v_0} \boxtimes
\omega _{v_0}.$

Denote the isomorphism of $\pi$ with $\otimes_v'\pi_v$ by the same symbol $\iota.$
This time, fix $\zeta^{v_0} \in \otimes_{v\neq v_0}' \pi_v,$ and define
$i_{\zeta^{v_0}}: ^{un}Ind_{P(F_{v_0})}^{G_{4n}(F_{v_0})} \tau_{v_0} \otimes |\det |^{\frac12}_{v_0} \boxtimes
\omega _{v_0} \to \pi.$  It follows easily from the definitions that
$$FC^\vartheta \circ  i_{\zeta^{v_0}}$$
factors through the Jacquet module
$
\mathcal{J}_{N_\ell, \vartheta}(\;^{un}
Ind_{P(F_{v_0})}^{G_{4n}(F_{v_0})} \tau_{v_0} \otimes |\det |^{\frac12}_{v_0} \boxtimes
\omega _{v_0} ).$
In appendix \ref{oddlocalresults} we show that this Jacquet module is zero.
The result follows.
\end{proof}

\begin{rmk} \label{oddr:characters}
A general character of $\descgroup{\ell}$  is of the form
$$\baseaddchar( c_1 u_{1,2} + \dots + c_{\ell-1} u_{\ell-1,\ell}
+ d_1 u_{\ell , \ell+1} + \dots + d_{4n-2\ell} u_{\ell, 4n-\ell}).$$
The Levi $\desclevi{\ell}$ acts on the space of characters (cf. section \ref{s:uniper}).  Over
an algebraically closed field there is an open orbit, which
consists of all those elements such that
$c_i \neq 0$ for all $i$ and $^t \underline{d} J \underline{d} \neq 0.$
Here,  $\underline{d}$ is the column vector  $^t(d_1, \dots, d_{4n-2\ell}),$
and $J$ is defined as in \ref{ss:GeneralNotation}.
Over a general field two such elements are in the same $F$-orbit
iff the two values of $^t \underline{d} J \underline{d}$ are in the
same square class.
\end{rmk}

Let $\Psi_\ell$ be the character of $\descgroup{\ell}$ defined
by
$$\Psi_\ell( u ) =\baseaddchar (  u_{12}+\dots+u_{\ell-1,\ell}+ u_{\ell, 2n}- u_{\ell, 2n+1} ).$$
It is not hard to see that
\begin{itemize}
\item the stabilizer
$L_\ell^{\Psi_\ell}$ (cf. $M^\vartheta$ in definition \ref{def:fc})
has two connected components,
\item  the one containing the identity
is isomorphic to $\ourgroup_{4n -2 \ell - 1},$
\item there is an ``obvious''
choice of isomorphism $inc:\ourgroup_{4n -2 \ell - 1} \to (L_{\ell}^{\Psi_\ell})^0$
having the following property: if $\{e_i^*: i=0$ to $2n\}$ is the basis for the cocharacter
lattice of $\ourgroup_{4n}$ as in section \ref{rootDataSection}, and
$\{\bar e_i^*, i = 0$ to $2n-\ell-1\}$ is the basis for that of $\ourgroup_{4n-2\ell-1},$ then
\begin{equation}\label{oddinc}
inc \circ \bar e_i^* = \begin{cases}
e_0^*, & i=0\\
e_{\ell+i}^*, &i=1 \text{ to } 2n-\ell-1.
\end{cases}
\end{equation}
\end{itemize}
In the case when $\ell = 2n-1,$
$N_\ell= \maxunip,$ and $\Psi_\ell$ is a generic character.  The above remarks
remain valid with the convention that $G_1=GL_1.$

Let $$DC_\omega ( \tau) = FC^{\Psi_{n-1}} \residuerep.$$
 It
 is a space of smooth functions $G_{2n+1}\quo \to \mathbb C,$ and
 affords a
representation of the group $G_{2n+1}(\A)$
acting by right translation, where we have identified
 $G_{2n+1}$ with  the identity component of $L_{n-1}^{\Psi_{n-1}}.$

\subsection{Main Result}
\begin{thm}\label{oddt:maintheorem}
Let $\omega$ be a Hecke character.
Let $\tau=\tau_1\boxplus\dots \boxplus\tau_r$ be an isobaric sum of distinct $\bar\omega$-symplectic irreducible cuspidal automorphic representations $\tau_1, \dots, \tau_r,$
of $GL_{2n_1}(\A),\dots GL_{2n_r}(\A),$ respectively.
\begin{enumerate}\item
The space $DC_\omega ( \tau)$ is a nonzero cuspidal representation
of $G_{2n+1}(\A).$
Furthermore, the representation $DC_\omega ( \tau)$ supports a nonzero Whittaker integral.
\item
If $\sigma$
 is any irreducible automorphic
representation contained
in $DC_\omega ( \tau),$ then $\sigma$ lifts weakly to $\tau$ under the
map
$r.$  Also, the central character of $\sigma$ is $\omega .$
\end{enumerate}
\end{thm}
\begin{rmk}
Since $DC_\omega(\tau)$ is nonzero and cuspidal, there exists at least one irreducible component
$\sigma.$  In the case of orthogonal groups, one may show (\cite{So-Paris}, pp. 342, item 4) that all of
the components are generic using the Rankin-Selberg integrals of \cite{G-PS-R},\cite{Soudry-SO2n+1Local}.
On the other hand, in the odd case, one may also show
(\cite{GRS4}, Theorem 8, p. 757, or \cite{So-Paris} page 342, item 6) using the results of
\cite{jiangsoudry} that $DC_\omega (\tau)$ is irreducible.
The extension of \cite{jiangsoudry} to $GSpin$ groups
is a work in progress of Takeda and Lau.
\end{rmk}
\subsection{Proof of main theorem}
\begin{proof}
The statements are proved by combining relationships between unipotent
periods and knowledge about $\residuerep.$

\begin{enumerate}
\item{\bf Nonvanishing and genericity}
For genericity, let $(U_1, \psi_1)$ denote the unipotent period obtained by composing the one
which defines the
descent with the one which defines the  Whittaker function
on $\ourgroup_{2n+1}$ embedded into $\ourgroup_{4n}$ as the stabilizer of the  descent
character.  Thus $U_1$ is the subgroup of the standard maximal unipotent defined by
the relations $u_{i,2n}=u_{i,2n+1}$ for $i = n$ to $2n-1,$ and
$$\psi_1( u ) = \psi_0( u_{1,2} + \dots +u_{n-2,n-1} + u_{n-1, 2n} - u_{n-1,2n+1}
+ u_{n, n+1} + \dots + u_{2n-1, 2n}).$$

Next, let $U_2$ denote the subgroup of the standard maximal unipotent defined
by $u_{i,i+1}=0$ for $i$ even and less than $2n.$  (One may also put $\leq 2n$:
the equation $u_{2n,2n+1}=0$ is automatic for any element of $\maxunip.$)
The character $\psi_2$ depends on whether $n$ is odd or even. If $n$ is even, it is
$$\psi_0( u_{1,3} + u_{2,4} + \dots + u_{2n-1, 2n+1} ),$$
while, if $n$ is odd, it is
$$\psi_0( u_{1,3} + u_{2,4} + \dots + u_{2n-3, 2n-1} + u_{2n-2, 2n+1} +u_{2n-1, 2n} ),$$

Finally, let $U_3$ denote the maximal unipotent, and $\psi_3$ denote
$$\psi_3( u) = \psi_0( u_{1,2} + \dots + u_{2n-1,2n}).$$  Thus $(U_3, \psi_3)$ is the
composite of the unipotent period defining the constant term along the Siegel
parabolic, and the one which defines the  Whittaker functional on the
Levi of this parabolic.  By Theorem \ref{oddt:EisensteinSeriesProperties}
\eqref{oddresidueSupportsPeriod} this
period is {\sl not} in $\uniper^\perp( \residuerep).$

In the appendices, we show
\begin{enumerate}
\item{$(U_1, \psi_1) | (U_2, \psi_2),$ in Lemma \ref{oddu1u2lemma}, and }
\item{$ (U_3, \psi_3) \in \langle
(U_2, \psi_2) , \{ (N_\ell, \vartheta): n \leq \ell < 2n \text{ and } \vartheta \text{ in general position.} \}
\rangle$ in Lemma \ref{oddu2u3deep}.
}
\end{enumerate}

By Theorem \ref{oddt:DeeperDescentsVanish} $(N_\ell , \vartheta) \in \uniper^\perp( \residuerep)$
for all $n \leq \ell < 2n$ and $\vartheta$ in general position.  It follows that $(U_1 , \psi_1)
\notin \uniper^\perp( \residuerep).$  This establishes genericity (and hence nontriviality)
of the descent.

\item{\bf Cuspidality}  Turning to cuspidality, we prove in the appendices an identity relating:
 \begin{itemize}
 \item{ Constant terms on $\ourgroup_{2n+1}$ embedded as $(L_{n-1}^{\Psi_{n-1}})^0,$}
 \item{ Descent periods in $\ourgroup_{4n},$}
 \item{ Constant terms on $\ourgroup_{4n},$}
 \item{ Descent periods on $\ourgroup_{4n - 2k },$ embedded in $\ourgroup_{4n}$ as a subgroup of
  a Levi.}
 \end{itemize}
To formulate the exact relationship we introduce some notation for the
maximal parabolics of GSpin groups.

 The group $\ourgroup_{2n+1}$ has one standard maximal parabolic having Levi
 $GL_i \times \ourgroup_{2n-2i+1}$ for each value of $i$ from $1$ to $n.$  Let
 us denote the unipotent radical of this parabolic by $V_i^{2n+1}.$  We denote the
 trivial character of any unipotent group by ${\bf 1}.$

 The group $\ourgroup_{4n}$ has one standard maximal parabolic having Levi
 $GL_k \times \ourgroup_{4n-2k}$ for each value of $k$ from $1$ to $2n-2.$  We denote
 the  unipotent radical of this parabolic by $V_k.$

 (The group $\ourgroup_{4n}$ also has two parabolics with Levi isomorphic to $GL_{2n} \times
 GL_1,$ but since they will not come up in this discussion, we do not need to bother over a
 notation to distinguish them.)

 We prove in  Lemma \ref{oddl:cuspidality-unip-id} that
$(V_k^{2n  + 1}, {\bf 1} ) \circ ( N_{n-1} , \Psi_{n-1})$
is contained in
$$
\langle
(N_{n+k-1} , \Psi_{n+k-1}) , \{
(N_{n+ j-1}, \Psi_{n+j-1})^{(4n-2k+2j)}
\circ ( V_{k-j} , {\bf 1} ): \; \; 1 \leq j < k
\}
\rangle,$$
where $(N_{n+ j-1}, \Psi_{n+j-1})^{(4n-2k+2j)}$ denotes the descent period, defined as above,
but on the group $G_{4n-2k+2j},$ embedded into $G_{4n}$ as a component of the
Levi with unipotent radical $V_{k-j}.$

By Theorem \ref{oddt:DeeperDescentsVanish}
$(N_{n+k-1} , \Psi_{n+k-1}) \in \uniper^\perp( \residuerep)$ for $k=1$ to $n.$
Furthermore, for   $k,j$ such that $1 \leq j < k \leq n,$
the function
$E(f)(s)^{(V_{k-j},{\bf 1})}$ may be expressed
 in terms of Eisenstein series on $GL_{k-j}$ and
$G_{4n-2k+2j}$ using Proposition II.1.7 (ii) of  \cite{MW1}.
What we require is the following:
\begin{lem}
For all $f \in V^{(2)}(\underline s,\bigotimes_{i=1}^r \tau \boxtimes \omega) $
$$\left.E_{-1}(f)^{(V_{k-j},{\bf 1})}\right|_{G_{4n-2k+2j}(\A)}
 \in \bigoplus_S
\mathcal{E}_{-1}( \tau_{S},\omega),$$
where the sum is over subsets $S$ of $\{1, \dots, r\}$ such that
$\sum_{i \in S }2n_i = 2n-k+j,$ and, for each such $S$,
$\mathcal{E}_{-1}( \tau_{S},\omega)$ is the space of functions  on
$G_{4n-2k+2j}(\A)$ obtained by applying the construction
of $\mathcal{E}_{-1}( \tau,\omega)$ to $\{ \tau_i : i \in S\},$
instead of $\{ \tau_i: 1 \le i \le r\}.$
\end{lem}
Once again, this is immediate from  \cite{MW1}
Proposition II.1.7 (ii).

Applying Theorem \ref{oddt:DeeperDescentsVanish}, with $\tau$
replaced by $\tau_S$ and $2n$ by $2n-k+j,$ we deduce
$$(N_{n+ j-1}, \Psi_{n+j-1})^{(4n-2k+2j)} \in \mathcal{U}^\perp\left(
\mathcal{E}_{-1}( \tau_{S},\omega)\right) \quad \forall S,$$
and hence $(N_{n+ j-1}, \Psi_{n+j-1})^{(4n-2k+2j)}
\circ ( V_{k-j} , {\bf 1} ) \in \uniper^\perp( \residuerep).$
This shows that any nonzero function appearing in the space $DC_\omega(\tau)$
must be cuspidal.
Such a function is also easily seen to be of uniformly moderate growth, being the
integral of an automorphic form over a compact domain.
In addition, such a function is easily seen to have central character $\omega,$ and
any function with these properties is necessarily square integrable modulo the
center (\cite{MW1} I.2.12).
It follows that the space $DC_\omega(\tau)$ decomposes discretely.

\item{\bf{The unramified parameters of descent}:}

Now, suppose $\sigma \iso \otimes_v' \sigma_v$ is an irreducible representation
which is {a constituent of } $DC_\omega (\tau).$
{Let $p_\sigma :DC_\omega (\tau) \to \sigma$ be the natural projection.}

Once again, by  Theorem \ref{oddt:EisensteinSeriesProperties}, \eqref{oddResidueIsL2} the
representation $\residuerep$ decomposes discretely.  Let $\pi$ be an irreducible
component of $\residuerep$ such that the restriction of $p_\sigma \circ FC^{\Psi_{n-1}}$ to
$\pi$ is nontrivial.  As discussed
previously in the proof of
Theorem \ref{oddt:DeeperDescentsVanish}, at all but finitely many $v,$
 $\tau$ is unramified at $v$ and
furthermore, $\pi_v$ is the unramified constituent
$^{un}Ind_{P(F_v)}^{G_{4n}(F_v)} \tau_v \boxtimes \omega _v \otimes |\det|_v^{\frac 12}$
of
$Ind_{P(F_v)}^{G_{4n}(F_v)} \tau_v \boxtimes \omega _v \otimes |\det|_v^{\frac 12}.$
If $v_0$ is such a place, the map $p_\sigma \circ FC^{\Psi_{n-1}} \circ i_{\zeta^{v_0}},$ with
$ i_{\zeta^{v_0}}$ defined as in Theorem \ref{oddt:DeeperDescentsVanish},
factors through
$\mathcal{J}_{N_{n-1},\Psi_{n-1}}
\left( \;
^{un}Ind_{P(F_v)}^{G_{4n}(F_v)} \tau_v\otimes |\det|_v^{\frac 12}\boxtimes \omega _v
\right),$
and gives rise
to a $\ourgroup_{2n+1}(F_{v_0})$-equivariant map from this Jacquet-module onto
$\sigma_{v_0}.$

 To pin things down precisely,
assume that $\tau_v$ is the unramified component of $Ind_{B(GL_{2n})(F_v)}^{GL_{2n}(F_v)}
\mu,$ and let $\mu_1 , \dots, \mu_{2n}$ be defined as in the proof of Lemma \ref{oddl:combolemma}.
By Lemma \ref{oddl:combolemma}, we may assume without loss of generality that
$\mu_{2n+1-i}= \omega  \mu_i^{-1}$ for $i = 1$ to $n.$

We also need to refer to the elements of the basis of the cocharacter lattice of
$\ourgroup_{2n+1}$ fixed in section \ref{rootDataSection}.  As in the
remarks preceding the definition of $DC_\omega (\tau),$ we denote
these $\bar e_0^* , \dots , \bar e_n^*.$

In the appendices, we show that
$$\mathcal{J}_{N_{n-1},\Psi_{n-1}}\left(
^{un}Ind_{P(F_v)}^{G_{4n}(F_v)} \tau_v \boxtimes \omega _v \otimes |\det|_v^{\frac 12}
\right)$$
is isomorphic as a $\ourgroup_{2n+1}(F_v)$-module to
$Ind_{B(\ourgroup_{2n+1})(F_v)}^{\ourgroup_{2n+1}(F_v)}\chi$
for $\chi$ the unramified character of $B(\ourgroup_{2n+1})(F_v)$ such that
$$\chi \circ \bar e_i^*= \mu_i, i=1 \text{ to }n, \chi \circ \bar e_0^*=\omega_v .$$

It follows that $\tau$ is a weak lift of $\sigma$ associated to the map $r.$
\end{enumerate}
\end{proof}

\section{Appendix I:  Local results on Jacquet Functors}\label{oddlocalresults}
In this appendix, $F$ is a non-archimedean local field, on which we place
the additional technical hypothesis
\begin{equation}\label{oddiwasawa}
B( G_{2n-1})(F)G_{2n-1}(\mathfrak o ) = G_{2n-1}(F),
\end{equation}
which is known (see \cite{tits}, 3.9, and 3.3.2) to hold at all but
finitely many non-Archimedean completions of a number field.
Here, $G_{2n-1}$ is identified with  $(L_{n-1}^\psi)^0$ is defined as in
\eqref{oddinc},
and $\mathfrak o$ denotes the ring of integers of $F.$

\begin{prop}
Let $\tau = Ind_{B(GL_{2n})(F)}^{GL_{2n}(F)} \mu,$
where $\mu$ satisfies $\mu \circ e_i^* = \omega  \mu \circ e_{ 2n+1-i}^*,$
and let $P$ denote the Siegel parabolic subgroup.
Then for $\ell \ge n$ and $\vartheta$ in general postion, the Jacquet module
$
\mathcal{J}_{N_\ell, \vartheta}(
^{un}Ind_{P(F)}^{G_{4n}(F)} \tau \otimes |\det |^{\frac12} \boxtimes
\omega )$
is trivial.
\end{prop}
\begin{proof}
First, let $\mu_i: F \to \C$ be the unramified character given by
$\mu_i = \mu \circ e_i^*.$  By induction in stages,
$$^{un}Ind_{P(F)}^{G_{4n}(F)} \tau \otimes |\det |^{\frac12} \boxtimes
\omega
=\; ^{un}Ind_{B(G_{4n})(F)}^{G_{4n}(F)} \tilde \mu,$$
where $\tilde \mu \circ e_i^*(x)  = |x|^{\frac12} \mu_i(x),$ for $i =1$ to $2n$ and
$\tilde \mu \circ e_0^* = \omega .$
If $\tilde\mu'$ is the character such that $\tilde \mu' \circ e^*_{2i-1}(x)  = \mu_i(x) |x|^{\frac12},$ and $\tilde \mu'\circ e^*_{2i}(x)
= \mu_i(x) |x|^{-\frac 12},$ for $i =1$ to $n,$  and $\tilde \mu' \circ e^*_0 = \omega,$  then it follows from lemma \ref{oddl:WeylAction} that $\tilde\mu'$ is in the Weyl orbit of $\tilde \mu.$
Hence, by the definition of the unramified constituent
$$^{un}Ind_{B(G_{4n})(F)}^{G_{4n}(F)} \tilde \mu =\; ^{un}Ind_{B(G_{4n})(F)}^{G_{4n}(F)} \tilde \mu'.$$
Now, it is well known that
$$^{un}Ind_{B(GL_2)(F)}^{GL_2(F)} \mu|\;|^{\frac12} \otimes  \mu|\;|^{-\frac12}
= \mu \circ \det.$$
It follows that
$$^{un}Ind_{B(G_{4n})(F)}^{G_{4n}(F)} \tilde \mu'=\;
^{un}Ind_{P_{2^{2n}}(F)}^{G_{4n}(F)} \hat \mu,
$$
where $P_{2^{2n}}$ is the parabolic of $G_{4n}$ having Levi isomorphic to
$GL_{2}^n \times GL_1,$ such that the roots of this Levi are
$e_1-e_2, e_3-e_4, \dots, e_{2n-1}-e_{2n},$ and
$\hat\mu$ is the character given by
$\hat \mu\circ  e_{2i-1}^* = \hat \mu\circ  e_{2i}^* = \mu_i, \hat\mu\circ e_0^* = \omega .$

The remainder of the proof of this lemma as well as
the next proposition may be viewed as a detailed worked example
of theorem 5.2 of \cite{BZ-ASENS}.

The space $Ind_{P_{2^{2n}}(F)}^{G_{4n}(F)} \hat \mu$ has
a filtration as a $\descentparabolic{\ell}(F)$-module, in terms of
$\descentparabolic{\ell}(F)$-modules indexed by the
elements of $(W \cap \middlestageparabolic) \backslash W / (W\cap \descentparabolic{\ell}).$
For any element $x$ of $\middlestageparabolic(F) w \descentparabolic{\ell}(F)$ the module
corresponding to $w$ is isomorphic to $c-ind_{ x^{-1} \middlestageparabolic(F) x
\cap \descentparabolic{\ell}(F) }^{\descentparabolic{\ell}(F)} \middlestagechar \circ Ad( x).$
Here $Ad(x)$ denotes the map given by conjugation by $x.$  It sends
$x^{-1} \middlestageparabolic(F) x
\cap \descentparabolic{\ell}(F)$ into $\middlestageparabolic(F).$
Also, here and throughout $c-ind$ denotes non-normalized compact induction.
(See \cite{cassnotes}, section 6.3.)

Recall from \ref{odds:WeylGroup} that the Weyl group of $G_{4n}$ is identified
(canonically after the choice of $\pr$) with
the set of permutations $w \in \mathfrak S_{4n}$ satisfying,

(1) $w( 4n+1-i) = 4n+1-w(i)$ and

(2)  $\det w = 1$ when $w$ is written as a $4n \times 4n$ permutation matrix.

As representatives for the double
cosets $(W \cap \middlestageparabolic) \backslash W / (W\cap \descentparabolic{\ell})$
we choose the element of minimal
length in each.  As permutations, these elements have the properties

(3) $w^{-1}( 2i ) > w^{-1}(2i-1)$ for $i=1$ to $2n, $ and

(4) If $\ell \leq i < j \leq 4n +1 - \ell$ and $w(i)>w(j),$ then $i = 2n$ and $j= 2n+1.$

Let
$I_ w$
 be the $\descentparabolic{\ell}(F)$-module obtained as $$c-ind_{ \dot w^{-1} \middlestageparabolic(F) \dot w
\cap \descentparabolic{\ell}(F) }^{\descentparabolic{\ell}(F)} \middlestagechar \circ Ad( \dot w)$$
using any element $\dot w$ of $\pr^{-1}( w).$

A function $f$ in $I_w$ will map to zero under the natural projection
to $\mathcal{J}_{N_\ell, \vartheta}(I_w)$ iff there exists a compact subgroup $N^0_\ell$ of $N_\ell(F)$
 such that
$$\int_{N^0_\ell} f( hn) \overline{ \vartheta ( n ) } dn = 0 \qquad \forall h \in \descentparabolic{\ell}(F).$$
(See \cite{cassnotes}, section 3.2.)
Let $\vartheta^h (n)= \vartheta( hn h^{-1}).$  It is easy to see that the integral above
vanishes for suitable $N^0_\ell$ whenever
\begin{equation}\label{oddvarthetah}
\vartheta^h |_{N_\ell(F) \cap w^{-1} P_{2^{2n}}(F) w } \text{ is nontrivial. }
\end{equation}
Furthermore, the function $h \mapsto \vartheta^h$ is continuous in $h,$
(the topology on the space of characters of $N_\ell(F)$ being defined by
identifying it with a finite dimensional $F$-vector space, cf. section \ref{s:uniper})
so if this condition holds for all $h$ in a compact set, then $N_0^\ell$
can be made uniform in $h.$

Now, $\vartheta$ is in general position.  Hence, so is $\vartheta^h$ for
every $h.$  So, if we write
$$\vartheta^h( u ) = \baseaddchar ( c_1 u_{1,2} + \dots + c_{\ell-1} u_{\ell-1,\ell}
+ d_1 u_{\ell , \ell+1} + \dots + d_{2m-2\ell} u_{\ell, 2m-\ell}),$$
we have that $c_i \neq 0$ for all $i$ and $^t \underline{d} J \underline{d} \neq 0.$

Clearly, the condition \eqref{oddvarthetah} holds for all $h$ unless

(5) $w(1) > w(2) > \dots > w( \ell).$

Furthermore, because $^t \underline{d} J \underline{d} \neq 0,$ there exists some
$i_0$ with $\ell+1 \leq i_0 \leq 2n$ such that $d_{i_0-\ell} \neq 0$ and $d_{4n+1+\ell - i_0} \neq 0.$
From this we deduce that the condition \eqref{oddvarthetah} holds for all $h$ unless
$w$ has the additional property

(6) There exists $i_0$ such that $w(\ell) > w(i_0)$ and $w(\ell) > w( 4n+1-i_0).$

However, if $\ell \geq n$ it is easy to check that no permutations with properties
(1),(3) (5) and (6) exist.

Thus $\mathcal{J}_{N_\ell, \vartheta}(I_w)=\{0\}$  for all $w$ and hence
$ \mathcal{J}_{N_\ell, \vartheta}(
^{un}Ind_{P(F)}^{G_{4n}(F)} \tau \otimes |\det |^{\frac12} \boxtimes
\omega )=\{0\}$ by exactness of the Jacquet functor.
\end{proof}

\begin{prop}Let $\tau = Ind_{B(GL_{2n})(F)}^{GL_{2n}(F)} \mu,$
where $\mu$ satisfies $\mu \circ e_i^* = \omega  \mu \circ e_{ 2n+1-i}.$
Then the  Jacquet module
$$\mathcal{J}_{N_{n-1},\Psi_{n-1}}\left(
^{un}Ind_{P(F)}^{G_{4n}(F)}\tau \otimes |\det |^{\frac12} \boxtimes
\omega \right)$$
is isomorphic as a $\ourgroup_{2n+1}(F)$-module to a subquotient of
$Ind_{B(\ourgroup_{2n+1})(F)}^{\ourgroup_{2n+1}(F)}\chi$
for $\chi$ the unramified character of $B(\ourgroup_{2n+1})(F)$ such that
$$\chi \circ \bar e_i^*= \mu_i, i=1 \text{ to }n, \chi \circ \bar e_0^*=\omega .$$
\end{prop}
\begin{proof}
As before, we have
$$^{un}Ind_{P(F)}^{G_{4n}(F)} \tau \otimes |\det |^{\frac12} \boxtimes
\omega =^{un}Ind_{P_{2^{2n}}}^{G_{4n}(F)} \hat \mu,$$
and we filter $Ind_{P_{2^{2n}}(F)}^{G_{4n}(F)} \hat \mu$ in terms of $Q_{n-1}(F)$-modules $I_w.$
This time, $\mathcal{J}_{N_{n-1}, \Psi_{n-1}}(I_w)=\{0\}$ for all $w$
except one.  This one Weyl element, which we denote $w_0,$
corresponds to
the unique permutation satisfying (1),(2),(3),(4)
of the previous result, together with
$w(i) = 4n-2i+1$ for $i=1$ to $n-1.$
Exactness yields
$$\mathcal{J}_{N_{n-1},\Psi_{n-1}}\left(
^{un}Ind_{P(F)}^{G_{4n}(F)}\tau \otimes |\det |^{\frac12} \boxtimes
\omega \right) \iso \mathcal{J}_{N_{n-1}, \Psi_{n-1}}(I_{w_0}).$$
(This is an isomorphism of $Q_{n-1}^{\Psi_{n-1}}(F)$-modules, where $Q_{n-1}^{\Psi_{n-1}}=N_{n-1} \cdot L_{n-1}^{\Psi_{n-1}} \subset Q_{n-1},$
is the stabilizer of $\Psi_{n-1}$ in $Q_{n-1}$ (cf. $L^\vartheta$ above).)

Now, recall that for each $h \in Q_{n-1}(F)$ the character
$\Psi_{n-1}^h ( u) = \Psi_{n-1}( h u h^{-1})$ is a character of $N_{n-1}$ in
general position, and as such determines coefficients
$c_1, \dots , c_{n-2}$ and $d_1, \dots, d_{2n+2}$ as in remark \ref{oddr:characters}.
Clearly,
$$Q_{n-1}^{o}:= \left\{ \left.h \in Q_{n-1}(F)\right|
d_i \neq 0 \text{ for some }i \neq n+1, n+2
\right\}$$
is open.  Moreover, one may see from the description of $w_0$ that
for $h$ in this set \ref{oddvarthetah} is satisfied.
We have an  exact sequence of $Q_{n-1}^{\Psi_{n-1}}(F)$-modules
$$0\to   I_{w_0}^* \to I_{w_0} \to \bar I_{w_0} \to 0,$$
where $I_w^*$ consists of those functions in $I_w$ whose
compact support happens to be contained in $Q_{n-1}^{o},$
and the third arrow is restriction to the complement of $Q_{n-1}^{o}.$
This complement is slightly larger than
$Q_{n-1}^{\Psi_{n-1}}(F)$ in that it contains the full torus
of $Q_{n-1}(F),$ but restriction of functions is an isomorphism
of $Q_{n-1}^{\Psi_{n-1}}(F)$-modules,
$$\bar I_{w_0} \to c-ind_{Q_{n-1}^{\Psi_{n-1}}(F) \cap w_0^{-1} P_{2^{2n}} (F)w_0 }^{Q_{n-1}^{\Psi_{n-1}}(F)}
\hat \mu\delta_{P_{2^{2n}}}^{\frac12} \circ Ad( w_0).$$

Clearly $\mathcal{J}_{N_{n-1},\Psi_{n-1}}\left(I_{w_0}^*\right) =\{0\},$
and hence $$\mathcal{J}_{N_{n-1},\Psi_{n-1}}\left(Ind_{P_{2^{2n}}(F)}^{G_{4n}(F)} \hat \mu
\right)
\iso
\mathcal{J}_{N_{n-1},\Psi_{n-1}}\left(
c-ind_{Q_{n-1}^{\Psi_{n-1}}(F) \cap w_0^{-1} P_{2^{2n}} (F)w_0 }^{Q_{n-1}^{\Psi_{n-1}}(F)}
\hat \mu\delta_{P_{2^{2n}}}^{\frac12} \circ Ad( w_0)
\right).
$$

Now let $\mathcal{W}$ denote
$$\left\{ f: Q_{n-1}^{\Psi_{n-1}}(F) \to \C\left|
\begin{array}{l}
f( uq) = \Psi_{n-1}( u ) f( q)\; \forall\; u \in N_{n-1}(F), \;q \in Q_{n-1}^{\Psi_{n-1}}(F),\\
f( b m ) = \chi( b )\delta_{B(G_{2n+1})}^{\frac12}
 f( m)\; \forall\; b \in B(L_{n-1}^{\Psi_{n-1}})(F), \;m \in L_{n-1}^{\Psi_{n-1}}(F)
 \end{array}\right.
 \right\}.$$
 Set $V=c-ind_{Q_{n-1}^{\Psi_{n-1}} \cap w_0^{-1} P_{2^{2n}} w_0 }^{Q_{n-1}^{\Psi_{n-1}}}
\hat \mu\delta_{P_{2^{2n}}}^{\frac12}  \circ Ad( w_0).$
For $f \in V,$
let
$$W(f)(q) = \int_{N_{n-1}(F) \cap w_0^{-1} \overline{\maxunip}(F) w_0 } f( u q) \bar \Psi_{n-1} ( u )
du.$$
\begin{lem}
The function
$W$ maps $V$ into $\mathcal{W}.$
\end{lem}
\begin{proof}
Note that an element of $V$ is left-invariant by
$N_{n-1} \cap w_0^{-1} U_{\max}(F) w_0,$
and that $\Psi_{n-1} \Big|_{N_{n-1} \cap w_0^{-1} U_{\max}(F) w_0}$
is trivial.
Given this, it easily follows
that for $f \in V,$ $W(f)(uq) = \Psi_{n-1}(u) W(f)(q)$
for all $u \in N_{n-1}(F)$ and $q \in Q_{n-1}(F).$
Further, $w_0$ conjugates $B(G_{2n+1})$ into $P_{2^{2n}},$
and so an element of $V$ is left-$B(G_{2n+1})$-equivariant
with respect to a certain quasicharacter.
The claim is then that the product of this character with
the Jacobian of $Ad(b),\; b \in B(G_{2n+1})(F),$
 acting on $N_{n-1}(F) \cap w_0^{-1} \overline{\maxunip}(F) w_0$
 is
 $\chi\delta_{B(G_{2n+1})}^{\frac12},$
 which is a straightforward calculation.
  \end{proof}

 Let us  denote by $V(N_{n-1}, \Psi_{n-1})$
the kernel of the linear
map $V \to \mathcal{J}_{N_{n-1}, \Psi_{n-1}}(V).$

It is easy to show that $V(N_{n-1}, \Psi_{n-1})$ is contained in the kernel of $W.$
In the next lemma, we show that in fact, they are equal.
Restriction from $Q_{n-1}^{\Psi_{n-1}}(F)$ to $L_{n-1}^{\Psi_{n-1}}(F)$ is clearly an isomorphism
$\mathcal{W} \to
Ind_{B(\ourgroup_{2n+1})(F)}^{\ourgroup_{2n+1}(F)}\chi.$
\end{proof}

\begin{lem}
With notation as in the previous proposition, we have
$Ker(W) \subset V(N_{n-1}, \Psi_{n-1}).$
\end{lem}
\begin{proof}  For this proof, we denote the Borel of $L_{n-1}^{\Psi_{n-1}}$ by $B.$
Also, let $N^{w_0}= N_{n-1} \cap w_0^{-1} P_{2^{2n}} w_0,$
and
$N_{w_0}= N_{n-1} \cap w_0^{-1} \overline{\maxunip} w_0,$

 We consider a smooth function
$f:Q_{n-1}^{\Psi_{n-1}}(F) \to \C$ which is compactly supported modulo
$Q_{n-1}^{\Psi_{n-1}}(F) \cap w_0^{-1} P_{2^{2n}}(F) w_0,$ and satisfies
$$f( bm ) = \chi \delta_B^{\frac12}(b) f(m) \qquad \forall \;b \in B(F),$$ and
$$f(uq ) = f(q) \qquad \forall \; u \in N^{w_0}(F) \text{ and }q \in Q_{n-1}^{\Psi_{n-1}}(F).$$
We assume that
$$\int_{N_{w_0}(F)} f( uq) \bar\Psi_{n-1}( u) du = 0,$$
for all $q\in Q_{n-1}^{\Psi_{n-1}}(F).$
What must be shown is that there is a compact subset $C$ of $N_{n-1}(F)$ such
that
$$\int_C f(gu) \bar\Psi_{n-1}(u) du = 0,$$
for all $q\in Q_{n-1}^{\Psi_{n-1}}(F).$

Consider first $m \in L_{n-1}^{\Psi_{n-1}}( \mathfrak o).$   Let $\mathfrak p$ denote the unique
maximal ideal in $\mathfrak o.$  If $U$ is a unipotent subgroup and $M$ an integer, we define
$$U(\primeideal^M)= \{u \in U(F): u_{ij} \in \primeideal^M \; \forall i \neq j \}.$$
Observe that for each $M\in \N,$ $N_{n-1}(\mathfrak p^M)$ is a subgroup of $N_{n-1}(F)$ which is
preserved by conjugation by  elements of $L_{n-1}^{\Psi_{n-1}}( \mathfrak o).$
 We may choose $M$ sufficiently large that
$supp( f) \subset N^{w_0} (F) N_{w_0}( \mathfrak{p}^{-M} ) L_{n-1}^{\Psi_{n-1}}(F).$  Then
we prove the desired assertion with $C= N_{n-1}( \mathfrak p^{-M}).$
Indeed, for  $m \in L_{n-1}^{\Psi_{n-1}}( \mathfrak o),$  we have
$$\int_{N_{n-1}( \mathfrak p^{-M})} f(mu) \overline\Psi_{n-1}(u) du
=\int_{N_{n-1}( \mathfrak p^{-M})} f(um)  \overline\Psi_{n-1}(u) du,$$
because $Ad(m)$ preserves the subgroup $N_{n-1}(\mathfrak p^{-M}),$
and has Jacobian 1.  Let $c=\Vol( N^{w_0}(\mathfrak p^{-M})),$ which is
finite.  Then by $N^{w_0}$-invariance of $f,$ the above equals
$$ =c \int_{N_{w_0}( \mathfrak p^{-M})} f(um)  \overline\Psi_{n-1}(u) du.$$
This, in turn, is equal to
$$=c \int_{N_{w_0}(F)} f(um)  \overline\Psi_{n-1}(u) du,
$$
since none of the points we have added to the domain of integration
are in the support of $f,$ and this last integral is equal to zero by hypothesis.

Next, suppose $q = u_1m$ with $u_1 \in N_{n-1}(F)$ and $m\in L_{n-1}^{\Psi_{n-1}}( \mathfrak o).$
If $u_1 \in  N_{n-1}(F)-N_{n-1}( \mathfrak p^{-M})$ then
$qu$ is not in the support of $f$ for any $u \in N_{n-1}( \mathfrak p^{-M}).$  On the other hand,
if $u_1 \in N_{n-1}( \mathfrak p^{-M}),$ then
$$\int_{N_{n-1}( \mathfrak p^{-M})} f(u_1mu)  \overline\Psi_{n-1}(u) du
=\int_{N_{n-1}( \mathfrak p^{-M})} f(u_1u m)  \overline\Psi_{n-1}(u) du$$
$$= \Psi_{n-1}( u_1 ) \int_{N_{n-1}( \mathfrak p^{-M})} f(u m)  \overline\Psi_{n-1}(u) du,$$
and now we continue as in the case $u_1=1.$

The result for general $q$ now follows from the left-equivariance properties of $f$
and \eqref{oddiwasawa}.
\end{proof}

\section{Appendix II:  Identities of Unipotent Periods}
\label{s:appII}

\subsection{A  lemma regarding the projection, and a remark}
\label{ss:projlemma}

\begin{lem} The action of $\ourgroup_m$ on itself by conjugation factors through
$\pr.$\end{lem}
\begin{proof}  One has only to check that the kernel of $\pr$ is in the center of
$\ourgroup_m.$  When we regard $\ourgroup_m$ as a quotient of $Spin_m \times GL_1,$
the kernel of $\pr$ is precisely the image of the $GL_1$ factor in the quotient.
\end{proof}

\begin{cor}\label{cor:remark}
 Let $u$ be a unipotent element of $\ourgroup_m(\A)$ and $g$ any element
of $\ourgroup_m(\A).$  Then $\pr( g u g^{-1})$ is a unipotent element of $SO_m(\A)$ and
$g u g^{-1}$ is the unique unipotent element of its preimage in $\ourgroup_m(\A).$
\end{cor}

\begin{rmk}\label{r:anypreimageof}
Recall that the projection $\pr: GSpin_m\to SO_m$ induces an
isomorphism between the unipotent subvarieties of the two groups.
Thus, the unipotent periods of $GSpin_m(\A)$ and $SO_m(\A)$
may be identified.  It follows from corollary \ref{cor:remark} that
any identity or relationship of unipotent periods which is proved
using only conjugation and swapping extends to $GSpin_m(\A).$
The bulk of this appendix may be viewed as a painstaking check
that nearly all the key identities in the descent
construction for special orthogonal groups may be proved
using only conjugation and swapping.
\end{rmk}

\subsection{Relations among Unipotent Periods used in Theorem \ref{oddt:maintheorem}}
Before we proceed with the proofs it will be convenient to formulate the statements
in a slightly different way, making use of the involution $\dagger,$ introduced in section \ref{section with definition of dagger}.

In section \ref{ss:ConjOfUniper},
we introduced the space $\mathcal{U}$
of unipotent periods attached to a reductive group $G(F),$
as well as an action of $G(F)$ on $\mathcal{U}$
by conjugation.
In the special case $G=G_{4n},$
it is convenient to allow ourselves to conjugate our unipotent periods by
elements of the slightly larger group $Pin_{4n}.$
We may  allow the involution $\dagger$ to act on unipotent periods by
$f^{^\dagger(U, \psi_U)}(g) = f^{(U, \psi_U)}(^\dagger g).$
 Denoting the action of $Pin_{4n}(F)$ on $\uniper$ by
$\gamma \cdot( U, \psi_U),$ we have
$$\gamma \cdot( U, \psi_U) \sim \begin{cases}
(U, \psi_U) \; \text{ when }\det \pr \gamma = 1, \\
^\dagger( U, \psi_U) \;  \text{ when }\det \pr \gamma = -1.
\end{cases}$$
Observe that in general $^\dagger(U, \psi_U)$ is {\it not } equivalent to $(U,\psi_U).$
For example, it is not difficult to verify that $^\dagger(\maxunip,\psi_{LW}) \in \uniper^\perp(\residuerep).$

We shall let $(U_1, \psi_1)$ and $(U_3,\psi_3)$ be defined as in
the proof of \ref{oddt:maintheorem}.  We also keep the definition of the group
$U_2.$  However, we now define the character $\psi_2$ by the formula
$$\psi_2( u) = \baseaddchar( u_{13} + \dots + u_{2n-1,2n+1} ),$$ regardless of the parity of
$n.$  (This agrees with the previous definition if $n$ is even; if $n$ is odd they differ
by an application of $\dagger.$)

\begin{lem}\label{oddu1u2lemma}
Let $(U_1, \psi_1)$ be defined as in Theorem \ref{oddt:maintheorem}, and $(U_2, \psi_2)$
defined as just above.
Then $(U_1, \psi_1)|(U_2, \psi_2)$ and $(U_1, \psi_1)|\;^\dagger(U_2, \psi_2).$
\end{lem}
\begin{proof}
We define some additional unipotent periods which appear at intermediate
stages in the argument.
Let $U_4$ be the subgroup defined by $u_{n-1, j}=0$ for $j = n$ to $2n-2$ and
$u_{2n-1,2n}=u_{2n-1,2n+1}.$
We define  a character $\psi_4$ of $U_4$ by the same formula as
$\psi_1.$  Then $(U_1,\psi_1)$ may be swapped for $(U_4, \psi_4).$
(See definition \ref{d:swap}.)

Now, for each $k$ from $1$ to $n,$ define $(U_5^{(k)},\psi_5^{(k)})$ as follows.
 First, for each $k,$ the group $U_5^{(k)}$ is contained in the subgroup of $\maxunip$
defined by, $u_{2n-1, 2n }= u_{2n-1,2n+1}.$
In addition,
$u_{n+k-2,j}= 0$ for $j < 2n-1,$ and $u_{i,i+1}=0$ if $n-k \leq i < n+k$ and $i \equiv n-k \mod 2,$
and
$\psi_5^{(k)}(u)$ equals
$$\baseaddchar\left(
\sum_{i=1}^{n-k-1} u_{i,i+1} + \sum_{i=n-k} ^{n+k-3} u_{i,i+2}
+ u_{n+k-2, 2n}+ u_{n+k-2, 2n+1} + \sum_{i=n+k-1}^{2n-1} u_{i,i+1}
\right).$$
(Note that one or more of the sums here may be empty.)

  Next, let $U_6^{(k)}$  be the subgroup of $\maxunip$ defined by the conditions
  $u_{2n-1, 2n }= u_{2n-1,2n+1},$
  $u_{n+k-2,j}= 0$ for $j < 2n-1,$ and $u_{i,i+1}=0$ if $n-k \leq i < n+k-2$ and $i \equiv n-k+1 \mod 2.$
    The same formula which defines $\psi_5^{(k)}$
  also defines a character of $U_6^{(k)}.$  We denote this character
  by $\psi_6^{(k)}.$

  We make the following observations:
  \begin{itemize}
  \item{$(U_5^{(1)},\psi_5^{(1)})$ is precisely $(U_4, \psi_4).$}
  \item{For each $k,$ $(U_5^{(k)},\psi_5^{(k)})$ is conjugate to $(U_6^{(k+1)}, \psi_6^{(k+1)}).$
  The conjugation is accomplished by any preimage of the
  permutation matrix which transposes $i$ and $i+1$
  for
  $n-k \leq i < n+k$ and $i \equiv n-k \mod 2.$
  }
  \item{
  $(U_6^{(k)} ,\psi_6^{(k)})$ may be swapped for  $(U_5^{(k)}, \psi_5^{(k)}).$  }
  \end{itemize}

 Thus $(U_4, \psi_4) \sim (U_5^{(n)}, \psi_5^{(n)}).$

 Now, let $\psi_2'$ be the character of $U_2$ which is defined by
$$\psi_2'(u)
 = \baseaddchar(u_{1,3} + \dots + u_{2n-2, 2n} -u_{2n-2, 2n+1} + u_{2n-1, 2n+1}).$$
  Then $U_5^{(n)}$ is the subgroup of $U_2$ defined by $u_{2n-1,2n}= u_{2n-1, 2n+1}$
 and  $\psi_5^{(n)}$ is the restriction of $\psi_2'$ to this group.
 Thus $(U_5^{(n)}, \psi_5^{(n)})|(U_2, \psi_2').$  (It is because of this step that
 $(U_1, \psi_1) \not\sim (U_2, \psi_2).$)

 Finally,
 $(U_2, \psi_2)$ and $(U_2, \psi_2')$ are conjugate by the unipotent element
 which projects to
 $I_{4n}-\sum_{i=2}^n e'_{2i-1,2i-2}$

 To obtain $^\dagger( U_2, \psi_2),$ we use
 $$\psi_2''(u)
 := \psi_0(u_{1,3} + \dots + u_{2n-2, 2n} -u_{2n-2, 2n+1} + u_{2n-1, 2n})$$
 instead of $\psi_2'.$
 \end{proof}

\begin{lem}\label{oddu2u3deep}
Let $(U_3, \psi_3)$ be defined as in Theorem \ref{oddt:maintheorem}, and let $(U_2, \psi_2)$
be defined as in the previous lemma.  Then
$$( U_3, \psi_3) \in \langle \; ^{\dagger^n}( U_2, \psi_2) ,
\{ (N_\ell, \vartheta): n \leq \ell < 2n \text{ and } \vartheta \text{ in general position.} \}\rangle.$$
Here $\dagger^n$ indicates that we apply $\dagger$ a total of
$n$ times, with the effect being $\dagger$ if
$n$ is odd and trivial if $n$ is even.
\end{lem}
\begin{proof}
To prove this assertion we introduce some additional unipotent periods.
For $k= 1$ to $2n-1$ let $U_7^{(k)}$ denote the subgroup of $\maxunip$ defined
by $u_{i,i+1}=0$ for $i>k$ and $i \equiv k+1 \mod 2.$  We use two
characters of this group:
 $$\tilde \psi_7^{(k)}= \baseaddchar
 \left( \sum_{1\leq i \leq k-1} u_{i,i+1} + \sum_{k \leq i \leq 2n-1} u_{i,i+2}\right),$$
$$\psi_7^{(k)}= \baseaddchar
\left( \sum_{1\leq i \leq k} u_{i,i+1} + \sum_{k+1 \leq i \leq 2n-1} u_{i,i+2}\right).$$
Then
$(U_7,\psi_7^{(k)})$ is conjugate to $(U_7, \tilde\psi_7^{(k)})$ by any preimage of the
permutation matrix which transposes $i$ and $i+1$ for $k < i < 4n-k$ and $i \equiv k+1 \mod 2.$
This matrix has determinant $-1$ iff $k$ is odd.

  If $k$ is odd then
$(U_7^{(k)}, \psi_7^{(k)})$ may be swapped for  $(U_7^{(k+1)}, \tilde \psi_7^{(k+1)}),$
while if $k$ is even, it may be swapped for
$ (U_{8}^{(k+1)},
\tilde \psi_{8}^{(k+1)}),$
where $U_{8}^{(k+1)}$ is the subgroup of $U_7^{(k+1)}$ defined by $u_{2n-1,2n}=0,$
and $\tilde \psi_{8}^{(k+1)}$ is the restriction of $\tilde \psi_7^{(k+1)}$ to this group.

Now, for $a \in F^\times$ define a character  $\tilde \psi_7^{(k+1, a)}$ of $U_7^{(k+1)}$ by
$$\tilde \psi_7^{(k+1, a)}
=\baseaddchar( u_{1,2} + \dots + u_{k-1, k} + u_{k, k+2} + \dots + u_{2n-1,2n+1}+a u_{2n-1,2n}).$$
Then a Fourier expansion along $U_{2n-1,2n}$ shows that
$$(U_{8}^{(k+1)}, \tilde \psi_{8}^{(k+1)}) \in
\langle (U_7^{(k+1)}, \tilde \psi_7^{(k+1)}),\{ (U_7^{(k+1)}, \tilde \psi_7^{(k+1,a)}):  a \in F^\times \}
\rangle.$$  Here $U_{ij}=\{u\in \maxunip: u_{k,\ell}=0,\, \forall\, (k,\ell) \ne (i,j)\}.$

In Lemma \ref{oddl:fexpdeeper} below we prove that for $k$ even and $a\in F^\times,$
$$(N_{n+\frac k2 } , \Psi_{n+ \frac k2 , a } ) | (U_7^{(k+1)}, \tilde \psi_7^{(k+1,a)}),$$
where
$$\Psi_{\ell, a}( u) = \baseaddchar( u_{1,2} + \dots + u_{\ell-1, \ell} + a u_{\ell, 2n} +u_{\ell, 2n+1} ).$$

The present lemma then follows from the following observations:
 \begin{itemize}
 \item{ $(U_7^{(1)}, \tilde \psi_7^{(1)} ) = (U_2, \psi_2),$  (with $\psi_2$ defined as at the
 beginning of this section). }
  \item{ $(U_7^{(2n-1)}, \psi_7^{(2n-1)}) = (U_3, \psi_3)$}
 \item{ If one applies $\dagger$ to both sides of a relation among unipotent
 periods, it remains valid.}
\item{ The character $\Psi_{n+ \frac k2 , a }$ of $N_{n+\frac k2 }$ is in
general position. (Cf.  remarks \ref{oddr:characters}) }
\item{The set $\{ (N_\ell, \vartheta): n \leq \ell < 2n \text{ and } \vartheta \text{ in general position.} \}$
is stable under $\dagger.$}
 \item{ The number of times we conjugate by the preimage of
 an element of determinant minus 1 in passing
 from $(U_7^{(k)}, \tilde \psi_7^{(k)} )$ back to
 $(U_7^{(k)},  \psi_7^{(k)} )$ is
 precisely $n.$}
 \end{itemize}
 \end{proof}

\begin{lem}\label{oddl:fexpdeeper}
 Let $(N_{n+\frac k2 } , \Psi_{n+ \frac k2 , a } )$ and $(U_7^{(k+1)}, \tilde \psi_7^{(k+1,a)})$
be defined as in the previous lemma.  Then
  $$(N_{n+\frac k2 } , \Psi_{n+ \frac k2 , a } ) | (U_7^{(k+1)}, \tilde \psi_7^{(k+1,a)}).$$
\end{lem}
\begin{proof}
We regard $a$ as fixed for the duration of this argument, and omit it from the notation.
We need still more unipotent periods.  Specifically,
for each $k, \ell$ define
$U_{9}^{(k,\ell)}$ to be the subgroup of $\maxunip$ defined by
requiring that $u_{ij}=0$ under any of the following conditions:
$$k < i \leq k+ 2\ell, \;i \equiv k+1 \mod 2 \text{ and } j=i+1$$
$$ i > k+ 2\ell$$
$$i =k+2\ell - 1, \text{ and } j \neq 4n +1 - k - 2 \ell,$$
$$i = k + 2\ell\text{ and }j <2n.$$
The formula
$$ \baseaddchar( u_{1,2} + \dots + u_{k-1,k}
+u_{k,k+2} + u_{k+1,k+3}
+ \dots + u_{k+2\ell -2, k_{k+2\ell}} + a u_{k+2\ell, 2n}+ u_{k+2\ell, 2n+1} )$$
defines a character of this group which we denote $\psi_{9}^{(k,\ell)}( u ).$
Also, let $U_{10}^{(k,\ell)}$ denote the subgroup of $\maxunip$
defined by
requiring that $u_{ij}=0$ under any of the following conditions:
$$k < i \leq k+ 2\ell, \;i \equiv k+1 \mod 2 \text{ and } j=i+1$$
$$ i > k+ 2\ell  - 1$$
$$i = k + 2\ell-1\text{ and }j > 2n, 2n+1.$$
The formula
$$ \baseaddchar( u_{1,2} + \dots + u_{k,k+1} + u_{k+1,k+3}
+ \dots + u_{k+2\ell -2, k_{k+2\ell}} + a u_{k+2\ell-1, 2n}+ u_{k+2\ell-1, 2n+1} )$$
defines a character of this group
which we denote $\psi_{10}^{(k,\ell)}( u ).$
The period
 $(U_{9},\psi_{9}^{(k,\ell)})$ is conjugate to
$(U_{10}, \psi_{10}^{(k,\ell)}).$

Let $U_{11}^{(k,\ell)}$ denote the subgroup of $\maxunip$
defined by
requiring that $u_{ij}=0$ under any of the following conditions:
$$k < i \leq k+ 2\ell, \;i \equiv k \mod 2 \text{ and } j=i+1$$
$$ i > k+ 2\ell  - 1$$
$$i = k + 2\ell-1\text{ and }j > 2n, 2n+1.$$

Then
$(U_{10}, \psi_{10}^{(k,\ell)})$ may be swapped for $(U_{11}, \psi_{11}^{(k,\ell)}),$
where $\psi_{11}^{(k,\ell)}$ is defined by the same formula as $\psi_{10}^{(k,\ell)}.$

Also, $(U_{11}, \psi_{11}^{(k,\ell)}),$
is clearly divisible by $(U_{9}, \psi_{9}^{(k+1,\ell-1)})$:  to pass from
the former to the latter one simply drops the integration over
$u_{k+2\ell -2,j},$ for $j \neq 4n-k-2\ell +2.$

To complete the argument:  for $k$ even the period
$(U_{9}^{(k+1, n-\frac k2-1)}, \psi_{9}^{(k+1, n-\frac k2-1)})$ divides
the period
$(U_{7}^{(k+1)}, \tilde \psi_7^{(k+1,a)}).$  Indeed the only difference
between the two is that in the former, we omit integration
over $u_{2n-2, 2n}.$

It follows that
$(U_{9}^{(k+1, n-\frac k2-1)}, \psi_{9}^{(k+1, n-\frac k2-1)})$ is divisible
by $(U_{10}^{n+\frac k2-1, 1 } , \psi_{10}^{n+\frac k2 -1, 1 }).$
Finally,  every extension of $\psi_{10}^{n+\frac k2 -1, 1 }$ to a character
of $N_{n+\frac k2 }$ is in the same orbit as $\Psi_{n+ \frac k2 , a }.$
(See Remarks \ref{oddr:characters}.)
Hence
$$(U_{10}^{n+\frac k2-1, 1 } , \psi_{10}^{n+\frac k2 -1, 1 })
 \sim (N_{n+\frac k2 } , \Psi_{n+ \frac k2 , a } ).$$
The result follows.
\end{proof}

\begin{lem}
 As in Theorem \ref{oddt:maintheorem}, let
 $V_i$ denote the unipotent radical of the standard parabolic of $G_{4n}$ having Levi
 isomorphic to $GL_i \times G_{4n-2i}$ (for $1 \leq i \leq 2n-2$).
 Let $V_i^{4n-2m-1}$ denote the unipotent radical of the standard maximal
 parabolic of $G_{2n+1}$ (embedded into $G_{4n}$ as $L_{n-1}^{\Psi_{n-1}}$)
 having Levi isomorphic to $GL_i \times G_{2n-2i+1}$ (for $1 \leq i \leq n$).
 Let  $(N_\ell, \Psi_\ell)$ be the period used to define the descent, as usual,
 and let
 $(N_\ell, \Psi_\ell)^{(4n-2k)}$  denote the analogue for $G_{4n-2k},$ embedded
 into $G_{4n}$ inside the Levi of a maximal parabolic.

 Then,
\label{oddl:cuspidality-unip-id}
$(V_k^{2n  + 1}, {\bf 1} ) \circ ( N_{n-1} , \Psi_{n-1})$
is an element of $$ \langle
(N_{n+k-1} , \Psi_{n+k-1}) , \{
(N_{n+ j-1}, \Psi_{n+j-1})^{(4n-2k+2j)}
\circ ( V_{k-j} , {\bf 1} ): \; \; 1 \leq j < k
\}
\rangle.$$
\end{lem}
\begin{proof}
In this proof, we shall not need to refer to any of the unipotent
periods defined previously.  On the other hand we will need to consider several
new unipotent periods.  For convenience, we start the numbering over from one.

Thus, let $(U_1,\psi_1)= (V_k^{2n  + 1}, {\bf 1} ) \circ ( N_{n-1} , \Psi_{n-1}).$
To describe this group and character in detail, $U_1$ is
 the subgroup defined by
$u_{ij} = 0 $ if $n-1 < i \leq n-1+k < j,$ or $ n-1+k < i$ and $u_{i,2n}=u_{i,2n+1}$ if
$n-1< i \leq n-1+k,$
and $\psi_1$ is given by
$$\psi_1(u)= \baseaddchar( u_{1,2} + \dots + u_{n-2,n-1} + u_{n-1, 2n}-u_{n-1,2n+1} ).$$

 Next, let $U_2$ denote the subgroup of $U_1$
 defined by the additional conditions $u_{ij}=0$ for $1 \leq i \leq n-1 < j \leq n-1+k.$
   Let $\psi_2$ denote the  restriction of $\psi_1$ to this subgroup.

 Next, let $U_3$ denote the subgroup defined by
 $u_{ij} = 0$ for $i \leq k , j \leq n-1+k,$ and $i > n-1+k,$
and $u_{i,2n}=u_{i,2n+1}$
 for $i \leq k.$  Let
 $$\psi_3(u) = \psi_0( u_{k+1, k+2} + \dots + u_{k+n-2,k+n-1} + u_{k+n-1, 2n}- u_{k+n-1, 2n+1}).$$
Then $(U_2, \psi_2)$ is conjugate to $(U_3, \psi_3),$ by any element of
$\ourgroup_{4n}(F)$ which projects to
$$\left(\begin{matrix}
&I_k &&&\\ I_{n-1} &&&& \\ && I_{4n-2m - 2k }  && \\ &&&& I_{n-1} \\ &&& I_k & \end{matrix} \right)
$$
(cf. subsection \ref{ss:projlemma}).

 Finally, let $U_4\supset U_3$ denote the subgroup of $\maxunip$
 given by
 $u_{ij}=0$ if
 $j \leq k+1,$ or $i \geq n+k.$
 Then take  $\psi_4$ defined by the same formula as $\psi_3$

  Certainly $(U_2, \psi_2) | (U_1, \psi_1),$  and $(U_2, \psi_2) \sim (U_3, \psi_3).$
 In Lemma \ref{oddl:u3u4} we prove  that $(U_3, \psi_3)\sim(U_4, \psi_4).$
It follows that $(U_4, \psi_4)| (U_1, \psi_1).$  In fact, one may prove by an
argument similar to the proof of Lemma \ref{oddl:u3u4} that in fact
$(U_2, \psi_2) \sim  (U_1, \psi_1)$ and hence  $(U_4, \psi_4)\sim (U_1, \psi_1).$
But this is not needed for our purposes.

Next, let $U^{(r)}$ denote the subgroup of
$\maxunip$ defined by
 $u_{ij} = 0$ for $j \leq r,$ or $i\geq n+k.$
 So,  $U_4 = U^{(k+1)},$ and $N_{n+k-1}= U^{(1)}.$

Let $\psi^{(r)}$ denote the character of $U^{(r)}$ defined by
$$\psi^{(r)}(u) = \baseaddchar\left(
\sum_{i=r}^{n-2+k} u_{i,i+1} + u_{n-1+k,2n}+u_{n-1+k,2n+1}\right).$$
Then $(U_4, \psi_4) = (U^{(k+1)},\psi^{(k+1)}),$
and $(N_{n+k-1}, \Psi_{n+k-1}) = (U^{(1)}, \psi^{(1)}).$
It is an easy consequence of Lemma \ref{uniperlemma} that
$$(U^{(r)},\psi^{(r)})\in\langle
(U^{(r-1)},\psi^{(r-1)}),(N_{n+ k-r}, \Psi_{n+k-r})^{(4n-2r+2)}
\circ (V_{r-1}, {\bf 1})\rangle.$$
The result follows.
\end{proof}
 \begin{lem}
 \label{oddl:u3u4}
 Let $(U_3, \psi_3)$ and $(U_4, \psi_4)$ be defined as in the previous lemma.
 Then
 $(U_4, \psi_4) \sim (U_3, \psi_3).$
 \end{lem}
  \begin{proof}
  It's clear that $(U_3, \psi_3) | (U_4, \psi_4),$ so we only need to
  prove that \\$(U_4, \psi_4) | (U_3, \psi_3).$  The proof involves a family of groups
  defining intermediate stages.  For $\ell$ such that $1 \leq \ell \leq n-1$ we define
  $U_4^{(\ell)}$ to be the subgroup of $U_4$ defined by the
  condition that for $i \leq k$ the coordinate $u_{ij}$ must be zero for $j \leq k+\ell.$
  Thus $U_4 = U_4^{(1)} \supset U_4^{(2)} \supset \dots \supset
  U_4^{(n-1)} \supset U_3.$
  For each of these groups we consider the period defined using the restriction of $\psi_4.$

  We must show that $(U_4^{(n-1)} ,\psi_4) | (U_3, \psi_3)$ and that
  $(U_4^{(i)}, \psi_4) | (U_4^{(i-1)}, \psi_4).$  In each case, all that is involved is
  an invocation of Lemma \ref{uniperlemma}.   For the first application, what must be checked is
  that the the normalizer of $U_4(F)$ in $G(F)$ permutes $\{\psi_4' : \psi_4' |_{U_3} = \psi_3 \}$
  transitively.  Let $y(\underline{r}) = y( r_1 , \dots, r_k)$ denote the unipotent element in
  $\ourgroup_{4n}(F)$
  which projects to $I+ r_1 e_{1,2n}' + \dots +r_k e_{k,2n}'.$  Then every element
  of $U_4^{(m)}$ is uniquely expressible as $u_3 y( \underline{r}),$
  for $u_3 \in U_3$ and $\underline{r} \in \G_a^k.$
    Hence a map $\psi_4'$
  as above is determined by its composition with $y,$ which defines a character of
  $\quo^k ,$ and hence is of the form
  $$(r_1 , \dots, r_k ) \mapsto \baseaddchar( a_1 r_1 + \dots +   a_k r_k)$$
  for some $a_1, \dots, a_k \in F.$  Consider the unipotent element $z(a_1, \dots , a_k)$ of
  $\ourgroup_{4n}$ which projects to $I+a_1 e'_{k+n-1,1}  + \dots + a_k e'_{k+n-1, k}.$
  We claim first that it normalizes $U_4^{(n-1)},$ and second that
  $\psi_4( z(a) y(r) z( a)^{-1} ) =  \baseaddchar( a_1 r_1 + \dots +   a_k r_k).$
  As noted in  \ref{ss:projlemma}
 this may be checked by a matrix multiplication
  in $SO_{4n}.$

  The proof that $(U_4^{(i)}, \psi_4) | (U_4^{(i-1)}, \psi_4)$ is entirely similar, with the
  role of $y(\underline{r})$ played by $y^{(i)}(\underline{r})$ which projects to
  $I+ r_1 e_{1,k+i+1}' + \dots +r_k e_{k,k+i+1}'$ and that of $z( \underline{a})$ played by
  $z^{(i)}(\underline{a})$ which projects to $I+a_1 e'_{k+i, 1} + \dots + a_k e'_{k+i, k}.$
\end{proof}

\part{Even case}
\section{Formulation of the main result in the  even case }\label{s:even}

Starting with this section, the ``even case'' of descent
from $GL_{2n}$ to $GSpin_{2n}$ will be treated.
This material depends on the general matters covered in part 1, but not on the odd case treated in part 2.

Recall the notion of a weak lift which was reviewed in subsection
\ref{ss:weak lift}.
For $\chi$ a nontrivial quadtratic
character, identify the $L$ group of
$GSpin_{2n}^\chi$
with $GO_{2n}(\C),$ and
consider the inclusion
\begin{equation}\label{functorialityeven}
r:\,^LG=\,^L(GSpin_{2n}^\chi)\to GO_{2n}(\C) \hookrightarrow
GL_{2n}(\C) ={}\; ^LGL_{2n}=^LH.
\end{equation}

We are now ready to formulate our main theorem.

\begin{theorem*}[MAIN THEOREM: EVEN CASE]\label{t:main-even}
For $r \in \N,$ take $\tau_1, \dots, \tau_r$ to be irreducible unitary
automorphic cuspidal representations of $GL_{2n_1}(\A), \dots, GL_{2n_r}(\A),$
respectively, and let $\tau=\tau_1\boxplus \dots \boxplus \tau_r$
be the isobaric sum (see section \ref{s:descrOfTau}).
Let $n=n_1+\dots +n_r,$ and assume that $n\ge 2.$
 Let $\omega$ denote a Hecke character, which is not the square of another Hecke character.
Suppose that
\begin{itemize}
\item $\tau_i$ is $\omega^{-1}$-orthogonal
 for each $i,$ and
 \item
 $\tau_i \iso \tau_j \Rightarrow i = j.$
 \end{itemize}
 For each $i,$ let $\chi_i=\omega_{\tau_i}/\omega^{n_i}$ (which is quadratic),
 and let $\chi=\prod_{i=1}^r \chi_i.$
Then there exists an irreducible generic cuspidal
automorphic representation $\sigma$ of
$GSpin_{2n}^\chi(\A)$ such that
\begin{itemize}
\item
$\sigma$ weakly lifts to $\tau,$  and
\item
the central character $\omega_{\sigma}$ of $\sigma$ is $\omega.$
\end{itemize}
\end{theorem*}
In fact, a refinement of this theorem with 
an explicit description of $\sigma$ is given in theorem 
\ref{t:maintheorem}, and proved in section 
\ref{s: main, even}.

\begin{rmk}  As a consistency check, we note that the case $n=1$ of theorem \ref{t:main-even}
follows from earlier work of
Labesse-Langlands \cite{Labesse-Langlands}. See also
\cite{Kazhdan}.

Indeed, when $n=1,$ the function $L(s, \tau, \sym^2\times
\omega^{-1})$ has a pole  iff $\chi$ is nontrivial, because $L(s,
\tau, \wedge^2 \times \omega^{-1})=L(s,\chi)$.
 In this case the representation $\tau$ that we consider is
  a cuspidal automorphic representation
of $GL(2, \A).$  It is known that in this case $\tilde \tau = \tau
\otimes \omega_\tau^{-1}$ (see, e.g., \cite{Bump-GreyBook},
Theorem 3.3.5, p. 305). It follows that our original $L$-function
on $\tau$ is, in this case, equivalent to requiring that $\tau =
\tau \otimes \chi$ for some nontrivial quadratic character $\tau.$
The automorphic representation obtained from the descent
construction in this case is simply
  a character of $\Res^E_F GL_1(\A),$ where $E$ is
the quadratic extension of $F$ corresponding to $\chi.$ Thus, we
have recovered proposition 6.5, p. 771 of
\cite{Labesse-Langlands}. We thank H. Jacquet for explaining this
to us.
\end{rmk}

\section{Notation}\label{s:even notation}

\subsection{Siegel parabolic}

In this case, we will construct an Eisenstein series on $\ourgroup_{2m+1}$ induced from
a standard parabolic $P=MU$ such that  $M$ is isomorphic to $GL_m\times GL_1.$
There is a unique such parabolic.
We shall refer to this parabolic as the ``Siegel.''\label{defOfSiegel}

\begin{rmk}
\begin{itemize}
\item
We can identify the based root datum of the Levi $M$ with that of
$GL_m\times GL_1$ in such a fashion that $e_0$ corresponds to
$GL_1$ and does not appear at all in $GL_m.$
We can then identify $M$ itself with $GL_m\times GL_1$ via
 a particular choice of isomorphism which is compatible with this
 and with the
usual usage of $e_i, e_i^*$ for characters, cocharacters of the standard
torus of $GL_m.$
\item
Having made this identification, a Levi $M'$ which is contained in $M$
will be identified with $GL_1 \times GL_{m_1} \times \dots GL_{m_k},$
(for some $m_1, \dots, m_k \in \N$ that add up to $m$) in the natural
way:  $GL_1$ is identified with the $GL_1$ factor of $M,$ and then
$GL_{m_1} \times \dots GL_{m_k}$ is identified with the subgroup of
$M$ corresponding to block diagonal elements with the specified
block sizes, in the specified order.      \label{r:identifications}
\item
The lattice of rational characters of $M$ is spanned by
the maps $(g,\alpha)\mapsto \alpha$ and $(g,\alpha)\mapsto \det g.$
Restriction defines an embedding $X(M)\to X(T(G_{2m+1})),$  which sends these
maps to $e_0$ and $(e_1+ \dots +e_m),$ respectively.  By
abuse of notation, we shall refer to the rational character of $M$
corresponding to $e_0$
as $e_0$ as well.
\item
The modulus of $P$ is $(g,\alpha)\to \det g^m.$
\end{itemize}
\end{rmk}

\subsection{Weyl group of $GSpin_{2m+1}$; it's action on standard Levis and their representations}
\label{s:WeylGroup}
Recall lemma \ref{lem: WeylGroupIso}, which establishes an
isomorphism between the Weyl groups of $G_m$ and $SO_m.$
One easily checks that every element of the Weyl group of $SO_{2n+1}$
is represented by a
matrix of the form $w= w_0 \det w_0,$ where $w_0$ is a
permutation matrix.  We denote the permutation
associated to $w_0$ also by $w_0.$
The set of permutations $w_0$ obtained is precisely the set of
permutations $w_0 \in \mathfrak S_{2n}$ satisfying,
  $w_0( 2n+2-i) = 2n+2-w_0(i)$
  It is well known that the Weyl group of $SO_{2n+1}$ (or $G_{2n+1}$)
is isomorphic to
$\mathfrak S_n \rtimes \{\pm 1 \}^{n}.$
To fix a concrete isomorphism, we identify $p \in \mathfrak S_n$ with an $n\times n$ matrix in
the usual way, and then with
$$\begin{pmatrix}p&&\\&1&\\& &_tp^{-1}\end{pmatrix} \in SO_{2n}.$$
We identify $\underline\epsilon=(\epsilon_1,\dots,\epsilon_n) \in \{\pm 1\}^n$ with the
permutation $p$ of $\{1,\dots, 2n+1\}$ such that
$$p(i)=\begin{cases}
i&\text{ if }\epsilon_i=1,\\
2n+2-i&\text{ if }\epsilon_i=-1.
\end{cases}$$
 We then identify $(p,\underline\epsilon)\in \mathfrak S_n \times \{\pm 1\}^n$
 (direct product of sets) with $p\cdot \underline\epsilon\in W_{SO_{2n+1}}.$

 With this identification made,
 \begin{equation}\label{e:so2nweyl}
 (p,\underline\epsilon)\cdot \begin{pmatrix}
 t_1&&&&&&\\
&\ddots&&&&&\\
&&t_n&&&&\\
&&&1&&&\\
&&&&t_n^{-1}&&\\
&&&&&\ddots&\\
&&&&&&t_1^{-1}\end{pmatrix}\cdot (p,\underline\epsilon)^{-1}
= \begin{pmatrix}
 t_{p^{-1}(1)}^{\epsilon_{p^{-1}(1)}}&&&&&&\\
&\ddots&&&&&\\
&&t_{p^{-1}(n)}^{\epsilon_{p^{-1}(n)}}&&&&\\
&&&1&&&\\
&&&&t_{p^{-1}(n)}^{-\epsilon_{p^{-1}(n)}}&&\\
&&&&&\ddots&\\
&&&&&&t_{p^{-1}(1)}^{-\epsilon_{p^{-1}(1)}}\end{pmatrix}.
\end{equation}
\begin{lem}\label{l:WeylAction}
Let $(p,\underline\epsilon)\in\mathfrak S_n\rtimes \{\pm1\}^{n-1}$ be idenified
with an element of $W_{SO_{2m}}=W_{G_{2m}}$ as above.  Then
the action on the character and cocharacter lattices of $G_{2m}$ is given as follows:
\begin{eqnarray*}
(p,\underline\epsilon)\cdot e_i&=&\begin{cases}
e_{p(i)}&i>0, \epsilon_{p(i)}=1,\\
-e_{p(i)}&i>0, \epsilon_{p(i)}=-1,\\
e_0+\sum_{\epsilon_{p(i)}=-1} e_{p(i)}& i=0.
\end{cases}\\
(p,\underline\epsilon)\cdot e_i^*&=&
\begin{cases}e_{p(i)}^*&i>0, \epsilon_{p(i)}=1,\\
e_0^*-e_{p(i)}^*&i>0, \epsilon_{p(i)}=-1,\\
e_0^*&i=0.
\end{cases}
\end{eqnarray*}
\end{lem}
\begin{rmk}
Much of this can be deduced from \eqref{e:so2nweyl}, keeping in mind that $w\in W_G$
acts on cocharacters by $(w\cdot\varphi)(t)=w\varphi(t)w^{-1}$ and on characters by
$(w\cdot \chi)(t) = \chi(w^{-1}tw).$  However, it is more convenient to give a different proof.
\end{rmk}
\begin{proof}
Let $\alpha_i=e_i-e_{i+1}, i=1$ to $n-1$ and $\alpha_n=e_n.$
Let $s_i$ denote the elementary reflection in $W_{G_{2n}}$ corresponding to
$\alpha_i.$
Then it
is easily verified that $s_1,\dots,s_{n-1}$
generate a group isomorphic to $\mathfrak S_n$ which acts on $\{e_1,\dots, e_n\}\in X(T)$
and $\{e_1^*,\dots, e_n^*\}\in X^\vee(T)$ by permuting the indices and acts trivially
on $e_0$ and $e_0^*.$
Also
\begin{eqnarray*}
s_{n}\cdot e_i&=&\begin{cases}
e_i&i\ne n,0\\
e_0+e_n&i=0\\
-e_n&i=n\\
\end{cases}\\
s_{n}\cdot e_i^*&=&\begin{cases}
e_i^*&i\ne n\\
e_0^*-e^*_n&i=n\\
\end{cases}
\end{eqnarray*}
If $\underline \epsilon \in \{\pm1\}^{n-1}$ is such that $\#\{i:\epsilon_i=-1\}=1,$
then $\underline \epsilon$ is conjugate to $s_n$ by an element of the subgroup
isomorphic to $\mathfrak S_n$ generated by $s_1,\dots,s_{n-1}.$  An arbitrary element
of $\{\pm 1\}^{n-1}$ is a product of elements of this form, so one is able to
deduce the assertion for general $(p,\underline\epsilon).$
\end{proof}

Observe that the $\mathfrak S_n$ factor in the semidirect product is precisely the Weyl
group of the Siegel Levi.

In the study of intertwining operators and Eisenstein series (e.g.,  section \ref{s:EisensteinSeries}
below),  one encounters a certain subset of the Weyl group associated to a standard
Levi, $M.$  Specifically,
$$W(M):= \left\{ w \in W_{G_{2n+1}}\left| \begin{array}{l}
w \text{ is of minimal length in }w\cdot W_M\\
wMw^{-1}\text{ is a standard Levi of }G_{2n+1}
\end{array}\right.\right\}.$$
For our purposes, it is enough to consider the case when $M$ is a subgroup of the
Siegel Levi.  In this case it is isomorphic to $GL_{m_1} \times \dots \times GL_{m_r}
\times GL_1$
for some integers $m_1, \dots, m_r$ which add up to $n,$ and we shall only need
to consider the case when $m_i$ is even for each $i.$
(This, of course, forces $n$ to be even as well.)
\begin{lem}\label{l:WeylActionii}
For each $w \in W(M)$ with $M$ as above, there exist a permutation
$p \in \mathfrak S_r$ and and element $\underline \epsilon \in \{\pm1\}^r$ such
that, if $m \in M =(g,\alpha)$ with $\alpha \in GL_1$ and
$$g=\begin{pmatrix} g_1 &&\\&\ddots&\\&&g_r\end{pmatrix} \in GL_n,$$
then
$$wmw^{-1} =(g',\alpha \cdot \prod_{\epsilon_i=-1}\det g_i)
\quad g'=\begin{pmatrix} g'_1&&\\&\ddots&\\&&g'_r\end{pmatrix},$$
where $$g'_i\approx\begin{cases}g_{p^{-1}(i)}&\text{ if }\epsilon_{p^{-1}(i)} = 1,\\
_tg_{p^{-1}(i)}^{-1}& \text{ if }\epsilon_{p^{-1}(i)} = -1.
\end{cases}$$
Here $\approx$ has been used to denote equality up to an inner automorphism.
The map $(p,\underline \epsilon)\mapsto w$ is a bijection between
$W(M)$ and $\mathfrak S_r\times \{\pm 1\}^r.$
(Direct product of sets:  $W(M)$ is not, in general, a group.)
\end{lem}
\begin{proof}
Since $wMw^{-1}$ is a standard Levi which does not contain any
short roots, it is again contained in the Siegel Levi.

The Levi $M$ determines an equivalence relation $\sim$ on the set of indices,
$\{1, \dots, n\}$
defined by the condition that $i\sim i+1$ iff $e_i-e_{i+1}$ is
an root of $M.$  When viewed as elements of $\mathfrak S_n \rtimes \{\pm 1\}^{n-1},$ the elements
of $W(M)$ are those pairs $(p,\underline\epsilon)$ such that $i\sim j\Rightarrow \epsilon_i=\epsilon_j,$ and
$i\sim i+1\Rightarrow p(i+1)=p(i)+\epsilon_i.$
This gives the identification with $\mathfrak S_r\times \{\pm 1\}^r.$

It is clear that the precise value of $g_i'$ is determined only up to conjugacy by an
element of the torus (because we do not specify a particular representative
for our Weyl group element).
By Theorem 16.3.2 of \cite{Springer}, it may be discerned, to this level of
precision, by looking at the effect of $w$ on the based root datum of $M.$
The result now follows from Lemma \ref{l:WeylAction}.
\end{proof}
\begin{cor}\label{c:weylAction}
Let $w\in W(M)$ be associated to $(p,\underline\epsilon) \in \mathfrak S_r\times \{\pm 1\}^r$
as above.  Let $\tau_1, \dots, \tau_r$ be irreducible
cuspidal representations of $GL_{m_1}(\A),\dots, GL_{m_r}(\A),$
respectively,
and let $\omega$ be a Hecke character.
Then our identification of $M$ with $GL_{m_1} \times \dots\times GL_{m_r}\times GL_1$
determines an identification of $\bigotimes_{i=1}^r \tau_i \boxtimes \omega$ with
a representation of $M(\A).$  Let $M'=wMw^{-1}.$  Then
$M'$ is also identified, via \ref{r:identifications} with
$GL_{m_{p^{-1}(1)}} \times \dots\times GL_{m_{p^{-1}(r)}}\times GL_1,$
and we have
$$\bigotimes_{i=1}^r \tau_i \boxtimes \omega\circ Ad(w^{-1})
=\bigotimes_{i=1}^r \tau_i' \boxtimes \omega,$$
where
$$\tau_i'=\begin{cases}\tau_{p^{-1}(i)}&\text{ if }\epsilon_{p^{-1}(i)} = 1,\\
\tilde \tau_{p^{-1}(i)} \otimes \omega & \text{ if }\epsilon_{p^{-1}(i)} = -1.\end{cases}$$
\end{cor}
\begin{proof}
 The
contragredient $\tilde \tau_i$ of $\tau_i$ may be realized as an action on
the same space of functions as $\tau_i$ via $g \cdot \varphi(g_1)
= \varphi ( g_1 \;_tg^{-1}).$
This follows from strong multiplicity one and the analogous statement
for local representations, for which see
 \cite{MR0404534} page 96, or
\cite{MR0425030} page 57.  Combining this fact with the Lemma, we obtain
the Corollary.
\end{proof}

\section{Unramified Correspondence}
\label{s:combolemma}
\begin{lem}
\label{l:combolemma}
Suppose that $\tau\iso \otimes_v' \tau_v$
is an $\omega^{-1}$-orthogonal irreducible cuspidal
automorphic representation of
$GL_{2n}(\A).$  Let $v$ be a place such that $\tau_v$
is unramified.  Let $t_{\tau,v}$ denote the semisimple conjugacy class in
$GL_{2n}(\C)$ associated to $\tau_v.$
Let $r: GO_{2n}(\C) \to GL_{2n}(\C)$ be the natural inclusion.
Then $t_{\tau,v}$ contains elements
of the image of $r.$
\end{lem}
\begin{proof}
For convenience in the application, we take $GL_{2n}$ to be identified
with a subgroup of the Levi of the Siegel parabolic as in section \ref{defOfSiegel}.
Since
$\tau_v$ is both unramified and generic, it is isomorphic to
$\Ind^{GL_{2n}(F_v)}_{B(GL_{2n})(F_v)} \mu$ for some unramified
character $\mu$ of the maximal torus $T(GL_{2n})(F_v)$ such that this induced representation
is irreducible.
(See \cite{cartiercorvallis}, section 4,  \cite{Zelevinsky}  Theorem 8.1, p. 195.)
Let $\mu_i = \mu \circ e_i^*.$

Since $ \tau \iso\tilde\tau \otimes \omega,$ it follows that $\tau_v \iso\tilde\tau_v \otimes \omega_v$
and from this we deduce that $\{\mu_i: 1 \leq i \leq 2n\}$ and $\{ \mu_i^{-1} \omega  : 1 \leq i
\leq 2n\}$ are the same set.
Hence $\prod_{i=1}^{2n} \mu_i = \chi\omega^n,$
where $\chi$ is quadratic.

Now, what we need to prove is the following:  if $S$ is a set of $2n$ unramified
characters of
$F_v,$ such that
for each $i$ there exists $j$ such that $\mu_i = \mu_j^{-1} \omega,$
then there is a permutation $\sigma: \{ 1, \dots, 2n \} \to \{1, \dots, 2n\}$ such that
$\mu_{\sigma(i) } = \omega  \mu_{2n - \sigma(i)}^{-1}$ for $i=1$ to $n-1.$
This we prove by induction on $n.$
If $n=1,$ there is nothing to be proved.

  If $n>1$ it is sufficient to show that
there exist $i\neq j$ such that $\mu_i =\mu_j^{-1} \omega .$  If there exists
$i$ such that $\mu_i \neq \mu_i^{-1} \omega $ then this is clear.  On the
other hand, there are
exactly
 two unramified characters $\mu$ such that
$\mu = \mu^{-1} \omega.$

Now, suppose that $\mu_1,\dots, \mu_{2n}$ have been renumbered according to
$\sigma$ as above.  Then $\mu_{n+1}\mu_n=\omega\chi.$  If
$\chi$ is trivial, it follows that $\mu_i=\omega\mu_{2n-i}^{-1}$ for all $i,$
and hence that the conjugacy class $t_{\tau,v}$ contains elements of the
maximal torus of $GSO_{2n}(\C).$

On the other hand,
if $\chi$ is nontrivial, then
$\mu_n \ne \omega\mu_{n+1}^{-1},$ from which it follows that
$\mu_n^2\mu_{n+1}=\omega$ and $\mu_{n+1}=\chi\mu_n.$
It follows that $t_{\tau,v}$ contains a diagonal element which is conjugate, in $GL_{2n}(\C),$
to an element of the connected component of $GO_{2n}(\C)$ which does not
contain the identity.
\end{proof}

\begin{cor}\label{c:UnramCorr-Isobaric}
Suppose $\tau= \tau_1 \boxplus \dots \boxplus \tau_r$
with $\tau_i$ an $\omega^{-1}$-orthogonal irreducible cuspidal
automorphic representation of
$GL_{2n_i}(\A),$ for each $i.$  Then the same conclusion holds.
\end{cor}

\begin{cortopf}\label{ctp:Unram}
Let $\tau$ be as in corollary \ref{c:UnramCorr-Isobaric}, and let $v$ be a place at which
$\tau$ and $\omega$ are unramified.
Let $\eta$ be one of the two unramified characters such that $\eta^2=\omega_v.$  Let $\chi_{un}$
denote  the unique nontrivial unramified quadratic character
of $F_v^\times.$
Then $\tau_v\iso\Ind_{B(GL_{2n})(F_v)}^{GL_{2n}(F_v)} \mu$
(normalized induction), for an unramified character $\mu$
of the torus of $GL_{2n}(F_v)$
which satisfies either
$$\mu\circ e^*_{2n+1-i}=\omega_v\cdot (\mu\circ e^*_i)^{-1} \quad \forall i=1 \text{ to }n,$$
or
$$\mu\circ e^*_{2n+1-i}=\omega_v\cdot (\mu\circ e^*_i)^{-1} \quad \forall i=1 \text{ to }n-1,
\quad \mu\circ e_n^* =\eta, \quad \mu\circ e_{n+1}^* =\chi_{un}\eta.
$$
\end{cortopf}

\section{Eisenstein series}
\label{s:EisensteinSeries}
The main purpose of this section is to construct, for each integer $n\ge 2$ and
Hecke character $\omega,$ a map from the set of all
isobaric representations $\tau$ satisfying the hypotheses of theorem
\ref{t:main-even}
into the residual spectrum of $G_{4n+1}.$  We use the same notation $\mathcal{E}_{-1}
(\tau, \omega)$ for all $n.$
The construction is given by a multi-residue of an Eisenstein series in
several complex variables, induced from the cuspidal representations
$\tau_1, \dots, \tau_r$ used to form $\tau.$  (Note that by
\cite{JacquetShalika-EPandClassnII}, Theorem 4.4, p.809, this data is
recoverable from $\tau.$)

Let $\omega$ be a Hecke character.
Let $\tau_1, \dots, \tau_r$  be a irreducible cuspidal automorphic representations of
$GL_{2n_1}, \dots, GL_{2n_r},$ respectively.

For each $i,$ let $V_{\tau_i}$ denote the space of cuspforms on which
$\tau_i$ acts.  Then pointwise multipication
$$\varphi_1 \otimes \dots \otimes \varphi_r
\mapsto \prod_{i=1}^r \varphi_i$$
extends to an isomorphism between the abstract tensor
product $\bigotimes_{i=1}^r V_{\tau_i}$ and
the space of all functions
$$\Phi(g_1, \dots, g_r) = \sum_{i=1}^N c_i \prod_{j=1}^r \varphi_{i,j} ( g_j)
\quad c_i \in \C, \; \varphi_{i,j} \in V_{\tau_j}\; \forall i,j.$$
(This is an elementary exercise.)
We consider the representation $\tau_1 \otimes \dots \otimes  \tau_r$
of $GL_{2n_1}\times  \dots\times  GL_{2n_r},$ realized on this latter space,
which we denote $V_{\otimes\tau_i}.$

Let $n=n_1+ \dots +n_r.$

We will construct an Eisenstein series on $\ourgroup_{4n+1}$ induced from
the subgroup $P=MU$ of the Siegel parabolic such that
$M \iso  GL_{2n_1} \times \dots \times GL_{2n_r}\times GL_1.$
 Let $s_1, \dots s_r$ be a complex variables.  Using the
 identification of $M$ with $GL_{2n_1} \times \dots \times GL_{2n_r}\times GL_1$ fixed in section
 \ref{defOfSiegel} above,
we define  an action of $M(\A)$
on the space of $\tau_1 \otimes \dots \otimes  \tau_r$ by
\begin{equation}\label{e:ProductAction}
( g_1, \dots, g_r, \alpha)\cdot
\prod_{j=1}^r \varphi_j(h_j) = \left( \prod_{j=1}^r \varphi( h_jg_j ) |\det g_j|^{s_j}
\right)
\omega  (\alpha).
\end{equation}
We denote this representation of $M(\A),$ by
$(\bigotimes_{i=1}^r  \tau_i\otimes |\det{}_i|^{s_i} ) \boxtimes \omega .$

To shorten the notation, we write $\underline g=(g_1, \dots, g_r).$
Then \eqref{e:ProductAction} may be shortened to
$$\underline g \cdot \Phi(\underline h) =\Phi(\underline h \cdot \underline g)
  \left( \prod_{j=1}^r  |\det g_j|^{s_j}\right)
\omega  (\alpha).$$
We shall also employ the shorthand
$\underline s = (s_1, \dots, s_r),$  and $\underline\tau=(\tau_1, \dots, \tau_r).$

For each $\underline s$ we have the induced representation
$\Ind_{P(\A)}^{\ourgroup_{4n+1}(\A)} (\bigotimes_{i=1}^r  \tau_i\otimes |\det{}_i|^{s_i} )  \boxtimes \omega ,$
(normalized induction)
of $\ourgroup_{4n+1}(\A).$  The standard realization of this representation
is action by right translation on the space \label{s:InducedReps}
$V^{(1)}(\underline s,\bigotimes_{i=1}^r \tau_i \boxtimes\omega )$ given by
$$\left\{ \tilde F: \ourgroup_{4n+1}(\A) \to V_\tau, \text{ smooth } \left|
\tilde F( (\underline g, \alpha) h )(\underline g') = \tilde F(h)(\underline g'\underline g)
\omega (\alpha)
\prod_{i=1}^r |\det g_i|^{s_i+n+\sum_{j=i+1}^r n_i - \sum_{j=1}^{i-1}n_i} \right.\right\}.$$
(The factor
\begin{equation}\label{deltaP-even}
\prod_{i=1}^r|\det g_i|^{n+\sum_{j=i+1}^r n_i - \sum_{j=1}^{i-1}n_i}
\end{equation}
is equal to $|\delta_P|^{\frac12},$ and makes the induction normalized.)
A second useful realization is action by right translation on
$$V^{(2)}(\underline s,\bigotimes_{i=1}^r \tau_i \boxtimes\omega)=
\left\{ f:\ourgroup_{4n+1}(\A) \to \C, \left|
f( h)= \tilde F(h)( id) ,\tilde F \in V^{(1)}(\underline s,\underline \tau,\omega ) \right.\right\}.$$
(Here $id$ denotes the identity element of $GL_{2n}(\A).$)

These representations fit together into a
fiber bundle over $\C^r.$  So a section of this bundle is a function $f$ defined
on $\C^r$ such that $f(\underline s) \in
V^{(i)}(\underline s,\bigotimes_{i=1}^r \tau_i \boxtimes\omega)$
($i=1$ or $2$)
for each $\underline s.$  We shall only require the use of flat, $K$-finite sections, which
are defined as follows.  Take $f_0 \in V^{(i)}(\underline 0,\bigotimes_{i=1}^r \tau_i \boxtimes\omega)$
$K$-finite,
and define $f(\underline s)(h)$ by
$$f(\underline s) (u (\underline g, \alpha) k ) = f_0 (u (\underline g, \alpha) k )\prod_{i=1}^r|\det g_i |^{s_i}
$$for  $u \in U(\A), \underline g \in GL_{2n_1}(\A)
\times \dots \times GL_{2n_r}(\A)
, \alpha \in \A^\times,
k \in K.$
This is well defined.  (I.e., although $g_i$ is not uniquely determined in the decomposition,
$|\det g_i|$ is.  Cf. the definition of $m_P$ on p.7 of \cite{MW1}.)

We begin with a flat $K$ finite section of the bundle of representations realized on the
spaces $V^{(2)}(\underline s,\bigotimes_{i=1}^r \tau_i \boxtimes\omega).$

\begin{rmk}
Clearly, the function $f$ is determined by $f(\underline s^*)$ for any choice
of base point $\underline s^*.$
In particular, any function
of $f$ may be regarded as a function of $f_{\underline s^*} \in
V^{(2)}(\underline s^*,\bigotimes_{i=1}^r \tau_i \boxtimes\omega),$
for any particular value of $\underline s^*.$  We have exploited this fact with $\underline s^*=0$ to streamline the definitions.
{\em A posteriori} it will become clear that the point
$\underline s^*=\poleloc:=(\frac12,\dots,\frac12)$
is of particular importance, and
we shall then switch to $\underline s^* = \poleloc.$
\end{rmk}

For such $f$ the sum
$$E(f)(g)(\underline s) : = \sum_{\gamma \in P(F) \backslash G(F) }
f(\underline s) (\gamma g )$$
converges for all $\underline s$ such that
$\Re(s_r), \Re(s_i-s_{i+1}), i=1$ to $r-1$ are all sufficiently large.
(\cite{MW1}, \S II.1.5, pp.85-86).   It has meromorphic continuation to
$\C^r$ (\cite{MW1} \S IV.1.8(a), IV.1.9(c),p.140).
These are our  Eisenstein series.  We collect some of their well-known
properties in the following theorem.
\begin{thm}
\label{t:EisensteinSeriesProperties}
We have the following:
\begin{enumerate}
\item The function
\begin{equation}\label{e:bigmultiresidue}
\prod_{i\ne j}(s_i+s_j-1)
\prod_{i=1}^r (s_i-\frac 12)E(f)(g)(\underline s)
\end{equation}
is  holomorphic at $s =\poleloc.$ (More precisely, while $E(f)(g)$ may have singularities,
there is a holomorphic function defined on an open neighborhood of $\underline s=\poleloc$
which agrees with \eqref{e:multiresidue} on the complement of the  hyperplanes
$s_i=\frac 12,$ and $s_i+s_j=1.$)
\label{poleIsSimple}
\item\label{poleConditions}
The function \eqref{e:bigmultiresidue} remains holomorphic (in the same sense) when
$s_i+s_j-1$ is omitted, provided $\tau_i \not \cong \omega \otimes \tilde \tau_j.$
It remains holomorphic when $s_i-\frac12$ is omitted, provided $\tau_i$ is not $\omega^{-1}$-
orthogonal.  Furthermore, each of these sufficient conditions is also necessary, in that
the holomorphicity conclusion will fail, for some $f$ and $g,$ if any of the factors is
omitted without the corresponding condition on $\underline\tau$ being satisfied.
From this we deduce that if
\begin{equation}\label{e:DistinctAndOmegaOrthogonal}
\text{ the representations }
\tau_1, \dots, \tau_r\text{ are
all distinct and }\omega^{-1}\text{-orthogonal,}
\end{equation}
then the function
\begin{equation}\label{e:multiresidue}
\prod_{i=1}^r (s_i-\frac 12)E(f)(g)(\underline s)
\end{equation}
is holomorphic at $s =\poleloc$  for all $f,g$ and nonvanishing at $s=\poleloc$ for
some $f,g.$
\item{
Let us now assume condition \eqref{e:DistinctAndOmegaOrthogonal} holds, and
regard $f$ as a function of  \linebreak $f_{\poleloc} \in
V^{(2)}(\poleloc,\bigotimes_{i=1}^r \tau_i \boxtimes\omega).$
Let $E_{-1}(f_{\poleloc})(g)$ denote the value of the function
\eqref{e:multiresidue} at $\underline s=\poleloc$ (defined by analytic continuation).
 Then $E_{-1}(f)$ is an $L^2$
function for all  $f_{\poleloc} \in
V^{(2)}(\poleloc,\bigotimes_{i=1}^r \tau_i \boxtimes\omega).$
\label{ResidueIsL2}
}
\item{
\label{ResidueMapIsAnIntOp}
The function $E_{-1}$ is an intertwining operator from
$\Ind_{P(\A)}^{\ourgroup_{4n+1}(\A)} (\bigotimes_{i=1}^r  \tau_i\otimes |\det{}_i|^{\frac12} )
 \boxtimes \omega$ into the space of $L^2$ automorphic forms.}
\item{
\label{residueSupportsPeriod}
If $\residuerep$ is the image of $E_{-1},$ and $\psi_{LW}$ is the character of
$\maxunip$ given by $\psi_{LW}(u) = \baseaddchar( \sum_{i=1}^{2n-1} u_{i,i+1} ),$
then $(\maxunip, \psi_{LW}) \notin \uniper^\perp( \residuerep).$
}
\item
\label{indepOfOrder}
The space of functions $\residuerep$ does not depend on the
order chosen on the cuspidal representations $\tau_1, \dots, \tau_r.$  Thus
it is well-defined as a function of the isobaric representation $\tau.$
\end{enumerate}
\end{thm}
\begin{rmk}
By induction in stages, the induced representation $\Ind_{P(\A)}^{\ourgroup_{4n}(\A)}
(\bigotimes_{i=1}^r  \tau_i\otimes |\det{}_i|^{\frac12} )  \boxtimes \omega,$
which comes up in part \eqref{ResidueMapIsAnIntOp} of the theorem can also
be written as $\Ind_{P_{\text{Sieg}}(\A)}^{G_{4n}(\A)} \tau\otimes|\det|^{\frac12}\boxtimes\omega,$
where $\tau=\tau_1\boxplus \dots\boxplus \tau_r$ as before, and $P_{\text{Sieg}}$ is the Siegel
parabolic.  (Cf. section \ref{s:descrOfTau}.)  Here, we also exploit the identification of the Levi
$M_{\text{Sieg}}$ of $P_{\text{Sieg}}$ with $GL_{2n}\times GL_1$ fixed in
\ref{r:identifications}.\label{r:InducedFromSiegel}
\end{rmk}
A detailed proof is given in appendices.
First, some preparations are made in section \ref{section: Preparations for the proof of Theorem t:EisensteinSeriesProperties}.  In section \ref{section:  proof of Theorem thm:Eis} theorem \ref{t:EisensteinSeriesProperties} is reduced to a number
of lemmas and propositions.  These lemmas and propositions are then proved in section \ref{s:appIII}.

\section{Descent Construction}\label{s: main, even}

\subsection{Vanishing of {\it deeper} descents and the descent representation}
\label{section:  vanishing of deeper descents (even)}
In this section, we shall make use of
remark \ref{r:InducedFromSiegel}, and regard  $\residuerep$ as affording an automorphic
realization of the representation induced from the representation
$\tau \otimes |\det|^{\frac12}\boxtimes \omega$ of the Siegel Levi.
Thus we may dispense with the smaller
Levi denoted by $P$ in the previous section, and in this section we denote the Siegel
parabolic more briefly by $P=MU.$

Next we describe certain unipotent periods of $\ourgroup_{2m}$
which play a key role in the argument.
For $1 \leq \ell < m,$ let $\descgroup{\ell}$ be the subgroup of $\maxunip$ defined
by $u_{ij}=0$ for $i > \ell.$  (Recall that according to the convention above, this
refers only to those $i,j$ with $i<j\leq m-i.$)
This is the unipotent radical of a standard parabolic $Q_\ell$ having Levi $\desclevi{\ell}$
isomorphic
to $GL_1^\ell \times \ourgroup_{2m -2\ell}.$

Let $\vartheta$ be a character of $\descgroup{\ell}$ then we may define
$$DC^\ell( \tau, \omega , \vartheta)=FC^\vartheta \residuerep.$$

\begin{thm}\label{t:DeeperDescentsVanish}
Let $\omega$ be a Hecke character.
Let $\tau=\tau_1\boxplus\dots \boxplus\tau_r$ be an isobaric sum of $\omega^{-1}$-orthogonal
 irreducible cuspidal automorphic representations $\tau_1, \dots, \tau_r,$
of $GL_{2n_1}(\A),\dots GL_{2n_r}(\A),$ respectively.
If $\ell > n,$ and $\vartheta$ is in general position, then
$$DC^\ell( \tau, \omega , \vartheta)=\{0\}.$$
\end{thm}
\begin{proof}
By Theorem \ref{t:EisensteinSeriesProperties}, \eqref{ResidueIsL2}
 the
representation $\residuerep$ decomposes discretely.
Let $\pi\iso \otimes_v'\pi_v$ be one
of the irreducible components, and $p_\pi: \residuerep \to \pi$
the natural projection.

Fix a place $v_0$ such which $\tau_{v_0}$ and $\pi_{v_0}$ are  unramified.
For any
$\xi^{v_0} \in \otimes_{v \neq v_0}'Ind_{P(F_v)}^{G_{4n}(F_v)} \tau_v \otimes
|\det |^{\frac12}_v  \boxtimes \omega _v $ we define a map
$$i_{\xi^{v_0}}:Ind_{P(F_{v_0})}^{G_{4n}(F_{v_0})} \tau_{v_0} \otimes
|\det |^{\frac12}_{v_0} \boxtimes
\omega _{v_0}
\to Ind_{P(\A)}^{G_{4n}(\A)} \tau   \otimes |\det |^{\frac12} \boxtimes \omega $$
by
$i_{\xi^{v_0}}( \xi_v) = \iota( \xi_{v_0}\otimes \xi^{v_0} ),$
where $\iota$ is an isomorphism of the restricted product
$\otimes_v' Ind_{P(F_v)}^{G_{4n}(F_v)} \tau_v  \otimes |\det |^{\frac12}_v \boxtimes \omega _v $
with the global induced representation
$Ind_{P(\A)}^{G_{4n}(\A)} \tau   \otimes |\det |^{\frac12} \boxtimes \omega .$
Clearly
$$\residuerep =
 E_{-1} \circ \iota( \otimes_v' Ind_{P(F_v)}^{G_{4n}(F_v)}
\tau_v  \otimes |\det |^{\frac12}_v \boxtimes \omega _v ).$$
For any decomposable vector $\xi = \xi_{v_0}\otimes \xi^{v_0},$
$$p_\pi \circ E_{-1} \circ \iota(\xi)
=p_\pi \circ   E_{-1} \circ i_{\xi^{v_0}}(\xi_{v_0}).$$

Thus, $\pi_{v_0}$ is a quotient of
$Ind_{P(F_{v_0})}^{G_{4n}(F_{v_0})} \tau_{v_0} \otimes |\det |^{\frac12}_{v_0} \boxtimes
\omega _{v_0},$ and hence (since we took $v_0$ such that $\pi_{v_0}$ is unramified)
it is isomorphic to the unramified constituent
$^{un}Ind_{P(F_{v_0})}^{G_{4n}(F_{v_0})} \tau_{v_0} \otimes |\det |^{\frac12}_{v_0} \boxtimes
\omega _{v_0}.$

Denote the isomorphism of $\pi$ with $\otimes_v'\pi_v$ by the same symbol $\iota.$
This time, fix $\zeta^{v_0} \in \otimes_{v\neq v_0}' \pi_v,$ and define
$i_{\zeta^{v_0}}: {}^{un}Ind_{P(F_{v_0})}^{G_{4n}(F_{v_0})} \tau_{v_0} \otimes |\det |^{\frac12}_{v_0} \boxtimes
\omega _{v_0} \to \pi.$  It follows easily from the definitions that
$$FC^\vartheta \circ  i_{\zeta^{v_0}}$$
factors through the Jacquet module
$
\mathcal{J}_{N_\ell, \vartheta}(\;^{un}
Ind_{P(F_{v_0})}^{G_{4n}(F_{v_0})} \tau_{v_0} \otimes |\det |^{\frac12}_{v_0} \boxtimes
\omega _{v_0} ).$
Propositions
\ref{p:JacVanishes1} and \ref{p:JacVanishes2} below each show that this Jacquet module
vanishes at approximately half of all places.  Inasmuch as vanishing at a single place
would suffice to prove global vanishing,
the result follows.
\end{proof}

A general character of $\descgroup{\ell}$  is of the form
\begin{equation}\label{e:GenChar}
\baseaddchar( c_1 u_{1,2} + \dots + c_{\ell-1} u_{\ell-1,\ell}
+ d_1 u_{\ell , \ell+1} + \dots + d_{4n+1-2\ell} u_{\ell, 4n+1-\ell}).\end{equation}
As described in section \ref{s:uniper}, the Levi $\desclevi{\ell}$ acts on the space of characters
of $N_\ell\quo.$  In order to define embeddings of the various forms of $G_{2n}$ into
$G_{4n+1},$ we need to make this more explicit.

First, we fix a specific isomorphism of $GL_1^\ell\times G_{4n-2\ell+1}$ with $L_\ell$
as follows.  As in section
\ref{rootDataSection},  let $e_0, \dots, e_{2n}$   and $e_0^*,\dots, e_{2n}^*$ denote the $\Z$-bases of $X(T(G_{4n+1}))$
and $X^\vee(T(G_{4n+1})),$ respectively.  Let $\hat e_0,\dots \hat e_{2n-\ell},$
and $\hat e_0^*, \dots, \hat e_{2n-\ell}^*$ denote the analogues for $G_{4n-2\ell+1}.$
We identify  $(\alpha_1, \dots, \alpha_\ell, \prod_{i=1}^{2n-\ell} \hat e_i^*(t_i) )
\in GL_1^\ell \times T(G_{4n-2\ell+1})$ with $\prod_{i=1}^\ell e_i^*(\alpha_i)\cdot \prod_{i=1}^{2n-\ell}
e_{i+\ell}^*(t_i) \in T(G_{4n+1}).$
In addition, we require that $g\in G_{4n-2\ell+1}$ be
identified with an element of $G_{4n+1}$ which projects to
$$\begin{pmatrix} I_\ell&&\\&\pr(g)&\\&&I_\ell\end{pmatrix}\in SO_{4n+1}.$$
Together, these requirements determine a unique identification.

Let $\underline d$ denote the column vector $^t(d_1, \dots, d_{4n+1-2\ell}).$
Suppose $\vartheta(u)$ is the character of $N_\ell$ given by
\eqref{e:GenChar},
and, for $h \in L_\ell,$ let
\begin{equation}\label{e:Action}
h\cdot \vartheta(u)= \vartheta(h^{-1}uh)
=\baseaddchar( ^hc_1 u_{1,2} + \dots + {}^hc_{\ell-1} u_{\ell-1,\ell}
+ {}^hd_1 u_{\ell , \ell+1} + \dots + {}^hd_{4n+1-2\ell} u_{\ell, 4n+1-\ell}).
\end{equation}
This is an  action of $L_\ell$ on the space of characters, and it
is easily verified that for $h$ identified with $(\alpha_1, \dots, \alpha_\ell, g),$
with $\alpha_1, \dots, \alpha_\ell \in GL_1(F)$ and $g \in G_{4n-2\ell+1}(F),$ we have
$$^hc_i = \frac{\alpha_{i+1}}{\alpha_{i}}\cdot  c_i,\;i=1,\dots, \ell-1, \quad \text{and}\quad
^h\underline d = \alpha_{\ell}^{-1}\cdot \pr(g)\cdot \underline d.$$

The above discussion amounts to an identification of the action of $L_\ell(F)$ on the space
of characters of $N_\ell\quo$ with a certain rational representation of
$L_\ell$ defined over $F,$ consisting of the direct sum of $\ell-1$ one dimensional
representations and a $(4n-2\ell+1)$-dimensional representation on which the
$G_{4n-2\ell+1}$ factor in $L_\ell$ acts via its ``standard'' representation.  We may consider this
rational representation over any field.
Over
an algebraically closed field there is an open orbit, which
consists of all those elements such that
$c_i \neq 0$ for all $i$ and $^t \underline{d} J \underline{d} \neq 0.$
Here,   $J$ is defined as in \ref{ss:GeneralNotation}.
Over a general field two such elements are in the same $F$-orbit
iff the two values of $^t \underline{d} J \underline{d}$ are in the
same square class.  Thus, this square class is an important invariant of the character
$\vartheta.$

\begin{defn} \label{d:Invt}
Let $\vartheta$ be the character of $N_\ell\quo$ given by
$$\vartheta(u)=\baseaddchar( c_1 u_{1,2} + \dots + c_{\ell-1} u_{\ell-1,\ell}
+ d_1 u_{\ell , \ell+1} + \dots + d_{4n+1-2\ell} u_{\ell, 4n+1-\ell}).$$
We denote the square class of $^t\underline d J \underline d$ by
$\Invt(\vartheta).$
We say that $\vartheta$ is in general position if $c_i\ne 0$ for $1\le i \le \ell-1$ and
$\Invt(\vartheta)\ne 0.$
 We denote the square class
consisting of the nonzero squares by $\square.$
\end{defn}

Clearly, a nonzero square class in $F$ may also be used to determine a quasi-split form of $G_{2n}.$
Indeed, the natural datum for determining a quasi-split group with $G$ such that $^LG^0=GSO_{2n}(\C)$ is a homomorphism  $\Gal(\bar F/F)\to \Aut(GSO_{2n}(\C))/\Inn(GSO_{2n}(\C)),$ which has
two elements.  Such homomorphisms are in one-to-one correspondence with quadratic characters
by class field theory, and this has been exploited in defining $G_{2n}^\chi$ above.  But they are
also in natural one-to-one correspondence with square classes in $F^\times,$
and this parametrization will be more convenient for the next part of the discussion.
\begin{defn}\label{d:G(sqClass)}
Let ${\bf a}$ be a square class in $F^\times.$  Let  $F(\sqrt{\bf a})$ denote the smallest extension of
$F$ in which the elements of ${\bf a}$ are squares, and for $a \in F^\times,$
let $F(\sqrt{a})=F(\sqrt{\{a\}}),$
where $\{a\}$ is the square class of $a.$  Let $G_{2n}^{\bf a}$ denote the
quasi-split form of $G_{2n}$  such that the associated
map $\Gal(\bar F/F)\to \Aut(GSO_{2n}(\C))/\Inn(GSO_{2n}(\C))$
factors through
$\Gal(F(\sqrt{\bf a})/F).$
\end{defn}
\begin{rmk}
Of course, if ${\bf a}=\square,$ then $F(\sqrt{\bf a})=F$ and $G_{2n}^{\bf a}$  is
just the split group $G_{2n}.$
\end{rmk}

\begin{lem}\label{l:embeddings}
\begin{enumerate}
\item\label{twoConnComps}
If $\vartheta$ is a character of $N_\ell$ in general position, then
the stabilizer $L_\ell^\vartheta$ (cf. $M^\vartheta$ in definition \ref{def:fc})
has two connected components
\item\label{stabIsoRightThing}
The identity component $(L_\ell^\vartheta)^0$ is isomorphic over $F$ to
$G_{4n-2\ell}^{\Invt(\vartheta)}.$
\end{enumerate}
\end{lem}
\begin{proof}
Identify $(\alpha_1, \dots, \alpha_\ell, g)\in GL_1^\ell\times G_{4n-2\ell+1}$ with
an element of $L_\ell$ as above.

 The identity component of $L_\ell^\vartheta$
 consists of those $(\alpha_1, \dots, \alpha_\ell, g)$
such that $\alpha_i=1$ for all $i$ and $g$ fixes the vector in the standard representation
obtained from $\vartheta.$  The other consists of those such that
$\alpha_i=-1$ for all $i,$ and $g$ maps the vector in the standard representation
obtained from $\vartheta$ to its negative (which is the only scalar multiple of the
same length).
This proves \eqref{twoConnComps}.

We turn to   \eqref{stabIsoRightThing}.  First suppose $\Invt (\vartheta)=\square.$
It suffices to consider the specific character
$\Psi_\ell$ defined by
$$\Psi_\ell( u) =
\baseaddchar(u_{12}+\dots+u_{\ell-1,\ell}+u_{\ell,2n+1}).$$
For this character, the column vector $\underline d$ is $v_1:=^t(0,\dots, 0,1,0,\dots, 0).$
 It is easily checked
that the stabilizer of this point in $SO_{4n-2\ell+1}$ is isomorphic to the split form
of $SO_{4n-2\ell}.$  In addition, the stabilizer in $G_{4n-2\ell+1}$ contains a split
torus of rank $2n-\ell+1,$ and hence is a split group.  An element of $\maxunip$
fixes $v_1,$ if and only if it satisfies $u_{i,2n-\ell+1}=0$
for $i=1$ to $2n-\ell.$  From this we easily compute the based root datum of the
stabilizer of $v_1$ and find that it is the same as that of $G_{4n-2\ell}.$

To complete the proof of \eqref{stabIsoRightThing}, let $a$ be a non-square in $F^\times,$ and let $v_a=^t(0,\dots, 0,1,0,\frac a2,0,\dots, 0)\in F^{4n-2\ell+1}$
(nonzero entries in positions $2n-\ell$ and $2n-\ell+2$ only).
Let $\Psi_\ell^a$ be the character of $N_\ell\quo$ corresponding to $c_i=1\forall i$ and
$\underline d =v_a.$
The stabilizers
of $\Psi_\ell^a$ and $\Psi_\ell$ are conjugate over the quadratic extension $E$ of
$F$ obtained by adjoining a square root of $a.$  Indeed, let $\sqrt{a}$ be an element of
$E$ such that $(\sqrt{a})^2=a.$
Suppose
$$\pr(h_a)=\begin{pmatrix}\sqrt{a}^{-1} I_{2n-1}&&\\&M_{\sqrt{a}}&\\&&\sqrt{a} I_{2n-1}\end{pmatrix},
\quad \text{where}\quad
M_{\sqrt{a}}=\begin{pmatrix}
-\frac 1{2\sqrt{a}}&\sqrt{a}^{-1}&\sqrt{a}^{-1}\\ \frac12 &0&1\\
\frac {\sqrt{a}} 4& \frac {\sqrt{a}} 2& -\frac {\sqrt{a}} 2\end{pmatrix}.$$
Then $h_a\cdot \Psi_\ell =\Psi_\ell^a.$
For each $a,$ fix an element $h_a$ as above for use throughout.

Clearly $(L_\ell^{\Psi_\ell^a})^0=h_a(L_\ell^{\Psi_\ell})^0h_a^{-1}.$  The image of this group under $\pr$
is isomorphic over $F$ to
the non-split quasisplit form of $SO_{4n-2\ell}$ corresponding the square class of $a.$  It follows
that $(L^{\Psi_\ell^a}_n)^0$ is isomorphic over $F$ to
the non-split quasisplit form of $G_{4n-2\ell}$ associated to the
square class of $a.$
\end{proof}

In the course of the preceding proof, we have seen that it is enough to consider one
conveniently chosen representative from each $F$-orbit of characters in general position.
However, it is generally more convenient to make definitions for general $a\in F^\times$
than it is to choose representatives for the square classes
in $F^\times.$
\begin{defn}\label{d:SpecificCharacters}
Take $a\in F^\times,$ and
let $\Psi_\ell^a$ be the character of $\descgroup{\ell}$ defined
by
$$\Psi_\ell^a( u ) =\baseaddchar (u_{12}+\dots+u_{\ell-1,\ell}+ u_{\ell, 2n} +
\frac a2 u_{\ell, 2n+2} ).$$
We also keep the notation
$$\Psi_\ell( u ) =\baseaddchar (u_{12}+\dots+u_{\ell-1,\ell}+ u_{\ell, 2n+1}).$$
Then the orbit of $\Psi_\ell^a$ is determined by the square class of $a.$  The character $\Psi_\ell$ is
in the same orbit as $\Psi_\ell^1.$
\end{defn}
Note that for any given square class ${\bf a}$ we have many conjugate embeddings of
$G_{2n}^{\bf a}$ into $G_{4n+1}$: one for each element $a$ of ${\bf a}.$
\begin{defn}\label{d:G(number)}
For each  element $a$ of $F^\times,$ we let  $G_{2n}^a$ denote $(L_n^{\Psi_n^a})^0.$
It is a subgroup of $G_{4n+1},$ which is isomorphic over $F$ to $G_{2n}^{\{a\}},$
where $\{a\}$ is the square class of $a.$
\end{defn}
\begin{lem}\label{l:G^aBorel}Assume $\{a\}\ne \square.$  Then,
\begin{enumerate}
\item\label{unip}
An element $u$ of $\maxunip$ is in $G_{2n}^a$ iff it satisfies $u_{ij}=0$ for $i\le n$ or
$i=2n,$ and $u_{i,2n}=-\frac a2 u_{i,2n+2}$ for $n< i<2n.$
The set of such elements $u$ is equal to
$h_a(\maxunip \cap (L_n^{\Psi_n})^0)h_a^{-1},$ and is a maximal unipotent subgroup of
$G_{2n}^a.$
\item \label{Fsplit}
An element $t=\prod_{i=0}^{2n}e_i^*(t_i)$ of $T(G_{4n+1})$ is in $G_{2n}^a$ iff it satisfies
$t_i=1$ for $0<i\le n,$ and $i=2n.$ The set of such  $t$ is a maximal $F$-split torus of $G_{2n}^a.$
\item \label{maxtor}
There is a maximal torus of $G_{2n}^a$ which contains the above maximal $F$-split torus and
is contained in the standard Levi of $G_{4n+1}$ whose unique positive root is the short simple
root $e_n.$  Its set of $F$ points is equal to
$$\left\{h_a t h_a^{-1}: t=\prod_{i=1}^{n-1}e_{n+i}^*(t_i) e_{2n}^*(x\cdot \bar x^{-1}) e_0^*(\bar x),\;
t_1, \dots, t_{n-1}\in F^\times, \; x \in F(\sqrt{a})^\times\right\},$$
where $\overline{\phantom{x}}$ denotes the action of the nontrivial element of $\Gal(F(\sqrt{a})/F).$
\end{enumerate}
If $\{a\}=\square,$ then \eqref{unip} remains true, while
$$\left\{h_a t h_a^{-1}: t=\prod_{i=0}^{n}e_{n+i}^*(t_i) \right\},$$
is a maximal torus, and is $F$-split, since $h_a$ has entries in $F.$
\end{lem}
\begin{rmk}\label{uniformMaxTor}
We may write an element of our maximal torus as
$$\left\{ h_a \prod_{i=1}^{n-1}e_{n+i}^*(t_i) \cdot e_{2n}^*\left((x+y\sqrt{a})\cdot (x-y\sqrt{a})^{-1}\right)
e_0(x-y\sqrt{a})h_a^{-1} : t_i \in F, \; x, y\in F, x^2-ay^2\ne 0\right\},$$
regardless of $\{a\}.$
\end{rmk}
\begin{proof}
Item \eqref{unip} is easily checked.  (Recall that $\pr$ is an isomorphism on $\maxunip.$)
Similarly, it is easily checked that an element $t$ of $T(G_{4n+1})$ stabilizes the
specified character iff
$t_1=\dots =t_n=t_{2n}=\pm 1.$  As noted in the proof of Lemma \ref{l:embeddings}, if they
are all minus 1, then this element is in the other connected component of $L_n^{\Psi_n^a}.$

Recall that $(L_n^{\Psi_n})^0,$ with $\Psi_n$ as in Definition \ref{d:SpecificCharacters}
is isomorphic to $G_{2n}.$  There is an ``obvious'' choice of isomorphism $\inc:G_{2n}\to
(L_n^{\Psi_n})^0,$ such that
$$\inc \circ \bar e_i^*=\begin{cases}e_0^*&i=0,\\ e_{n+i}^*&1\le i \le n,
\end{cases}\qquad \text{and}\qquad
\inc(u)_{ij}=\begin{cases}
0&i\le n, \text{ or }j=2n+1,\\
u_{i-n,j-n}&i>n, \; j< 2n+1,\\
u_{i-n,j-n-1}& i>n, \; j>2n+1.
\end{cases}
$$
Here, we have used $e_i^*$ for elements of the $\Z$-basis of the cocharacter lattice of
$G_{4n+1}$ and $\bar e_i^*$ for elements of that of $G_{2n}.$
It follows from the definitions that conjugation by $h_a$ is an isomorphism of $G_{2n}^a$
with $(L_n^{\Psi_n})^0,$ which is defined over $F(\sqrt{a}).$  This yields an identification of
the maximal $F$-split torus of $G_{2n}^{\{a\}}$ as computed in section
\ref{rootDataSection} with the $F$-split torus in item \eqref{Fsplit}.

Clearly $h_a\cdot \inc(T(G_{2n}))\cdot h_a^{-1}$ is a maximal torus of $G_{2n}^a.$ The fact that
an element is of the form specified in item \eqref{maxtor} of the present lemma follows from the
action of $\Gal(\bar F/F)$ on the lattice of cocharacters computed in section \ref{rootDataSection}.
\end{proof}

Finally, following the works of Ginzburg-Rallis-Soudry we arrive
at the {\it descent construction}.

\begin{defn}
Let $$DC_\omega^a(\tau) = FC^{\Psi_{n}^a} \residuerep.$$
 It
 is a space of smooth functions $G_{2n}^a\quo \to \mathbb C,$ and
 affords a
representation of the group $G_{2n}^a(\A)$
acting by right translation, where we have identified
 $G_{2n}^a$ with  the identity component of $L_{n}^{\Psi_{n}^a}.$
\end{defn}
\begin{defn}
We say that a square class ${\bf a}$ in $F^\times$ and a character $\chi$ are {\rm compatible}
if they correspond to the same homomorphism from $\Gal(\bar F/F)$ to the
group with two elements.  We say that an element $a$ of $F^\times$ and a character $\chi$
are compatible if $\chi$ is compatible with the square class of $a.$
\end{defn}

\subsection{Vanishing of incompatible descents}

\begin{thm}\label{t:DescentsWithWrongSqClassVanish}
Let $\omega$ be a Hecke character.
Let $\tau=\tau_1\boxplus \dots \boxplus \tau_r$ be the isobaric sum of
distinct $\omega^{-1}$-orthogonal unitary cuspidal automorphic representations
of $GL_{2n_1}(\A),$ \dots, $GL_{2n_r}(\A),$ respectively.  For $i=1$ to $r$ let $\omega_{\tau_i}$
denote the central character of $\tau_i$ and let
$\chi_i=\omega_{\tau_i}/\omega^{n_i},$ which is quadratic.  Let $\chi=\prod_{i=1}^r \chi_i.$
Suppose that $\chi$ and $a$ are not compatible.  Then $DC_\omega^a(\tau)=\{0\}.$
\end{thm}
\begin{proof}
As in Theorem \ref{t:DeeperDescentsVanish}, it suffices to prove the vanishing of the corresponding
twisted Jacquet module
of $\Ind_{P(F_v)}^{G_{4n+1}(F_v)}\tau_v\otimes|\det|^{\frac 12} \boxtimes \omega_v$
 at a single unramified place $v.$    The vanishing follows from Proposition
\ref{p:JacVanishes1},
if there is an unramified  place $v$ such that $\chi_v$ is trivial and $a$ is not a square, and from
Proposition
\ref{p:JacVanishes2}
if there is an unramified place $v$ such that $\chi_v$ is nontrivial and $a$ is a square.
If $\chi$ and $a$ are incompatible, then there is at least one unramified place at which one
of these cases occurs.
 \end{proof}

\subsection{Main Result}

\begin{thm}\label{t:maintheorem} Let $\omega$ be a Hecke character.
Let $\tau=\tau_1\boxplus \dots \boxplus \tau_r$ be the isobaric sum of
distinct $\omega^{-1}$-orthogonal unitary cuspidal automorphic representations
of $GL_{2n_1}(\A),$ \dots, $GL_{2n_r}(\A),$ respectively.  For $i=1$ to $r$ let $\omega_{\tau_i}$
denote the central character of $\tau_i$ and let
$\chi_i=\omega_{\tau_i}/\omega^{n_i},$ which is quadratic.  Let $\chi=\prod_{i=1}^r \chi_i.$
Then
\begin{enumerate}
\item
$DC_\omega^a(\tau)$ is nontrivial if and only if $\chi$ and $a$ are compatible,
\item
 If $\chi$ and $a$ are compatible then
 the space $DC_\omega^a(\tau)$ is a
nonzero,  cuspidal representation
of $G_{2n}^a(\A),$ with central character $\omega.$ Furthermore, the representation
$DC_\omega^a(\tau)$  supports a nonzero Whittaker integral
for the generic character of  $\maxunip(\A)\cap G_{2n}^a(\A)$ given by
$$u \mapsto \baseaddchar\left(
\sum_{i=1}^{2n-2} u_{i,i+1} + u_{2n-1,2n+2} \right).$$
\item If $\sigma$
 is any irreducible automorphic
representation contained
in $DC_\omega^a(\tau),$ then $\sigma$ lifts weakly to $\tau$ under the
map
$r.$
\end{enumerate}
\end{thm}
\begin{rmk}
Since $DC_\omega(\tau)$ is nonzero and cuspidal, there exists at least one irreducible component
$\sigma.$  In the case of special orthogonal groups, one may show (\cite{So-Paris}, p. 342, item 4) that
the descent module is {\it in the $\psi$-generic spectrum}
for a suitable choice of $\psi$ (cf. section \ref{s:NotionsOfGenericity}).
It follows that {\it all} of the irreducible components are
distinct and  globally $\psi$-generic.
This is done using the Rankin-Selberg integrals of \cite{G-PS-R},\cite{Soudry-SO2n+1Local}.
In the odd case, one may also show
(\cite{GRS4}, Theorem 8, p. 757, or \cite{So-Paris} page 342, item 6) using the results of
\cite{jiangsoudry} that the descent module is irreducible.  This does {\bf not}
extend to the even case, even for special orthogonal groups, because the construction
actually yields a representation of the full stabilizer-- which is
isomorphic to the full orthogonal group.   (Cf. Proposition 7.0.20.)
\end{rmk}

\subsection{Proof of the main theorem (Even case)}
The statements are proved by combining relationships between unipotent
periods and knowledge about $\residuerep.$

\begin{enumerate}

\item {\bf Genericity and non-vanishing}
For $a\in F^\times,$ we let $(U_1^a, \psi_1^a)$ denote the unipotent period obtained by composing the
period
$(N_n,\Psi_n^a),$ used in defining the descent to $G_{2n}^a,$
(embedded into $\ourgroup_{4n+1}$ as the stabilizer of $\Psi_n^a$)
 with a period which defines a  Whittaker integral on this group.
 Specifically,  $U_1^a$ is the subgroup of the standard maximal unipotent defined by
the relations $u_{i,2n}=-\frac a 2u_{i,2n+2}$ for $i = n+1$ to $2n-1,$
as well as $u_{2n,2n+1}=0,$
and
$$\psi_1^a( u ) = \baseaddchar( u_{1,2} + \dots +u_{n-2,n-1} + u_{n-1, 2n} + \frac a2 u_{n-1,2n+2}
+ u_{n, n+1} + \dots + u_{2n-1, 2n}).$$
The definitions of $U_1^a$ and $\psi_1^a$ make sense also in the case when $a=0,$
although in that case there is no interpretation in terms of a descent.  We use this period
in that case also.

Next, let $U_2$ denote the subgroup of the standard maximal unipotent defined
by $u_{2n,2n+1}=0,$ and $u_{12}=u_{34}=\dots=u_{2n-1,2n}.$  For all $a \in F,$
we may define a character of this group by the formula
$$\psi_2^a(u)=\baseaddchar\left(\sum_{i=1}^{2n-2}u_{i,i+2}+u_{2n-1,2n+2}+\frac a2
u_{2n-1,2n}\right).$$

Finally, let $U_3$ denote the maximal unipotent, and $\psi_3$ denote
$$\psi_3( u) = \baseaddchar( u_{1,2} + \dots + u_{2n-1,2n}).$$  Thus $(U_3, \psi_3)$ is the
composite of the unipotent period defining the constant term along the Siegel
parabolic, and one which defines a Whittaker integral on the
Levi of this parabolic.  By Theorem \ref{t:EisensteinSeriesProperties}
\eqref{residueSupportsPeriod} this
period is {\sl not} in $\uniper^\perp( \residuerep).$

In the appendices, we show
\begin{enumerate}
\item{$(U_1^a, \psi_1^a) \sim (U_2, \psi_2^a),$ for all $a\in F,$ in Lemma \ref{u1u2lemma}, }
\item{$(U_2,\psi_2^0)\in \langle \{(U_2,\psi_2^a):a\in F^\times \}\rangle,$ in Lemma \ref{thetaLemma},
and }
\item{$ (U_3, \psi_3) \in \langle
(U_2, \psi_2^0) , \{ (N_\ell, \vartheta): n <\ell < 2n \text{ and } \vartheta \text{ in general position.} \}
\rangle$ in Lemma \ref{u2u3deep}.
}
\end{enumerate}

By Theorem \ref{t:DeeperDescentsVanish} $(N_\ell , \vartheta) \in \uniper^\perp( \residuerep)$
for all $n < \ell < 2n$ and $\vartheta$ in general position.  It follows that at least
one of the periods $(U_1^a , \psi_1^a)$ is not in
$\uniper^\perp( \residuerep).$
It follows from theorem \ref{t:DescentsWithWrongSqClassVanish} that
$(U_1^a , \psi_1^a)$ vanishes for $a$ incompatible
with $\chi,$
so it must not vanish for some $a$ compatible with
$\chi.$
This establishes genericity (and hence nontriviality)
of the corresponding descent module $DC_\omega^a(\tau).$  The spaces $DC_\omega^a$ for $a$
compatible with $\chi$ may all be identified with one
another via suitable isomorphisms among the groups
$G_{2n}^a,$ and so $DC_\omega^a$ is nonzero
and generic for all $a$ compatible with $\chi.$

\item {\bf cuspidality:}
\\
Turning to cuspidality, we prove in the appendices an identity relating:
 \begin{itemize}
 \item{ Constant terms on $\ourgroup_{2n}^a,$ }
 \item{ Descent periods in $\ourgroup_{4n+1},$}
 \item{ Constant terms on $\ourgroup_{4n+1},$}
 \item{ Descent periods on $\ourgroup_{4n - 2k+1 },$ embedded in $\ourgroup_{4n}$ as a subgroup of
  a Levi.}
 \end{itemize}
To formulate the exact relationship we introduce some notation for the
maximal parabolics of GSpin groups.

 The group $\ourgroup_{4n+1}$ has one standard maximal parabolic having Levi
 $GL_i \times \ourgroup_{4n-2i+1}$ for each value of $i$ from $1$ to $2n.$  Let
 us denote the unipotent radical of this parabolic by $V_i.$  We denote the
 trivial character of any unipotent group by ${\bf 1}.$

 For any square class ${\bf a},$ the group $\ourgroup_{2n}^{\bf a}$ has one standard
 maximal parabolic having Levi
 $GL_k \times \ourgroup_{2n-2k}^{\bf a}$ for each value of $k$ from $1$ to $n-2.$  We denote
 the  unipotent radical of this parabolic by $V^{2n}_k.$
The split group $\ourgroup_{2n}=\ourgroup_{2n}^\square$
also has two parabolics with Levi isomorphic to $GL_{n} \times
 GL_1.$  One has the property that $e_{n-1}-e_n$ is a root of the Levi, and the other does not.
 Let us denote the unipotent radical of this first parabolic by $V^{2n}_n.$  Then the
 unipotent radical of the other is $^\dagger V_n^{2n},$ where $\dagger$ is the
 outer automorphism of $G_{2n}$ which reverses the last two simple roots while fixing the
 others.  In a nonsplit quasisplit form of $G_{2n},$ there is a parabolic subgroup with Levi
 isomorphic to the product of $GL_{n-1}$ and a nonsplit torus which is maximal.
 (The corresponding parabolic in the split case is not maximal.)  We denote
 its unipotent radical by $V_{n-1}^{2n}.$

  We prove in  Lemma \ref{l:cuspidality-unip-id} that, for $1 \le k \le n-1,$
$(V_k^{2n}, {\bf 1} ) \circ (N_{n} , \Psi_{n}^a)$
is contained in
$$
\langle
(N_{n+k} , \Psi_{n+k}) , \{
(N_{n+ j}, \Psi_{n+j}^a)^{(4n-2k+2j+1)}
\circ ( V_{k-j} , {\bf 1} ): \; \; 1 \leq j < k
\}
\rangle,$$
where $(N_{n+ j}, \Psi_{n+j}^a)^{(4n-2k+2j+1)}$ denotes the descent period, defined as above,
but on the group $G_{4n-2k+2j+1},$ embedded into $G_{4n+1}$ as a component of the
Levi with unipotent radical $V_{k-j}.$

Now suppose that $a$ is a square.  Then both $(V_n^{2n}, {\bf 1} ) \circ ( N_{n} , \Psi_{n}^a)$
and $(^\dagger V_n^{2n}, {\bf 1} ) \circ ( N_{n} , \Psi_{n}^a)$ are in
$$
\langle
(N_{2n} , \Psi_{2n}) , \{
(N_{n+ j}, \Psi_{n+j}^a)^{(2n+2j+1)}
\circ ( V_{n-j} , {\bf 1} ): \; \; 1 \leq j < n
\}
\rangle.$$
Indeed, the two periods are actually conjugate in $G_{4n+1},$ so it suffices to consider only one
of them.

By Theorem \ref{t:DeeperDescentsVanish}
$(N_{n+k} , \Psi_{n+k}^a) \in \uniper^\perp(\residuerep)$ for $k=1$ to $n.$
Furthermore, for   $k,j$ such that $1 \leq j < k \leq n$
the function $E(f)(s)^{(V_{k-j},{\bf 1})}$ may be expressed
 in terms of Eisenstein series on $GL_{k-j}$ and
$G_{4n-2k+2j},$
using Proposition II.1.7 (ii) of  \cite{MW1}.
What we require is the following:
\begin{lem}
For all $f \in V^{(2)}(\underline s,\bigotimes_{i=1}^r \tau \boxtimes \omega) $
$$\left.E_{-1}(f)^{(V_{k-j},{\bf 1})}\right|_{G_{4n-2k+2j+1}(\A)}
 \in \bigoplus_S
\mathcal{E}_{-1}( \tau_{S},\omega),$$
where the sum is over subsets $S$ of $\{1, \dots, r\}$ such that
$\sum_{i \in S }2n_i = 2n-k+j,$ and, for each such $S$,
$\mathcal{E}_{-1}( \tau_{S},\omega)$ is the space of functions  on
$G_{4n-2k+2j+1}(\A)$ obtained by applying the construction
of $\mathcal{E}_{-1}( \tau,\omega)$ to $\{ \tau_i : i \in S\},$
instead of $\{ \tau_i: 1 \le i \le r\}.$
\end{lem}
Once again, this is immediate from  \cite{MW1}
Proposition II.1.7 (ii).

Applying Theorem \ref{t:DeeperDescentsVanish}, with $\tau$
replaced by $\tau_S$ and $2n$ by $2n-k+j,$ we deduce
$$(N_{n+ j}, \Psi_{n+j})^{(4n-2k+2j+1)} \in \mathcal{U}^\perp\left(
\mathcal{E}_{-1}( \tau_{S},\omega)\right) \quad \forall S,$$
and hence $(N_{n+ j-1}, \Psi_{n+j-1})^{(4n-2k+2j)}
\circ ( V_{k-j} , {\bf 1} ) \in \uniper^\perp( \residuerep).$
This shows that any nonzero function appearing in any of the spaces $DC^a_\omega(\tau)$
must be cuspidal.
Such a function is also easily seen to be of uniformly moderate growth, being the
integral of an automorphic form over a compact domain.
In addition, such a function is easily seen to have central character $\omega,$ and
any function with these properties is necessarily square integrable modulo the
center (\cite{MW1} I.2.12).
It follows that each of the spaces $DC^a_\omega(\tau)$ decomposes discretely.

\item {\bf Verification of weak lifting: unramified parameters:}

Now, suppose $\sigma \iso \otimes_v' \sigma_v$ is an irreducible representation
which is contained in $DC_\omega^a(\tau).$
Let $p_\sigma$ denote the natural projection $DC_\omega^a (\tau) \to \sigma.$
Once again, by Theorem \ref{t:EisensteinSeriesProperties} \eqref{ResidueIsL2}, the
representation $\residuerep$ decomposes discretely.  Let $\pi$ be an irreducible
component of $\residuerep$ such that the restriction of $p_\sigma \circ FC^{\Psi_n^a}$ to
$\pi$ is nontrivial.  As discussed
previously in the proof of
Theorem \ref{t:DeeperDescentsVanish}, at all but finitely many $v,$
 $\tau$ is unramified at $v$ and
furthermore, $\pi_v$ is the unramified constituent
$^{un}Ind_{P(F_v)}^{G_{4n+1}(F_v)} \tau_v \boxtimes \omega _v \otimes |\det|_v^{\frac 12}$
of
$Ind_{P(F_v)}^{G_{4n+1}(F_v)} \tau_v \boxtimes \omega _v \otimes |\det|_v^{\frac 12}.$
If $v_0$ is such a place, the map $p_\sigma \circ FC^{\Psi_n^a} \circ i_{\zeta^{v_0}},$ with
$ i_{\zeta^{v_0}}$ defined as in Theorem \ref{t:DeeperDescentsVanish},
factors through
$\mathcal{J}_{N_{n},\Psi_{n}^a}
\left( \;
^{un}Ind_{P(F_{v_0})}^{G_{4n+1}(F_{v_0})} \tau_v\otimes |\det|_v^{\frac 12}\boxtimes \omega _v
\right),$
and gives rise
to a $\ourgroup_{2n}^a(F_{v_0})$-equivariant map from this Jacquet-module onto
$\sigma_{v_0}.$

 To pin things down precisely,
assume that $\tau_v$ is the unramified component of $Ind_{B(GL_{2n})(F_v)}^{GL_{2n}(F_v)}
\mu,$ and let $\mu_1 , \dots, \mu_{2n}$ be defined as in the proof of Lemma \ref{l:combolemma}.
By Lemma \ref{l:combolemma}, we may assume without loss of generality that
$\mu_{2n+1-i}= \omega  \mu_i^{-1}$ for $i = 1$ to $n-1,$ and that either
$\mu_n=\omega\mu_{n+1}^{-1},$ or $\mu_n^2=\mu_{n+1}^2=\mu_n\mu_{n+1}\chi_{un}=\omega$
(with $\chi_{un}$ defined as in the lemma).
Furthermore, suppose that $\chi_v$ is the local component at $v$ of the global quadratic character
obtained from $\tau$ and $\omega$ as in
the statement of the theorem.  Then either $\chi_v$ is trivial and  $\mu_n=\omega\mu_{n+1}^{-1},$ or
$\chi_v=\chi_{un}$ and
$\mu_n^2=\mu_{n+1}^2=\mu_n\mu_{n+1}\chi_{un}=\omega.$

Recall that a basis for the lattice of $F$-rational cocharacters of the maximal torus of
$G_{2n}^a$ fixed in Lemma \ref{l:G^aBorel} is given by
$$\{ e_{n+i}^*:1\le i <n\} \cup \{e_0^*\} \cup \{ e_n^*, \text{\it if  }a\text{ is a square}\}.$$
Observe that when $a$ is not a square in $F,$ it is a square in $F_v$ for many
unramified $v,$ and that the cocharacter $e_n^*$ is $F_v$-rational at such $v.$

In proposition \ref{p:JFcompNonSplit}, we show that in the nonsplit case
$$\mathcal{J}_{N_{n},\Psi_{n}}\left(
^{un}Ind_{P(F_v)}^{G_{4n+1}(F_v)} \tau_v \boxtimes \omega _v \otimes |\det|_v^{\frac 12}
\right)$$
is isomorphic as a $\ourgroup_{2n}^a(F_v)$-module to
a subquotient of a principal series representation $\pi_v$ of $\ourgroup_{2n}^a(F_v)$
such that the corresponding parameter $t_{\pi,v}$ maps to the parameter $t_{\tau,v}$
under $r.$  In the split case (proposition \ref{p:JFcompSplit}), we obtain instead a direct sum of two principal series
representations, but {\em both} have parameters which map to $t_{\tau,v}.$  It follows
that $\tau$ is the weak lift of $\sigma$ associated to the map $r.$
\end{enumerate}

\section{Appendix III: Preparations for the proof of Theorem \ref{t:EisensteinSeriesProperties}}
\label{section: Preparations for the proof of Theorem t:EisensteinSeriesProperties}

In this section we review some standard arguments
by which the presence or
absence of a singularity of an Eisenstein series reduces to the
presence or absence of a singularity of a relative rank one
intertwining operator.

To do so, we recall the set
$$W(M):= \left\{ w \in W_{G_{4n+1}}\left| \begin{array}{l}
w \text{ is of minimal length in }w\cdot W_M\\
wMw^{-1}\text{ is a standard Levi of }G_{4n+1}
\end{array}\right.\right\}.$$
It will be convenient and harmless to treat the elements of $W(M)$ as though they
were elements of $G_{4n+1}(F),$ rather than repeatedly choose representatives
and remark the independence of the choice.

\subsection{Intertwining operators}

For each $w \in W(M),\;\underline s\in \C^r,$ we define
$P^w$ to be the standard parabolic with Levi $wMw^{-1}.$
For $\underline s$ such that $s_r$ and $s_i-s_{i+1}, i =1$ to $r-1$
are all sufficiently large, the integral
$$M(w, \underline s)f(g):= \int_{\maxunip \cap w \overline \maxunip w^{-1} \quo}
 f(\underline s)(w^{-1} ug)\; du$$
 converges (\cite{MW1}, II.1.6), defining an operator
 $M(w, \underline s)$ from  $V^{(2)} (\underline s, \bigotimes_{i=1}^r  \tau_i \boxtimes \omega)$
 to a space of functions which is easily verified to afford a realization of
 $$\Ind_{P^w(\A)}^{\ourgroup_{4n+1}(\A)} \left((\bigotimes_{i=1}^r  \tau_i\otimes |\det{}_i|^{s_i} )
 \boxtimes \omega\right)\circ Ad(w^{-1}).$$
 Here, $\left((\bigotimes_{i=1}^r  \tau_i\otimes |\det{}_i|^{s_i} )
 \boxtimes \omega\right)\circ Ad(w^{-1}),$ denotes the representation of $wMw^{-1}$
 obtained by composing the representation $\left(\bigotimes_{i=1}^r  \tau_i\otimes |\det{}_i|^{s_i}
 \boxtimes \omega\right)$ of $M$ with conjugation by $w^{-1}.$
 We denote this latter space of functions by
 $V_w^{(2)}(\underline s, \bigotimes_{i=1}^r  \tau_i \boxtimes \omega).$
 Then $M(w, \underline s)f(g)$ has meromorphic continuation to $\C^r.$
 (\cite{MW1},IV.1.8(b).)

 It may be helpful also to review the sorts of singularities which Eisenstein series and
 intertwining operators have-- lying along so-called ``root hyperplanes.''
 (cf. \cite{MW1}, IV.1.6)
 We defer
 the notion of ``{\it root} hyperplane'' until later.  For now, we allow arbitrary hyperplanes in
 $\C^r,$ defined by equations of the form $l(\underline s)=c,$ with $l$ a linear functional
 $\C^r\to \C$ and $c$ a constant.  Then for any bounded open set $U\subset \C^r,$ there exist
 a finite number of distinct hyperplanes $H_1, \dots, H_N,$ which ``carry'' the singularities
 of the Eisenstein series and intertwining operators in $U,$ in the following sense.
 For each $i$ fix $l_i, c_i$ such that $H_i=\{ \underline s\in \C^r\;|\; l_i(\underline s)=c_i\}.$
 Then for each $i$ there is a non-negative integer $\nu(H_i)$ such that
 \begin{equation}\label{e:hyperplanes}
 \prod_{i=1}^N (l_i(\underline s)-c_i)^{\nu(H_i)}E(f)(g)(\underline s)
 \end{equation}
 continues to a function holomorphic on all of $U.$  Covering $\C^r$ with bounded open
 sets and taking a union, we obtain an infinite, but {\it locally } finite, set of hyperplanes which
 carry all the singularities of the Eisenstein series and intertwining operators.  The same hyperplane
 $H$ will of course occur more than once.  It is easily verified that the minimal exponent $\nu(H)$
 appearing in  \eqref{e:hyperplanes}  is the same each time.   Thus we may speak of whether
 an Eisenstein series or intertwining operator does or does not have a pole along $H,$
 and of the order of the pole.

 One may define ``analytic/meromorphic continuation'' for functions taking values in
 Fr\'echet spaces  of locally $L^2$ functions and the like
  (\cite{MW1} I.4.9, IV.1.3)
 of  functions and operators.  In this case, outside of the domain of convergence, one's functions are defined only up to $L^2$ equivalence.  However, in view of \cite{MW1}, I.4.10, one has a unique smooth
 representative for the class.  For us it will be more convenient simply to adopt the convention
 that when we say the Eisenstein series has a pole along $H,$ we mean for some $f,g.$

 Now let us state the relationship between poles of Eisenstein series and intertwining
 operators, which we prove in section \ref{ss: proof of ES to IO}.
\begin{prop}\label{p:EStoIO}
 For $f \in V^{(2)}(\underline s,  \bigotimes_{i=1}^r  \tau_i
 \boxtimes \omega),$
 there exists
 $g \in G_{4n+1}(\A)$ such that $E(f)(g)$ has a pole along $H$ if and only
 if there exist $w \in W(M), g' \in G_{4n+1}(\A)$ such that $M(w, \underline s)f(g')$ has a pole
 along $H$.
 \end{prop}

 The same construction can be performed with the Levi $M$ replaced by $wMw^{-1},$
 yielding an operator
 $$M_w(w',w\cdot \underline s):
  V_w^{(2)}(\underline s, \bigotimes_{i=1}^r  \tau_i \boxtimes \omega)
\to  V_{w'w}^{(2)}(\underline s, \bigotimes_{i=1}^r  \tau_i \boxtimes \omega),$$
for each $w' \in W(wMw^{-1}).$  Furthermore, one has for all $f,g,$
 the equality of meromorphic functions
 $$M_w( w' , w\cdot \underline s)
 \circ M(w, \underline s) f(g) = M(w'w,\underline s)f(g)$$
 (\cite{MW1}, II.1.6, IV.4.1).  (For now, the reader may think of ``$w\cdot \underline s$''
 simply as a notational contrivance.  We shall give it a precise meaning below.)

\subsection{Reduction to relative rank one situation}

 Next we wish to describe the decomposition of $w\in W(M)$ as a product of elementary
 symmetries, as in \cite{MW1} I.1.8.    The lattice
 $X(Z_M)$ of rational characters of the center of $M$ has a unique basis
 $\{ e_0, \varepsilon_1, \dots, \varepsilon_r\},$ with the property that for each $i=1, \dots m,$
 there exists $j \in \{1, \dots , r\}$ such that the restriction of $e_i$ as in \ref{rootDataSection}
 to $Z_M$ is $\varepsilon_j.$
The set of restrictions of positive roots of $G_{4n+1}$ to $Z_M$ is
\begin{equation}\label{restns of pos roots}
\{0\} \cup \{\varepsilon_i-\varepsilon_j: 1\le i<j\le r\}\cup
\{\varepsilon_i: 1\le i \le r\}\cup
\{\varepsilon_i+\varepsilon_j: 1\le i<j\le r\}\cup \{2\varepsilon_i: 1\le i\le r\}.
\end{equation}
Let $\Phi^+(Z_M)$
 denote the set obtained by excluding zero.  For $\alpha
\in \Phi^+(Z_M),$ and $w \in W(M),$ one may say ``$w\alpha >0$'' or ``$w\alpha <0$''
without ambiguity.  We say an element of $\Phi^+(Z_M)$ is indivisible if it is not of the
form $2\varepsilon_i.$

Each element $w\in W(M)$ can be decomposed as a product $s_{\alpha_1}
\dots s_{\alpha_\ell}$ of elementary symmetries as in \cite{MW1} I.1.8.
The element $s_{\alpha_\ell}$ will be in $W(M),$ while $s_{\alpha_{\ell-1}}$
will be in $W(s_{\alpha_\ell}Ms_{\alpha_\ell}^{-1})$ and so on.
Each is labeled with the unique indivisible restricted root (for the operative Levi) which
it reverses.  That is $\{ \alpha \in \Phi^+(Z_M): s_{\alpha_\ell} \alpha < 0\}
=\{\alpha_\ell\},$ or $\{\alpha_\ell,2\alpha_\ell\}$ and in the
latter case $\alpha_\ell = \varepsilon_r.$
(Cf. \cite{MW1} I.1.8.)

Let $w=s_{\alpha_1}\dots s_{\alpha_\ell}$ be a minimal-length decomposition
 into elementary symmetries, and put $w_i =s_{\alpha_{i+1}}\dots s_{\alpha_\ell}.$
Then
$$\{ \alpha \in \Phi^+(Z_M), \text{ indivisible }\;|\; w\alpha <0\} = \{ w_i^{-1} \alpha_i| \; 1\le i \le \ell\}$$
and $\ell$ is the cardinality of this set (i.e., there is no repetition).
Combining this discussion with that of the previous paragraphs, we obtain a
decomposition of $M(w, \underline s)$ as a composite of intertwining
operators $M_{w_i}(s_{\alpha_i}, w_i\cdot \underline s),$ each
corresponding naturally to one of the elements
of $\{ \alpha \in \Phi^+(Z_M), \text{ indivisible }\;|\; w\alpha <0 \}.$

Let $\det{}_i$ denote the rational character $(\underline g, \alpha) \mapsto \det g_i$
of $M.$  Then
$\{e_0,\det{}_1, \dots,\det{}_r \}$ is a basis for the lattice $X(M)$ of rational characters of $M.$
Here, the character $e_0$ of $T$ introduced in section \ref{rootDataSection}
has been identified with a character of $M$ as in \ref{r:identifications}.
Let $\{e_0^*,\det{}_1^*, \dots,\det{}_r^* \}$ be the dual basis of the
dual lattice. Again, $e_0^*$
is the same as in section \ref{rootDataSection}.
  Elements of $X(M)$ may be paired with elements
of $X^\vee(T)$ defining a projection from $X^\vee(T)$ onto the dual lattice.
For each $i=1, \dots m,$ there
exists unique $j\in \{ 1, \dots, r\}$ such that $e_i^*$ maps to $\det{}_j^*.$  If
$\alpha$ is a root, then the projection of the coroot $\alpha^\vee$ to the
dual lattice of $X(M)$ depends only on the restriction of $\alpha$ to $Z_M,$
and the correspondence is as follows:
$$0\leftrightarrow 0,$$
$$\varepsilon_i -\varepsilon_j \leftrightarrow \det{}_i^*-\det{}_j^*,$$
$$\varepsilon_i+\varepsilon_j \leftrightarrow \det{}_i^*+\det{}_j^*-e_0^*,$$
$$\varepsilon_i\leftrightarrow 2\det{}_i^*-e_0^*$$
$$2\varepsilon_i \leftrightarrow 2\det{}_i^*-e_0^*.$$
We denote the element corresponding to $\alpha \in \Phi^+(Z_M)$ by $\alpha^\vee$
(in agreement with \cite{MW1}, I.1.11).

We may identify $\underline s\in \C^r$ with
$$\sum_{i=1}^r \det{}_i \otimes s_i \in X(M)\otimes_\Z \C.$$
This is compatible with \cite{MW1}, I.1.4.  Restriction of functions gives a natural injective map
$X(M) \to X(T),$ and hence $X(M)\otimes_\Z \C\to X(T)\otimes_\Z\C,$ which we use to
identify the first space with a subspace of the second.  This gives the notation $w\cdot \underline s$
a precise meaning, as an element of $X(wMw^{-1}) \otimes_\Z\C,$ which is compatible
with the usage above.  In addition, it gives a ``meaning'' to the set
$$\{ s_i-s_j \}\cup \{s_i+s_j\} \cup \{2s_i\},$$
 of linear functionals
on $\C^r,$ identifying each with an element of $\Phi^+(Z_M).$  Formally,

\begin{defn}
A root hyperplane (relative to the Levi $M$)
is a hyperplane  of the form
$$H=\{ s \in \C^r\; | \; \langle \alpha^\vee, \underline s\rangle = c\}$$
for some $\alpha \in \Phi^+(Z_M)$ which is indivisible, and some $c \in \C.$
We say that the hyperplane
$H$ is associated to the root $\alpha,$ which is uniquely determined.
\end{defn}

 The next main statement is
 \begin{prop}\label{genIOtoRelRk1}
 Let $w=s_{\alpha_1} \dots s_{\alpha_\ell}$ be any decomposition
 of minimal length, and for each $i$ let $w_i=s_{\alpha_{i+1}} \dots s_{\alpha_\ell}.$
 Then the set of poles of $M(w,\underline s)$ is the disjoint union of the
 sets of poles of the operators $M_{w_i}(s_{\alpha_i},w_i \cdot \underline s).$
 A pole of $M(w,\underline s)$ comes from $M_{w_i}(s_{\alpha_i},w_i \cdot \underline s)$
 if and only if it is associated to $w_i^{-1}\alpha_i.$  Furthermore,
 if $\{ \underline s \in \C^r | \langle \alpha^\vee, \underline s \rangle = c\}$
 is a pole of $M(w, \underline s),$ then $c \ne 0.$
 \end{prop}

 This is proved in section \ref{genIOtoRelRk1 proof}.

 \section{Appendix IV: proof of Theorem \ref{t:EisensteinSeriesProperties}}
\label{section:  proof of Theorem thm:Eis}
\begin{enumerate}
\item
 We now prove Item \eqref{poleIsSimple}.  A root hyperplane
passing through
 $\poleloc$ is defined by an equation of one of three forms:
 $s_i = \frac 12, s_i+s_j=1,$ or $s_i-s_j=0.$  The third kind can not support
 singularities of the Eisenstein series.
 The first two can, but by \cite{MW1}IV.1.11 (c), they will be
 without multiplicity, and so the factor of
$$ \prod_{i\ne j}(s_i+s_j-1)
\prod_{i=1}^r (s_i-\frac 12)$$
will take care of them.

 The operators corresponding to elementary symmetries are called
 relative rank one because they could be defined without reference $G_{4n+1},$ considering
 $M$ instead as a maximal Levi of another Levi subgroup $M_\alpha$ of $G_{4n+1},$
 having semisimple rank one greater than that of $M.$   Furthermore, in a suitable
 sense, the relative rank one operator only ``lives on one component of $M_\alpha,$''
 which will allow us to deduce  the general case of {\bf \eqref{poleConditions}}
 from the case $r=1$ and a similar fact about intertwining operators on $GL_n.$
Let us make this more precise.

Fix $\alpha \in \Phi^+(Z_M).$  There is a minimal Levi subgroup $M_\alpha$
of $G_{4n+1}$
containing $M$ such that $\alpha$ is the restriction of a root of $M_\alpha.$
(It is standard iff $\alpha$ is the restriction of a simple root.)
  Fix $w\in W(M)$ such that $w\alpha<0,$
 and a decomposition
$w=s_{\alpha_1}\dots s_{\alpha_\ell}$ of $w$ as into elementary symmetries, which
is of minimal length.  For some unique $i,$ we have $\alpha = w_i^{-1}\alpha_i,$
where $w_i$ is as above.  Then $w_iM_\alpha w_i^{-1}$ is a standard Levi
of $G_{4n+1}.$
  Different choices of decomposition give different (even conjugate)
  embeddings of the same
reductive group into $G_{4n+1}$ as a standard Levi.

If $\alpha = \varepsilon_j-\varepsilon_k,$ or $\varepsilon_j+\varepsilon_k,$
then $M_{\alpha_i}$ is isomorphic to
$GL_{2(n_j+n_k)} \times \prod_{l\ne j,k} GL_{2n_l} \times GL_1.$
while if $\alpha = \varepsilon_j,$ it is isomorphic to
$G_{4n_j+1}\times \prod_{k\ne j}GL_{2n_k}.$  Let $G'$ denote
$GL_{2(n_j+n_k)}$ or $G_{4n_j+1}$ as appropriate and let $\iota$ be a choice
of isomorphism with the ``new'' factor.  Then $\iota^{-1}(\iota(G')\cap P^{w_i})$
is a maximal parabolic subgroup $P'=M'U'$ of $G',$ and
$\sigma:=(\bigotimes_{i=1}^r \tau \otimes \omega)\circ Ad(w_i)\circ \iota,$
is an  irreducible unitary
cuspidal automorphic representation of $M'(\A).$
The map $\iota$ also induces a linear projection
$$\iota_*:X(w_iMw_i^{-1})\otimes_\Z\C \to X(M')\otimes_\Z\C.$$
(Recall that we have agreed to think of $w_i\cdot \underline s$
as an element of the former space.)

Following, \cite{MW1} I.1.4, define $m^\mu$
for $m\in M'(\A)$ and $\mu$ in $X(M')\otimes_\Z\C,$
 by stipulating that
$m^\mu=|\chi(m)|^s$ if $\mu=\chi\otimes s$ and $m^{\mu_1+\mu_2}=m^{\mu_1}m^{\mu_2}.$

The set $W_{G'}(M'),$ defined
analogously to $W(M)$ above,
contains a unique nontrivial element.  It is the elementary symmetry $s_\beta$
associated to the restriction to $Z_{M'}$ of the unique simple root  of $G'$ which
is not a root of $M'.$
The map $\iota$ identifies $s_\beta$ with $s_{\alpha_i}.$

For $\mu \in X(M')\otimes_\Z\C,$
let   $V^{(1)}(\mu, \sigma)$ denote
$$\{h:G'(\A) \to V_\sigma, \text{ smooth}\; |\; h(mg')(m')=h'(g')(m'm)
m^{\mu+\rho_{P'}} \quad m,m'\in M'(\A), g'\in G'(\A)\},$$
$$V^{(2)}(\mu, \sigma)=
\{h:G'(\A) \to \C, \text{ smooth} | h(g')(e)
\in V^{(1)}(\mu, \sigma)
\}.$$

 There is a standard intertwining operator
$M(s_\beta,\mu): V^{(2)}(\mu, \sigma)\to V^{(2)}_{s_\beta}(\mu, \sigma).$
One has the identity
$$M_{w_{i-1}}(s_{\alpha_i}, w_i \cdot \underline s)
f(\iota(h)g) =
M(s_\beta, \mu)f(\iota(h)g).$$
That is, if $p_g$ denotes the map
$$V^{(2)}_{w_i}(\underline s, \bigotimes_{i=1}^r  \tau_i \boxtimes \omega)
\to V^{(2)}(\mu, \sigma)$$
corresponding to evaluation at $\iota(h)g$ for a fixed $g,$ then, for all $g,$
the following diagram commutes:
  $$
  \begin{CD}
  V^{(2)}_{w_i}(\underline s, \bigotimes_{i=1}^r  \tau_i \boxtimes \omega)
  @>M_{w_i}(s_{\alpha_i},w_i\cdot\underline s)>>
   V^{(2)}_{w_{i-1}}(\underline s, \bigotimes_{i=1}^r  \tau_i \boxtimes \omega)
      \\
@V{p_g}VV @V{p_g}VV\\
 V^{(2)}(\iota_*(w_i\cdot \underline s+\rho_{P_{\alpha_i}}), \sigma)@>{M(s_\beta,
 \iota_*(w_i\cdot \underline s+\rho_{P_{\alpha_i}}))}>>
  V^{(2)}_{s_\beta}(\iota_*(w_i\cdot \underline s+\rho_{P_{\alpha_i}}), \sigma).
  \end{CD}
  $$
  Hence $M_{w_i}(s_{\alpha_i},w_i\cdot\underline s)$ has a pole along a root hyperplane
  associated to $\alpha$ iff $M(\iota_*(w_i \cdot s +\rho_{P_{\alpha_i}}),\sigma)$
  does.

Since the set of poles of $M_{w_i}(s_{\alpha_i},w_i\cdot \underline s)$
is equal to the set of poles of $M(w,\underline s)$ along hyperplanes associated
to $\alpha,$ it is independent of the choice of decomposition
$w=s_{\alpha_1} \dots s_{\alpha_\ell}.$
Hence, for each
$\alpha \in \Phi^+(Z_M),$ we may use a decomposition  tailored to that $\alpha.$

First suppose $\alpha = \varepsilon_j-\varepsilon_k.$  One may choose a decomposition
so that $w_i$ corresponds to the permutation matrix in $GL_{2n}$ (identified with a subgroup
of the Siegel Levi) which moves the $j$th block of $M$ up so that it is immediately after the
$i$th, and otherwise preserves order.  It is then easily verified that
$\sigma = \tau_i\otimes\tau_j$ and
$$\begin{pmatrix}h_1&\\ &h_2\end{pmatrix}^{ \iota_*( w_i\cdot s )+\rho_{P_{\alpha_i}}}
= |\det h_1|^{s_i+\kappa}|\det h_2|^{s_j+\kappa},$$
where $\kappa = \sum_{k>i,k\ne j}n_k- \sum_{k<i}n_k+n.$

Next suppose $\alpha=2\varepsilon_j.$  Then we choose a decomposition so that
$w_i$ is in the Weyl group of $GL_{2n},$ and moves the $j$th block to be last,
otherwise preserving order.
Then one easily verifies that $\sigma$ is the representation $\tau_j \boxtimes \omega$
of the Siegel Levi of $G_{4n_j},$ and
that, for $(g',\alpha)$ in the Siegel Levi of $G_{4n_j},$
$$(g',\alpha)^{\iota_*(w_i\cdot \underline s+\rho_{P_{\alpha_i}})}= |\det g'|^{s_j}.$$

Finally, suppose $\alpha= \varepsilon_j+\varepsilon_k.$  Then we choose a
decomposition so that $w_i$ that projects to a permutation matrix in $SO_{4n+1}$ of the form
$$\begin{pmatrix}
I&&&&\\
&&&I&\\
&&I&&\\
&I&&&\\
&&&&1
\end{pmatrix},$$
with the off-diagonal blocks being $2n_j\times 2n_j,$ and the first block being
$\sum_{k=1}^i 2n_k.$
We deduce  from Corollary \ref{c:weylAction} that $\sigma = \tau_i \otimes (\tilde\tau_j\otimes \omega),$
and from Lemma \ref{l:WeylActionii} that
$$\begin{pmatrix}h_1&\\ &h_2\end{pmatrix}^{ \iota_*( w_i\cdot s )+\rho_{P_{\alpha_i}}}
= |\det h_1|^{s_i+\kappa}|\det h_2|^{-s_j+\kappa},$$
where $\kappa$ is as before.

\item

Item \eqref{poleConditions} follows from
\begin{prop}\label{relRk1toLFcn-GSpin}
 Let $w$ denote the unique nontrivial element of $W(M),$ in the
 case when $M$ is the Levi of the Siegel parabolic of $G_{2m+1}.$
 Let  $\tau$ be a cuspidal
 representation of $GL_{m}.$
  Then
 $M(w, s)f(g)$ has a pole at $s=\frac 12$
 for some
 $f \in
 \Ind_{P(\A)}^{G_{2m+1}(\A)} (\tau \otimes |\det|^s)\boxtimes \omega,$
and $g \in G_{2m+1}(\A)$
  if and only if
  $\tau$ is $\omega^{-1}$-orthogonal
 \end{prop}
 \begin{rmk}\label{blindtoparity}
 Of course we are only interested in the case $m=2n.$  Furthermore, since we
 assume $\omega$ is not the square of another Hecke character, it follows that
 $\tau$ can be $\omega^{-1}$-orthogonal only if $m$ is even.  However,  the proof
 of this proposition is ``blind to'' the parity of $m.$
 \end{rmk}
\begin{prop}\label{relRk1toLFcn-GL}
Let $P=MU$ be a maximal standard parabolic of $GL_n$ such that $M\iso GL_k \times GL_{n-k}.$
Let $f$ be an element of $\Ind_{P(\A)}^{GL_n(\A)}
(\tau_1 \otimes |\det |^{s_1} ) \bigotimes (\tau_2 \otimes |\det|^{s_2}).$
Let $w$ be the unique nontrivial element of $W(M).$  Then
$M(w, \underline s)f(g)$ is singular along the hyperplane $s_1-s_2=1$ for some
$f,g$ iff $n=2k$ and $\tau_2 \iso \tau_1.$
\end{prop}
We defer the proofs to the section \ref{s:appIII}.

Now, we assume \eqref{e:DistinctAndOmegaOrthogonal}
 holds and prove the remaining part of the theorem.
 Let
$N(\underline s)= \prod_{i=1}^r (s_i-\frac 12).$

\item Item \eqref{ResidueIsL2} follows from \cite{MW1} I.4.11.  The constant
term of $E(f)$ along a parabolic $P'=M'U'$ has nontrivial cuspidal
component iff $M'$ is conjugate to $M.$ (\cite{MW1} IV.1.9 (b)(ii)).
For such $P'$ it is equal to
$$\sum_{w\in W(M),\; wMw^{-1}=M'} M(w,\underline s)f(g).$$
Take $w\in W(M),$ such that  $wMw^{-1}=M'.$  If $w\cdot\varepsilon_i>0$
for some $i,$ then
$M(w,\underline s)f(g)$ does not have a pole at $s_i-\frac 12,$
and hence
$N(\underline s)M(w,\underline s)f(g)$
vanishes at $\poleloc.$  On the other hand, if
$w\cdot\varepsilon_i<0$ for all $i,$  then $M(w,\underline s)f(g)$ satisfies the criterion
of \cite{MW1} I.4.11.

 It  follows from \cite{MW1} IV.1.9 (b)(i) applied to $N(\underline s)E(f)$ (which is valid
by \cite{MW1} IV.1.9 (d)) that the residue is an automorphic form.

\item
To complete the proof of Item \eqref{ResidueMapIsAnIntOp}, let $\rho(g)$ denote
right translation.  It is clear that for values of $s$ in the domain of convergence,
$N(\underline s)E(\rho(g) f)(\underline s) = N(\underline s) \rho(g)(E(f)(\underline s)).$
By uniqueness of analytic
continuation, the equality also holds at values of $s$ where both sides are defined
by analytic continuation, including $\poleloc$.  The action of the Lie algebra at
the infinite places is handled similarly.

Next we consider the constant term of $E(f)$ along the Siegel parabolic.  By
\cite{MW1} II.1.7(ii) it may be expressed in terms of $GL_{2n}$ Eisenstein
series, formed using the functions $M(w,\underline s)f,$ corresponding
to those $w\in W(M)$ such that $w^{-1}(e_i-e_{i+1})>0$ for all $i.$
(Note: we proved in Lemma \ref{l:WeylAction} that $wMw^{-1}$ is contained
in the Siegel Levi for every $w\in W(M).$)  When we pass to $E_{-1}(f),$
the term corresponding to $w$ only survives if $w\cdot \varepsilon_i<0$
for all $i.$  This condition picks out a unique element, $w_0.$  It is the shortest element of
$W_{GL_{2n}}\cdot w_\ell \cdot W_{GL_{2n}},$ where $w_\ell$ is the longest
element of $W_{G_{4n+1}},$ and we have identified $GL_{2n}$ with a subgroup
of the Siegel Levi as usual.  Via corollary \ref{c:weylAction} one finds
that
$$(\bigotimes_{i=1}^r\tau_i \boxtimes \omega)\circ Ad(w_0)
= (\bigotimes_{i=1}^r (\tilde \tau_{r+1-i}\otimes \omega)\boxtimes \omega)
=(\bigotimes_{i=1}^r \tau_{r+1-i}\boxtimes \omega).$$
For $f\in V^{(2)}(\bigotimes_{i=1}^r \tau_i \boxtimes\omega, \poleloc),$
$M(w_0,\poleloc)f|_{GL_{2n}(\A)}$ is an element of the analogue of
$V^{(2)}(\bigotimes_{i=1}^r \tau_i \boxtimes\omega, \underline s),$
for the induced representation
$$\Ind_{\bar P^0(\A)}^{GL_{2n}(\A)} (\bigotimes \tau_{r+1-i}\otimes |\det{}_i|^{n-\frac 12})
=|\det|^{n-\frac 12} \otimes\tau
$$
of $GL_{2n}.$
Here $\bar P^0=GL_{2n} \cap P^{w_0},$ and $\tau=\tau_1\boxplus\dots\boxplus\tau_r.$
Furthermore, since this representation is irreducible, it may be regarded as
an arbitrary element.   Also, we may regard this representation as
induced from  $\tau_1, \dots, \tau_r$ in the usual order.  Let $\bar P$ denote
the relevant parabolic of $GL_{2n}.$

The representation $\tau$ sits inside a fiber bundle of induced representations
$\Ind_{\bar P(\A)}^{GL_{2n}(\A)} (\bigotimes_{i=1}^r \tau_i\otimes |\det{}_i|^{s_i}).$
For a flat, $K$-finite section $f$ let
 $E^{GL_{2n}}(f)(g)(\underline s)$ be the $GL_{2n}$ Eisenstein series
defined by
$$\sum_{\bar P(F)\backslash GL_{2n}(F)}
f(\underline s)(\gamma g)$$
when  $s_i-s_{i+1}$ is sufficiently large for each $i,$ and by meromorphic continuation
elsewhere.

Let $\maxunip^{GL_{2n}}$ denote the usual maximal unipotent
subgroup of $GL_{2n},$ consisting of all upper triangular unipotent matrices.  Let
$\psi_W(u)=\baseaddchar(u_{1,2}+\dots+u_{m-1,m})$ be the usual
generic character.

\item To complete the proof of Item \eqref{residueSupportsPeriod}, we must prove that
\begin{equation}\int_{\maxunip^{GL_{2n}}\quo}\label{e:GLmEisSerGeneric}
E^{GL_{2n}}(f)(ug)(\underline 0)\psi_W(u) \; du \ne 0\end{equation}
for some $f\in \Ind_{\bar P(\A)}^{GL_{2n}(\A)}
\bigotimes_{i=1}^r \tau_{r+1-i}, g \in GL_{2n}(\A),$
i.e., that the space of $GL_{2n}$ Eisenstein series   $E^{GL_{2n}}(f)$ is
globally $\psi_W$-generic.
Granted this, {\bf \eqref{residueSupportsPeriod}} follows from \cite{MW1}II.1.7(ii)
and the discussion just above.

The following proposition follows from work of Shahidi.
\begin{prop}\label{FourierCoeffOfEisSer}
We have the following:
\begin{enumerate}
\item
$$\int_{\maxunip^{GL_{2n}}\quo}
E^{GL_{2n}}(f)(ug)(\underline s)\psi_W(u) \; du \\
=\prod_{v\in S}W_v(g_v) \cdot \prod_{v \notin S} W_v^\circ
(g_v)
\cdot \prod_{i<j}L^S(s_i-s_j+1,\tau_i\times \tilde \tau_j)^{-1},
$$
where,
\begin{itemize}\item
for each $v,$ $W_v$ is a Whittaker function in the $\psi_{W,v}$-Whittaker
model of $\Ind_{\bar P(F_v)}^{GL_{2n}(F_v)}(\bigotimes_{i=1}^r\tau_{i,v} \otimes |\det{}_i|_v^{s_i}),$
\item
$S$ is a finite set of places, depending on $f,$
\item for $v \notin S$, $\tau_v$ is unramified
\item for $v \notin S$,  $W_v^\circ$ is
the normalized spherical vector in the $\psi_{W,v}$-Whittaker
model of $\Ind_{\bar
P(F_v)}^{GL_{2n}(F_v)}(\bigotimes_{i=1}^r\tau_{i,v} \otimes
|\det{}_i|_v^{s_i}).$

\end{itemize}
\item A flat, $K$-finite section $f$ may be chosen so that, for all $v\in S,$ the function $W_v$
is not identically zero at $\underline s=\underline 0.$
\end{enumerate}
\end{prop}
We briefly review the steps of the proof in  section \ref{s:InducedRepFactors}.

It follows from \cite{JacquetShalika-EPandClassnII}
Propositions 3.3 and 3.6 that the product of partial $L$ functions appearing
in Proposition \ref{FourierCoeffOfEisSer}
does not have a pole at $\underline s=\underline 0$ provided the representations $\tau_1, \dots, \tau_r$
are distinct.  This completes the proof of {\bf \eqref{residueSupportsPeriod}}.

\item Finally, Item \eqref{indepOfOrder} follows from the functional equation of the Eisenstein series
(\cite{MW1}IV.1.10(a)), and the
fact that $\tau$ is equal to an irreducible full induced representation (as opposed
to a constituent of a reducible one).
\end{enumerate}

\section{Appendix V: Auxilliary results used to prove
Theorem \ref{t:EisensteinSeriesProperties}
}\label{s:appIII}
In this appendix we complete the proofs of several intermediate statements used in
the proof of Theorem \ref{t:EisensteinSeriesProperties}.  As far as we know, all of these results are well-known to the
experts, but do not appear in the literature in the precise form we need.
\subsection{Proof of Proposition \ref{p:EStoIO}}
\label{ss: proof of ES to IO}
First, suppose that a set $D$ of hyperplanes carries all the singularities of
 all the intertwining operators $M(w, \underline s)f.$   Then it follows from
 \cite{MW1} II.1.7, IV.1.9 (b) that it carries all the singularities of the cuspidal
 components of all the constant terms of $E(f)(g)(\underline s).$  By I.4.10, it
 therefore carries the singularities of the Eisenstein series itself.

On the other hand, it is clear that a set which carries the singularities of the
Eisenstein series carries those of all of its constant terms.  Thus, what we
need to prove is:
\begin{lem}  Fix  $M'$ a standard Levi which is conjugate to
$M$ and $\alpha \in \Phi^+(Z_M).$
Let $H$ be the root hyperplane given by $\langle \alpha^\vee , \underline s\rangle =c,$
$c\ne 0.$
Consider the family of
 functions
$M(w,\underline s)f$ corresponding to
$\{w\in W(M)| wMw^{-1} = M'\}.$
If any one or them has a pole along $H,$ then the constant term of the Eisenstein series along $P'$
does as well.  In other words, it is not possible for two poles to cancel one another.
\end{lem}
\begin{proof}
Clearly, it is enough to prove this under the additional hypothesis that $M'=M.$

Let $A_M^+$ denote the group isomorphic to  $(\R^\times_+)^{r+1},$
embedded diagonally at the infinite places,  which is inside the center of $M.$

The Lie algebra of $A_M^+$ is naturally identified with the real dual of $X(M)\otimes_\Z \R.$
Recall that above we identified $\underline s$ with an element of
$X(M)\otimes_\Z \C.$
So, there is a natural pairing $\langle X, \underline s\rangle,$ $X \in \mathfrak a_M^+,$
given as follows.  Write $\det_i$ for the determinant of the $i$th block of
an element of $M,$ regarded as a $2n\times 2n$ matrix via the identification
with $GL_m \times GL_1$ fixed above.  Then we have
 $$\prod_{i=1}^r|\det{}_i \exp(\log y \cdot X)|^{s_i}= y^{\langle X,\underline s\rangle}.$$
It follows that
$$|M(w,\underline s)f(\exp( \log y \cdot X) g)| = y^{\Re(\langle w^{-1} X, \underline s\rangle)}
\cdot \delta_P^{\frac 12}(w^{-1}\exp( \log y\cdot X)w) \cdot
|M(w,\underline s)f(g)|.$$
Here  $\delta_P$ is the modular quasicharacter of $P.$

Let
$$W_{sing}(M,H)= \{ w\in W(M), wMw^{-1} =M,  M(w,\underline s) \text{ has a pole along }H\}.$$
Suppose that this set is nonzero.  Choose $w_0\in W_{sing}(M,H)$ such that the
order of the pole of $M(w_0,\underline s)$ is of maximal order.
Let $\nu(H)$ denote the order.
Choose $X \in \mathfrak a_M^+$ such that the points $w^{-1}\cdot X, w \in W_{sing}(M,H)$
are all distinct.
Consider the family of functions
$$(\langle \alpha^\vee,\underline s\rangle -c)^{\nu(H)} M(w,\underline s)f(\exp(\log y\cdot X) g),
\quad
w \in W_{sing}(M,H).$$
They have singularities carried by a locally finite set of root hyperplanes not containing $H.$
Assume $g$ has been chosen so that
$(\langle \alpha^\vee,\underline s\rangle -c)^{\nu(H)} M(w_0,\underline s)f(g) \ne 0.$
For $\underline s$ restricted to an open subset of $H$ not intersecting any of the
singular hyperplanes we obtain a family of holomorphic functions, at least
one of which is nonzero.
 If we further exclude the intersection of $H$ with the hyperplanes
$$\langle w_1^{-1} X-w_2^{-1}X,\underline s\rangle = 0, \quad w_1,w_2\in  W_{sing}(M,H),$$
(which can not  coincide with $H$ because $c\ne 0$),
then at every point $\underline s,$ those  functions which are nonzero all have distinct orders of magnitude
as functions of $y.$  Hence they can not possibly cancel one another.
  \end{proof}

\subsection{Proof of Lemma \ref{genIOtoRelRk1}}\label{genIOtoRelRk1 proof}
Regarding $w_i\cdot \underline s+\rho_{P_{\alpha_i}}$
as an element of $X(w_iMw_i^{-1})\otimes_\Z\C,$
we may decompose it as $\mu_1+ \langle \alpha_i^\vee, w_i\cdot \underline s\rangle\tilde \alpha_i,$
where $\tilde \alpha_i$ is defined by the property that
$$\langle \alpha^\vee ,\tilde\alpha_i\rangle =\delta_{\alpha, \alpha_i},\quad
\text{for }\alpha \in \Phi^+(Z_{w_iMw_i^{-1}}).$$
Then it follows easily from the definitions that $\mu_1$ is in the image of the natural projection
$X(M_{\alpha_i})\otimes_\Z\C\to X(w_iMw_i^{-1})\otimes_\Z\C$ corresponding to
restriction of characters of $M_{\alpha_i}(\A)$ to $w_iMw_i^{-1}(\A).$

Take $f$ a $K$-finite flat section of $\Ind_{P^{w_i}(\A)}^{G_{4n+1}(\A)}
(\bigotimes_{j=1}^r \tau_j\otimes |\det{}_j|^{s_j}\boxtimes\omega)\circ Ad(w_i^{-1}).$
Then $M_{w_i}(s_{\alpha_i},w_i\cdot \underline s)f$ resides in a finite
dimensional subspace of $\Ind_{P^{w_{i-1}}(\A)}^{G_{4n+1}(\A)}
(\bigotimes_{j=1}^r \tau_j\otimes |\det{}_j|^{s_j}\boxtimes\omega)\circ Ad(w_{i-1}^{-1}),$
corresponding to a finite set of $K$-types determined by $f.$  Write
$M_{w_i}(s_{\alpha_i},w_i\cdot \underline s)f$ in terms of a basis of flat $K$-finite sections.
The coefficients are functions of $\underline s,$ but it follows easily from the
integral definition where this is valid, and by meromorphic continuation elsewhere, that
in fact they are independent of $\mu_1$ (which corresponds to a character of
$M_{\alpha_1}(\A)$ and may be pulled out of the integration).  Thus,
they depend only on $\langle w_i \cdot \underline s , \alpha_i^\vee\rangle=
\langle \underline s, w_i^{-1}\alpha_i^\vee\rangle.$

The first two assertions are now clear.  A proof that $c\ne 0$ is obtained by a straightforward
modification of the opening paragraph of \cite{MW1}, IV.3.12.

\subsection{Proof of Proposition \ref{relRk1toLFcn-GSpin}}
In this section, we denote by $V^{(i)}(s,\tau,\omega), \; i=1,2,$ the spaces of functions
previously introduced in section \ref{s:InducedReps} as
$V^{(i)}(\underline s, \bigotimes_{i=1}^r \tau_i \boxtimes \omega),$
in the special case when $r=1.$

Let $\tilde M(s)$ denote the analogue of $M(w,s)$ defined using
$V^{(1)}( s , \tau , \omega ).$  It maps into the space
$V^{(3)}( -s , \tilde\tau\otimes\omega , \omega )$ given by
$$\left\{ \tilde F: \ourgroup_{2m+1}(\A) \to V_{\tau}, \text{ smooth } \left|
\tilde F( (g, \alpha) h )(g_1)  = \omega(\alpha \det g )|\det g|^{-s+\frac m2}
\tilde F(h)(g_1 \;_tg^{-1}) \right.\right\}.$$
Fix realizations of the local induced representations $\tau_v$ and
an isomorphism
$\iota: \otimes_v'\tau_v \to \tau.$
Define, for each $v,$
$V^{(1)}(s, \tau_v, \omega _v)$ to be
$$\left\{ \tilde F_v: \ourgroup_{2m+1}(F_v) \to V_{\tau_v}, \text{ smooth } \left|
\tilde F_v( (g, \alpha) h ) = \omega_v (\alpha)|\det g|_v^{s+\frac m2}
\tau_v(g) \tilde F_v(h) \right.\right\},$$ and
$V^{(3)}(s, \tilde\tau_v\otimes\omega_v, \omega _v)$ to be
$$\left\{ \tilde F_v: \ourgroup_{2m+1}(F_v) \to V_{\tau_v}, \text{ smooth } \left|
\tilde F_v( (g, \alpha) h ) = \omega_v (\alpha\det g)|\det g|_v^{s+\frac m2}
 \tau_v(_tg^{-1}) \tilde F_v(h) \right.\right\}.$$

Then the formula
$$\tilde \iota (\otimes_v \tilde F_v )( g) =  \iota( \otimes_v ' \tilde F_v (g_v) )$$
defines   maps
$$\otimes_v' V^{(1)}(s, \tau_v, \omega _v)
\to V^{(1)}( s , \tau , \omega ),
$$
$$\otimes_v' V^{(3)}(s, \tilde\tau_v\otimes\omega_v, \omega _v)
\to V^{(3)}( s , \tilde\tau\otimes\omega_v , \omega ),
$$
both of which we denote  by $\tilde \iota.$

It is known \label{s:InducedRepFactors}
that each map is, in fact, an isomorphism.  For the benefit of the reader we
sketch an argument.  On pp. 307 of \cite{Shahidi-ParkCity}
certain explicit elements of (a generalization of) $V^{(1)}(s,\tau,\omega )$
are constructed as integrals involving matrix coefficients.  Using Schur orthogonality,
one may check
that $\tilde F$ is expressible in this form iff both the $K$-module it
generates and the $K \cap M(\A)$-module it generates are irreducible.
It is clear that such vectors span the space of all $K$-finite vectors.
On the other hand the (finite dimensional)
space of matrix coefficients of this irreducible representation of $K$
is spanned by those that factor as a product of matrix coefficients of
local representations, and these are clearly in the image of $\tilde\iota.$

For $\tilde F_v \in V^{(1)}(s, \tau_v, \omega _v),$ let
$$A_v( s) \tilde F_v ( g)
= \int_{U_w(F_v) } \tilde F_v ( \dot w u g ) du.$$
  Then the following diagram commutes
  $$
  \begin{CD}
   \otimes_v' V^{(1)}(s, \tau_v, \omega _v)@>A(s)>>\otimes_v' V^{(1)}(-s, \tau_v, \omega _v)\\
@V{\tilde \iota}VV @V{\tilde \iota}VV\\
 V^{(1)}( s , \tau , \omega )
@>{\tilde M(s)}>>
 V^{(1)}( -s , \tau , \omega )
  \end{CD}
  $$
with  $A(s) := \otimes_v A_v(s).$

Now, $M(w,s)f(s)$ has a pole (i.e., there exists $g\in G_{2m+1}(\A)$ such that $M(w,s)f(s)(g)$ has a pole)
if and only if $\tilde M(s) \tilde F(s)$ has a pole (i.e., there exist $g\in G_{2m+1}(\A)$ and
$m\in M(\A)$ such that
$\tilde M(s) \tilde F(s) (g)(m)$ has a pole), where $\tilde F$ is the element of
$V^{(1)}( s , \tau , \omega )$ such that $f(g)=\tilde F(g)(id).$

We wish to show that there exists $\tilde F$ such that this is the case iff $\tau$ is $\omega^{-1}$-orthogonal.
Clearly, we may restrict attention to $\tilde F$ of the form $\tilde \iota( \otimes_v \tilde F_v).$

Recall that for all but finitely many non-archimedean $v,$ the space $V_{\tau_v}$ comes equipped with a choice of
$GL_{m}( \mathfrak{o}_v)$-fixed vector $\xi_v^\circ$ used to define the restricted tensor
product.

If  $\tilde F= \tilde \iota( \otimes_v \tilde F_v)
\in V^{(1)}( s , \tau , \omega ),$
then there is a finite set
$S$ of places, such that if $v\notin S$
then $v$ is non-archimedean,
$\tau_v$ is unramified,
and $\tilde F_v (s)= \tilde F_{(s, \tau_v,\omega _v)}^\circ$ is the unique element of
$V^{(1)}(s, \tau_v, \omega _v)$ satisfying $\tilde F_{(s, \tau_v,\omega _v)}( k) = \xi_v^\circ$ for all $k \in
G_{2m+1}(\mathfrak{o}_v).$

Now
$$A_v(s) \tilde F_{(s, \tau_v,\omega _v)}^\circ
= \frac{ L_v( 2s, \tau_v, sym^2 \times \omega^{-1} _v)}
{L_v( 2s+1 , \tau_v, sym^2 \times \omega^{-1} _v)}
 \tilde F_{(-s, \tilde\tau _v\otimes\omega_v,\omega_v)}^\circ
.$$
(A proof of this appears in \cite{eulerproducts}, albeit not in this precise language.
See especially pp. 25-27.)
Thus,
$$A(s) \tilde \iota( \otimes_v \tilde F_v)
=
\frac{ L^S( 2s, \tau, sym^2 \times \omega^{-1} )}
{L^S( 2s+1 , \tau, sym^2 \times \omega^{-1} )}
\tilde \iota\left(
\left(
\bigotimes_{v \in S}A_v(s) \tilde F_v(s)\right) \otimes
\left( \bigotimes_{v \notin S} \tilde F_{-s , \tilde\tau _v\otimes\omega_v , \omega _v} \right) \right).$$
To complete the proof we must show:
\begin{description}
\item[(i)]{\label{locintop} $A_v(s)$ is holomorphic and nonvanishing (i.e., not the zero operator)
on $Ind_{P(\A)}^{G_{2m}(\A)}\tau \otimes |\det |^s \boxtimes \omega$
at $s= \frac 12,$ for all $\tau.$  }
\item[(ii)]{ \label{loclLfunction}
$L_v( s, \tau_v, sym^2 \times \omega^{-1} _v)$ is holomorphic and nonvanishing
at $s=1,$ for all $\tau_v.$  }
\item[(iii)]{\label{partialLfunction} $L^S( s, \tau , sym^2 \times \omega^{-1} )$ is holomorphic and nonvanishing at $s=2.$}
\end{description}

Item (iii) is  covered by Proposition
7.3 of \cite{Kim-Shahidi-simplicity}.
Items (i) and (ii) are  essentially contained in Proposition 3.6, p. 153 of
\cite{Asg-Sha1}.   Since what we need is {\sl part} of the same information,
presented differently, we repeat the part of the arguments we are using.

The nonvanishing part of (i) is a completely general fact (i.e., holds
at least for any Levi of any split reductive group).
For example, the only element of the arguments made on p. 813 of \cite{GRS3} which
is particular to the situation they consider there (the Siegel of $Sp_{4n}$) is the precise
ratio of $L$ functions appearing in the constant term.

Similarly, local $L$ functions never vanish.  At a finite prime the local $L$ function is
$P(q_v^{-s})^{-1}$ for some polynomial $P,$ while at an infinite prime it is given in
terms of the $\Gamma$ function and functions of exponential type.

We turn to holomorphicity.
\begin{lem}\label{l:tadic}
Let $\pi$ be any representation of $GL_m(F_v),$ which is irreducible, generic, and unitary.
Then there exist
\begin{itemize}
\item
integers $k_1, \dots, k_r$ of such that $k_1+ \dots +k_r=m,$
\item real numbers $\alpha_1, \dots , \alpha_r \in ( -\frac 12 , \frac 12),$
\item discrete series representations $\delta_i$ of $GL_{k_i}(F_v)$ for $i=1$ to $r$
\end{itemize}
such that
$$\pi\cong
Ind_{P_{(k)}(F_v)}^{GL_m(F_v)}
\bigotimes_{i=1}^r
 ( \delta_i \otimes |\det{}_i |^{\alpha_i}).$$
 Here $P_{(k)}$ denotes the standard parabolic of $GL_m$ with Levi consisting of block diagonal
 matrices with the block sizes $k_1, \dots ,k_r$ (in that order),
  and  $\det_i$  denotes the determinant of
  the $i$th block.
\end{lem}
\begin{rmk}In fact, one may prove a much more precise statement, but the above is
what is needed for our purposes.
\end{rmk}
\begin{proof}
This follows from \cite{JacquetShalika-The Whittaker models of induced representations}
and proposition 3.3
of \cite{MullerSpeh}.
\end{proof}

Continuing with the proof of Proposition \ref{relRk1toLFcn-GSpin}, let $(k)= (k_1, \dots, k_r),$
$\delta = (\delta_1, \dots, \delta_r)$ and $\alpha = (\alpha_1 , \dots , \alpha_r)$
be
such that
$$\tau_v\cong
Ind_{P_{(k)}(F_v)}^{GL_m(F_v)}
\bigotimes_{i=1}^r
 ( \delta_i \otimes |\det{}_i |^{\alpha_i}),$$
 and let
$\tilde P_{(k)}$ denote the standard parabolic of $G_{2m}$ which is contained in
the Siegel parabolic $P$
such that $\tilde P_{(k)}\cap M = P_{(k)}.$
Then
$$Ind_{P(F_v)}^{G_{2m}(F_v)} \tau_v \otimes |\det|_v^s \boxtimes \omega_v^s
\iso
Ind_{\tilde P_{(k)}(F_v)}^{G_{2m}(F_v)}
\bigotimes_{i=1}^r
 ( \delta_i \otimes |\det{}_i |_v^{s+\alpha_i}) \boxtimes \omega_v.$$
 This family (as $s$ varies) of representations lies inside the
 larger family,
 $$Ind_{\tilde P_{(k)}(F_v)}^{G_{2m}(F_v)}
\bigotimes_{i=1}^r
 ( \delta_i \otimes |\det{}_i |^{s_i}) \boxtimes \omega_v \qquad
 s = (s_1, \dots, s_r) \in \C^r,$$
 and our intertwining operator $A_v(s)$ is the restriction, to the line
 $s_i = s + \alpha_i$ of the standard intertwining operator for this induced
 representation, which we denote $A_v( \underline{s}).$  This operator
 is defined, for all $\Re(s_i)$ sufficiently large, by the same integral
 as $A_v(s).$

A result of Harish-Chandra says
 that ``$\Re(s_i)$ sufficiently large'' can
 be sharpened to ``$\Re(s_i) > 0.$''  (This is because all $\delta_i$  are
 discrete series, although tempered would be enough.)
 This result is given in the $p$-adic case as   \cite{Silberger} Theorem 5.3.5.4,
and in the Archimedean case, \cite{Knapp} Theorem 7.22, p. 196.

 Hence, the integral defining  $A_v(s)$ converges for $s > \max_i(- \alpha_i),$
 and in particular converges at $\frac 12.$

From the relationship between the local $L$ functions and the so-called
local coefficients, it follows that the local $L$ functions
are also holomorphic in the same region.  For this relationship see
\cite{Shahidi-Artin} for the Archimedean case and
\cite{Shahidi-Plancharel}, p. 289 and p. 308 for the non-Archimedean case.

This completes the proof of (i) and (ii).

\subsection{Proof of Proposition \ref{relRk1toLFcn-GL}}
The proof is the same as the previous proposition, except that the ratio of partial $L$ function
which emerges from the intertwining operators at the unramified
places is $$\frac{L^S(s_1-s_2,\tau_1\times \tilde\tau_2)}{L^S(s_1-s_2+1,\tau_1\times\tilde \tau_2)}.$$
  Convergence of local $L$ functions and intertwining operators at $s_1-s_2=1$ follows again
  from Lemma \ref{l:tadic}.  The only difference is the reference for (iii), which in this case is
 Theorem   5.3 on p. 555 of  \cite{JS2}.

  \subsection{Proof of \ref{FourierCoeffOfEisSer}}
As noted, this material is mostly due to Shahidi.

Since the statement is true (with the same proof) for general $m,$ not only $m=2n,$
we prove it in that setting.

In this subsection only, we write $\tau$ for the irreducible
unitary cuspidal representation $\bigotimes_{i=1}^r\tau_i$ of $M(\A)$
(as opposed to the isobaric representation $\tau_1\boxplus\dots\boxplus\tau_r$).

First, observe that the integral in question is clearly absolutely and uniformly convergent,
and as such defines a meromorphic function of $\underline s$ for each $g$ with poles
contained in the set of poles of the Eisenstein series itself.

For $\underline s$ in the domain of convergence
\begin{equation}\label{step1}
\int_{\maxunip^{GL_m}\quo}
E^{GL_m}(f)(ug)(\underline s)\psi_W(u) \; du
=\int_{U_{w_1}(\A)\cdot U^{w_1}\quo}
f(\underline s)(w_1^{-1}ug)\psi_W(u) \; du,
\end{equation}
where $w_1$ is the longest element of $W_{GL_m}(\bar M)$  (defined analogously
to $W(M)$ above),
$U_{w_1}=\maxunip^{GL_m} \cap w_1\overline\maxunip^{GL_m}w_1^{-1}$
and $U^{w_1}=\maxunip^{GL_m} \cap w_1\maxunip^{GL_m}w_1^{-1}.$

Indeed,
$$\bar P(F)\backslash GL_m(F)= \coprod_w w^{-1} U_w(F),$$
where the union is over $w$ of minimal length in $wW_{\bar M}.$
Telescoping, we obtain a sum of terms similar to the right hand side
of \eqref{step1} for these $w.$
Let $\maxunip^M=M\cap \maxunip.$  Observe that $w\maxunip^Mw^{-1} \subset
\maxunip$ for all such $w.$  The restriction of $\psi_W$ to $w\maxunip^Mw^{-1}$
is a generic character iff $wMw^{-1}$ is a standard Levi.  If it is not, the term
corresponding to $w$ vanishes by cuspidality of $\tau.$

On the other hand, $f(w^{-1} ug)$ vanishes if $w^{-1}U_\alpha w$ is
contained in the unipotent radical of $\bar P$ (which we denote $U_{\bar P}$)
for any simple root $\alpha.$  Here $U_\alpha$ denotes the one-dimensional
unipotent subgroup corresponding to the root $\alpha.$
The element $w_1$ is the only element of $W_{GL_m}(\bar M)$ such that
this does not hold for any $\alpha.$

 Let $\lambda$ denote the Whittaker functional on $V_\tau$ given by
 $$\varphi \mapsto \int_{\maxunip^M\quo}\varphi(u) \; \psi_W(w_1uw_1^{-1}) \; du.$$
 Then \eqref{step1} equals
\begin{equation}\label{step2}
\int_{U_{w_1}(\A)}
\lambda(\tilde f(\underline s)(ug))\psi_W(u) \; du,
\end{equation}
where $\tilde f:GL_m(\A) \to V_{\otimes\tau_i}$ is given by
$\tilde f(g)(m) = f(mg)\delta_{\bar P}^{-\frac12}.$   (I.e., $\tilde f$ is the element
of the analogue of $V^{(1)}(\bigotimes_{i=1}^r \tau_i \boxtimes\omega, \underline s),$
corresponding to $f.$)

For each place $v$ there exists a Whittaker functional $\lambda_{v}$
on the local representation $\tau_{v}$ such that
$\lambda(\otimes_v\xi_{v})=\prod_v \lambda_{v}(\xi_{v}).$
(A finite product because $\lambda_v(\xi_v^\circ)=1$ for almost all $v.$
Cf. \cite{Shahidi-ParkCity}, \S1.2.)
The induced representation
$\Ind_{\bar P(\A)}^{GL_m(\A)} (\bigotimes_{i=1}^r \tau_i |\det{}_i|^{s_i}$
is isomorphic to a restricted tensor product of local induced
representations
$\otimes_v{}'\Ind_{\bar P(F_v)}^{GL_m(F_v)} (\bigotimes_{i=1}^r \tau_{i,v} |\det{}_i|_v^{s_i}.$
(Cf. section \ref{s:InducedRepFactors}.)
Consider an element $\tilde f$ which corresponds to a pure tensor $\otimes_v\tilde f_v$ in this
factorization.  So $\tilde f_v(\underline s)$ is a smooth function $GL_m(F_v)\to V_{\otimes\tau_{i,v}}$
for each $\underline s.$)
Then \eqref{step2} equals
 \begin{equation}\label{step3}\prod_v
\int_{U_{w_1}(F_v)}
\lambda_v(\tilde f(\underline s)(u_vg_v))\psi_W(u_v) \; du_v,
\end{equation}
whenever each of the local integrals is convergent, and the infinite
product is convergent (cf \cite{Tate-Thesis} Theorem 3.3.1).
By Propositions 3.1 and 3.2 of \cite{Shahidi-OnCertainLFcns}, all of the
local integrals are always convergent.  (See also Lemma 2.3 and the remark at the
end of section 2 of \cite{Shahidi-Artin}.)

It is an application of Theorem 5.4 of \cite{CS} that the term corresponding to
an unramified nonarchimedean place $v$
in \eqref{step2} is equal to
$W_v^\circ(g_v) \cdot \prod_{i<j}L_v(s_i-s_j+1,\tau_{i,v} \otimes \tilde\tau_{j,v})^{-1}.$
The convergence of the infinite product is then an elementary exercise, as is the main
equation in the statement of our present theorem.

The fact that $f$ may be chosen so that the local Whittaker functions at the places in
$S$ do not vanish follows again from Propositions 3.1 and 3.2 of \cite{Shahidi-OnCertainLFcns}
(see also the remark at the
end of section 2 of \cite{Shahidi-Artin}).

\section{Appendix VI:  Local results on Jacquet Functors}\label{localresults}
In this appendix, $F$ is a non-archimedean local field of characteristic zero
We denote the ring of integers and its unique maximal ideal by $\mathfrak o,$ and $\mathfrak p,$
respectively, and let $q_F:=\#\mathfrak o/\mathfrak p.$  The absolute value on $F$ is normalized so that
its image is $\{q_F^j:\; j\in \Z\}.$
Also,  $\omega$ is an unramified character of $F^\times,$  $\tau$ is an irreducible
unramified principal
series representation of $GL_{2n}(F)$ such that $\tau\iso \tilde \tau \otimes \omega,$
and $\baseaddchar$ is a nontrivial additive character of $F.$

For simplicity, we assume that the characteristic of the residue field $\mathfrak o/\mathfrak p$
is not equal to two.
Hence there are four square classes in $F,$ of which two contain units.
If $\vartheta$ is a character of $N_\ell(F)$ for $1\le \ell \le 2n,$ then we
may define the square class $\Invt(\vartheta)$ as in Definition \ref{d:Invt}
and it is an invariant which separates orbits of characters in general position.
Where convenient, we may restrict attention to those $\vartheta$ such  that
$\Invt(\vartheta)$ contains units, as this condition is
satisfied at almost all places by any global character.  We also
define abstract $F$-groups
$$G_{2n}^{\bf a} \qquad {\bf a} \in F^\times/{(F^\times)^2},$$
and concrete subgroups
$$G_{2n}^a\subset G_{4n+1} \qquad a \in F^\times,$$
such that $G_{2n}^a \iso G_{2n}^{\bf a} \; \forall a \in {\bf a},$
as in Definitions \ref{d:G(sqClass)}, and \ref{d:G(number)}.  The latter is
defined using a character $\Psi_n^a$
given by the same formula as in Definition \ref{d:SpecificCharacters}.

We require
the additional technical hypothesis
\begin{equation}\label{iwasawa}
(B(G_{4n+1})\cap G^{a}_{2n})(F)G_{2n}^{a}(\mathfrak o) = G_{2n}^{a}(F),
\end{equation}
which is known (see \cite{tits}, 3.9, and 3.3.2) to hold at all but
finitely many non-Archimedean completions of a number field.

Throughout this section we shall express certain characters of reductive $F$-groups as
complex linear combinations of rational characters.  The identification is such that
$$\left(\sum_{i=1}^r s_i\chi_i \right)(h) := \prod_{i=1}^r |\chi_i(h)|^{s_i}.$$  Clearly, the
coefficients $s_1, \dots s_r$ appearing in this expression are determined by the
character at most up to $(2\pi i)/\log q_F.$
If $M$ is a Levi, then  restriction gives an injective map $X(M)\to X(T).$  We shall frequently abuse
notation and denote an element of $X(M)$ by the same symbol as its restriction to $T.$
Finally, we let $\Omega$
denote a complex number such that $\omega(x) =|x|^\Omega.$

Lemma \ref{l:combolemma} may be reformulated as stating that  $\tau \iso\Ind_{B(GL_{2n})(F)}^{GL_{2n}(F)}\mu$
for an unramified character $\mu,$ which is of one of the the following two forms:
\begin{equation}\label{e:muFormOne}
\mu_1e_1+\dots+\mu_n e_n +(\Omega-\mu_n)e_{n+1}+\dots +(\Omega-\mu_1)e_{2n}
\end{equation}
\begin{equation}\label{e:muFormTwo}\mu_1e_1+\dots+\mu_{n-1} e_{n-1} + \frac\Omega2e_n +
\left(\frac\Omega2+\frac{\pi i}{\log q_F}\right)e_{n+1}+
(\Omega-\mu_{n-1})e_{n+2}+\dots +(\Omega-\mu_1)e_{2n}.\end{equation}
In either case, by induction in stages,
$$^{un}\Ind_{P(F)}^{G_{4n+1}(F)}\tau\otimes |\det|^{\frac12}\boxtimes \omega
\iso\;{}^{un}\Ind_{B(G_{4n+1})(F)}^{G_{4n+1}(F)}\mu+\frac 12(e_1+\dots +e_{2n})+\Omega e_0.$$
(Here $^{un}$ indicates the unramified constituent, and $P$ the Siegel parabolic of $G_{4n+1}.$)
\begin{rmk}
Because every unramified character is the square of another unramified character, it is possible
to express $\tau$ as a twist of a self-dual representation, and deduce essentially all the results
of this section from the ``classical,'' self-dual case.
\end{rmk}
\begin{lem}\label{l:idOfUnIndsCaseOne}
If $\mu$ is of the form \eqref{e:muFormOne}, then
$$^{un}\Ind_{B(G_{4n+1})(F)}^{G_{4n+1}(F)}\mu+\frac 12(e_1+\dots +e_{2n})+\Omega e_0
\iso \;{}^{un}\Ind_{P_1(F)}^{G_{4n+1}(F)}
\mu'$$
where $P_1$ is the standard parabolic with Levi isomorphic to $GL_2^n\times GL_1,$ such that
the roots of the Levi are $e_{2i-1}-e_{2i}, \; i=1$ to $n,$ and
$\mu'$ is the rational character of this Levi given by
$$\mu':=
\mu_1\det{}_1+\dots +\mu_n\det{}_n+\Omega e_0.
$$
Here
$\det{}_i$ denotes the determinant
of the $GL_2$-factor with unique root $e_{2i-1}-e_{2i}.$
\end{lem}
\begin{proof}
Let $$\tilde \mu= \mu+\frac 12(e_1+\dots +e_{2n})+\Omega e_0=
\sum_{i=1}^n \left(\mu_i+\frac 12\right)e_i+\sum_{i=1}^n \left(\Omega-\mu_i+\frac 12\right)e_{n+i}
+\Omega e_0.$$
Using the description of the Weyl action in
Lemma \ref{l:WeylAction}
it is easily verified that this is in the same orbit as
$$\tilde \mu' := \sum_{i=1}^n \left[ \left(\mu_i+\frac 12\right)e_{2i-1}
+\left(\mu_i-\frac 12\right)e_{2i} \right] +\Omega e_0.$$
By the definition of the unramified constituent, then,
$$^{un}\Ind_{B(G_{4n+1})(F)}^{G_{4n+1}(F)}\tilde \mu=\;{}
^{un}\Ind_{B(G_{4n+1})(F)}^{G_{4n+1}(F)}\tilde \mu'.$$
The lemma now follows from the well known (and easily verified) fact that
\begin{equation}\label{e:GL2UnramFact}
^{un}Ind_{B(GL_2)(F)}^{GL_2(F)} (\mu+\frac12)e_1'+ (\mu-\frac12)e_2'
=^{un}Ind_{B(GL_2)(F)}^{GL_2(F)} (\mu-\frac12)e_1'+ (\mu+\frac12)e_2'
= \mu \det,\end{equation}
where $e_1'$ and $e_2'$ are the usual basis for the lattice of rational
characters of the torus of diagonal elements of $GL_2.$
\end{proof}
The next lemma is similar, but slightly more complicated.  It makes use of alternative
$\Z$-bases of the lattices of characters and cocharacters.  Specifically,
$\{e_1, \dots e_{2n-2}, f_1,f_2,f_0\}, \{e_1^*, \dots e_{2n-2}^*, f_1^*, f_2^*, f^*_0\},$ where
\begin{eqnarray*}
e_0=-f_1&\qquad& e_0^*=-2f_0^*-f_1^*-f_2^*\\
e_{2n-1}=-f_0+f_1+f_2&\qquad& e_{2n-1}^*=-f_0^*\\
e_{2n}=f_1-f_2&\qquad& e_{2n}^*=-f_0^*-f_2^*.
\end{eqnarray*}
The key feature of these $\Z$-bases is as follows.  Recall that the group $G_{4n+1}$ has
a unique standard Levi isomorphic to $GL_2^{n-1}\times G_5,$ with the based root datum
of the $G_5$ component lying in the sublattices spanned by $\{e_{2n-1},e_{2n},e_0\},$
$\{e_{2n-1}^*,e^*_{2n},e^*_0\}.$  Now, $G_5$ and
$GSp_4$ are the same $F$-group.
When we write the based root datum of this Levi with respect to the new basis, the expression
for the $G_5$ component matches the ``standard form'' for the based root datum of $GSp_4$
as in section \ref{rootDataSection}.  In particular, the character $f_0$ is the restriction to the
torus of $GSp_4$ of the similitude factor (which is a generator for the rank-one lattice of rational
characters of $GSp_4$), and there is a standard Levi, isomorphic to  $GL_2$ such that its unique root
is $f_1-f_2.$
\begin{rmks}
To avoid confusion, let us draw attention the following tricky point:  we have defined a notion of
``Siegel parabolic'' and ``Siegel Levi'' for $G_{2n+1},$ any $n.$  There is also a well known
notion of ``Siegel parabolic'' and ``Siegel Levi'' for $GSp_{2n},$ any $n,$ which is very
widespread in the literature.  The two groups $G_5$ and $GSp_4$ happen to coincide,
{\em and the two notions of ``Siegel parabolic'' and ``Siegel Levi'' do {\bf not.}}
\end{rmks}

\begin{lem}\label{l:idOfUnIndsCaseTwo}
If $\mu$ is of the form \eqref{e:muFormTwo},
and $\tilde \mu$ is defined in terms of $\mu$ as in the proof of Lemma \ref{l:idOfUnIndsCaseOne},
then
$$^{un}\Ind_{B(G_{4n+1})(F)}^{G_{4n+1}(F)}\tilde \mu
\iso \;{}^{un}\Ind_{P_2(F)}^{G_{4n+1}(F)}
\mu''$$
where
\begin{eqnarray*}
\mu''&=&\sum_{i=1}^{n-1}\mu_i\det{}_i-\frac{\Omega-1}{2}f_0
+\left(-\frac 12+\frac{\pi i}{\log q_F}\right) \det{}_0\\
&=&\sum_{i=1}^{n-1}\mu_i\det{}_i+\frac{\Omega-1}{2}e_{2n+1} - \left(
\frac \Omega 2+\frac{\pi i }{\log q_F}\right) \det{}_0
\end{eqnarray*}
where the notation is as follows:  $P_2$ is the standard parabolic with Levi isomorphic to
$GL_2^n\times GL_1,$ such that
the roots of the Levi are $e_{2i-1}-e_{2i}, \; i=1$ to $n-1,$ and $e_{2n}.$
(One might also describe this Levi as
$GL_2^{n-1}\times GL_1 \times GSpin_3.$)
As in Lemma  \ref{l:idOfUnIndsCaseOne}
$\det{}_i$ denotes the determinant
of the $GL_2$-factor with unique root $e_{2i-1}-e_{2i},$ for $i=1$ to $n-1,$ while
$\det{}_{0}$ denotes the determinant of the $GL_2$ with unique root $e_{2n}=f_1-f_2.$
\end{lem}
\begin{proof}
This time
$\tilde \mu$ is in the same Weyl orbit as
\begin{eqnarray*}
\tilde \mu'' &:=& \sum_{i=1}^{n-1} \left[ \left(\mu_i+\frac 12\right)e_{2i-1}
+\left(\mu_i-\frac 12\right)e_{2i} \right] +\left(\frac{\Omega-1}2\right)e_{2n-1}
+\left(\frac{\Omega-1}2+\frac{\pi i}{\log q_F}\right)e_{2n}+\Omega e_0\\
&=& \sum_{i=1}^{n-1} \left[ \left(\mu_i+\frac 12\right)e_{2i-1}
+\left(\mu_i-\frac 12\right)e_{2i} \right]
-\left(\frac{\Omega-1}2\right)f_0+\left(-1+\frac{\pi i}{\log q_F}\right)f_1-\frac{\pi i}{\log q_F}f_2.
\end{eqnarray*}

Using \eqref{e:GL2UnramFact} again, in conjunction with the fact that $-\frac{\pi i}{\log q_F}f_2$
and $\frac{\pi i}{\log q_F}f_2$ are the same character, we obtain the lemma.
\end{proof}
Next, we need a slight extension of this.
Let
$P_3$ be the  standard parabolic of $G_{4n+1}$ with Levi isomorphic to $GL_2^{n-1}\times GSp_4.$
Identify $GSp_4$ with the component of this Levi, and let $R=GSp_4 \cap P_2.$  This is the subgroup
known in the literature as the ``Siegel'' parabolic of $GSp_4.$  When regarded as a parabolic of $GSpin_5,$ it is the one for which we have introduced the notation $Q_1=L_1N_1.$
Its lattice of rational characters is
spanned by $f_0$ and $\det{}_0,$ defined as in Lemma \ref{l:idOfUnIndsCaseTwo}.
Let $\pi_0= \;^{un}\Ind_{R(F)}^{GSp_4(F)} \left(\frac 12+\frac{\pi i}{\log q_F}\right)\det{}_0.$
Extend $\pi_0$ trivially to a representation of the Levi of $P_3.$
\begin{cor}\label{c:indInStages}
$$^{un}\Ind_{B(G_{4n+1})(F)}^{G_{4n+1}(F)}\tilde \mu''
\iso
  \;{}^{un}\Ind_{P_3(F)}^{G_{4n+1}(F)}
\mu'''\otimes \pi_0,$$
where
$$\mu''':=\left(  \sum_{i=1}^{n-1}\mu_i\det{}_i-\frac{\Omega-1}{2}f_0\right).$$
\end{cor}
\begin{proof}
Induction in stages and the definition of the unramified constituent.
\end{proof}
An important fact about $\pi_0$ is the following:
\begin{lem}\label{l:pi0asubrep}
The representation $\pi_0$ may be realized as a subrepresentation of
$$\Ind_{R(F)}^{GSp_4(F)} \left(-\frac 12+\frac{\pi i}{\log q_F}\right)\det{}_0.$$
\end{lem}
\begin{proof}
In fact, it is one of the
spaces $R_2(V)$ introduced on p. 223 of \cite{kudlarallis}.  This can be checked by
direct computation.  It also follows from Proposition 5.5 of \cite{kudlarallis},
in that the intertwining operator
is easily seen not to vanish on the spherical vector.
\end{proof}
\begin{cor}
The representation $\Ind_{P_3(F)}^{G_{4n+1}(F)}
\mu'''\otimes \pi_0$
may be realized as a subrepresentation of
$$\Ind_{P_2(F)}^{G_{4n+1}(F)}
\left[\mu'''
+\left(-\frac 12+\frac{\pi i}{\log q_F}\right) \det{}_0\right].$$
\end{cor}
A second important fact about the representation $\pi_0$ is the following:
\begin{lem}\label{l:kudlarallis}
Let $\vartheta$ be a character of the unipotent radical of $R$ in general position.  Regarding $R$
as the parabolic $Q_1$ of $G_5,$ the square class $\Invt(\vartheta)$ is defined.  A sufficient condition
for the vanishing of the twisted Jacquet module $\mathcal{J}_{N_1,\vartheta}(\pi_0)$ is that
the Hilbert symbol $(\cdot, \Invt(\vartheta))$ not equal the unique nontrivial unramified quadratic
character.
\end{lem}
\begin{proof}
This follows from \cite{kudlarallis}, Lemma 3.5 (b), p. 226.  (Here, we again use the fact
 that the unramified
constituent of $\Ind_{R(F)}^{GSp_4(F)} \left(\frac 12+\frac{\pi i}{\log q_F}\right)\det{}_0$ is one
of the spaces $R_2(V)$ introduced on p. 223 of \cite{kudlarallis}.)
\end{proof}

\begin{prop}\label{p:JacVanishes1}
Let $\tau = Ind_{B(GL_{2n})(F)}^{GL_{2n}(F)} \mu,$
with $\mu$ of the form \eqref{e:muFormOne}, and let $P$ denote the Siegel parabolic subgroup.
Then for $\ell >n$ and $\vartheta$ in general postion,
the Jacquet module
$\mathcal{J}_{N_\ell, \vartheta}(
^{un}Ind_{P(F)}^{G_{4n+1}(F)} \tau \otimes |\det |^{\frac12} \boxtimes
\omega )$
is trivial.
The same is true if $\ell =n$ and $\Invt(\vartheta) \ne \square.$
\end{prop}
\begin{proof}
By Lemma \ref{l:idOfUnIndsCaseOne}, it suffices to prove that the corresponding
Jacquet module of $\Ind_{P_1(F)}^{G_{4n+1}(F)}  \mu'$ vanishes.
This is in essence an application of theorem 5.2
of \cite{BZ-ASENS}.
The space $Ind_{P_1(F)}^{G_{4n+1}(F)}  \mu' $ has
a filtration as a $\descentparabolic{\ell}(F)$-module, in terms of
$\descentparabolic{\ell}(F)$-modules indexed by the
elements of $(W \cap \middlestageparabolic) \backslash W / (W\cap \descentparabolic{\ell}).$
For any element $x$ of $\middlestageparabolic(F) w \descentparabolic{\ell}F)$ the module
corresponding to $w$ is isomorphic to $c-ind_{ x^{-1} \middlestageparabolic(F) x
\cap \descentparabolic{\ell}(F) }^{\descentparabolic{\ell}(F)} (\mu'+\rho_{P_1}) \circ Ad( x).$
Here $Ad(x)$ denotes the map given by conjugation by $x.$  It sends
$x^{-1} \middlestageparabolic(F) x
\cap \descentparabolic{\ell}(F)$ into $\middlestageparabolic(F).$
Also, here and throughout $c-ind$ denotes non-normalized compact induction.
(See \cite{cassnotes}, section 6.3.)

Recall from \ref{s:WeylGroup} that the elements of
the Weyl group of $G_{4n+1}$ are
(after the choice of $\pr$) in natural one-to-one correspondence  with
the set of permutations $w \in \mathfrak S_{4n+1}$ satisfying,

(1) $w( 4n+1-i) = 4n+1-w(i)$

As representatives for the double
cosets $(W \cap \middlestageparabolic) \backslash W / (W\cap \descentparabolic{\ell})$ we choose the element of minimal
length in each.  The permutations corresponding to  these elements satisfy

(2) $w^{-1}( 2i ) > w^{-1}(2i-1)$ for $i=1$ to $2n,$ and

(3) $\ell < i <j < 4n+2-\ell \; \Longrightarrow w(i) < w(j).$

Let
$I_ w$
 be the $\descentparabolic{\ell}(F)$-module obtained as $$c-ind_{ \dot w^{-1} \middlestageparabolic(F) \dot w
\cap \descentparabolic{\ell}(F) }^{\descentparabolic{\ell}(F)} \left(\mu'+\rho_{P_1}\right)\circ Ad( \dot w)$$
using any element $\dot w$ of $\pr^{-1}( \det w \cdot w).$ (Cf. section \ref{s:WeylGroup}.)

A function $f$ in $I_w$ will map to zero under the natural projection
to $\mathcal{J}_{N_\ell, \vartheta}(I_w)$ iff there exists a compact subgroup $N^0_\ell$ of $N_\ell(F)$
 such that
$$\int_{N^0_\ell} f( hn) \overline{ \vartheta ( n ) } dn = 0 \qquad \forall h \in \descentparabolic{\ell}(F).$$
(See \cite{cassnotes}, section 3.2.)
Let $h\cdot \vartheta (n)= \vartheta(h^{-1}n h).$  It is easy to see that the integral above
vanishes for suitable $N^0_\ell$ whenever
\begin{equation}\label{varthetah}
h\cdot \vartheta |_{N_\ell(F) \cap w^{-1} P_1(F) w } \text{ is nontrivial. }
\end{equation}
Furthermore, the function $h \mapsto h\cdot \vartheta$ is continuous in $h,$
(the topology on the space of characters of $N_\ell(F)$ being defined by
identifying it with a finite dimensional $F$-vector space, cf. section \ref{s:uniper})
so if this condition holds for all $h$ in a compact set, then $N^0_\ell$
can be made uniform in $h.$

Now, $\vartheta$ is in general position.  Hence, so is $h\cdot \vartheta$ for
every $h.$  So, if we write
$$h\cdot \vartheta( u ) = \baseaddchar ( c_1 u_{1,2} + \dots + c_{\ell-1} u_{\ell-1,\ell}
+ d_1 u_{\ell , \ell+1} + \dots + d_{4n-2\ell+1} u_{\ell, 4n-\ell+1}),$$
we have that $c_i \neq 0$ for all $i$ and $^t \underline{d} J \underline{d} \neq 0.$

Clearly, the condition \eqref{varthetah} holds for all $h$ unless

(4) $w(1) > w(2) > \dots > w( \ell).$

Furthermore, because $^t \underline{d} J \underline{d} \neq 0,$ there exists some
$i_0$ with $\ell+1 \leq i_0 \leq 2n$ such that $d_{i_0-\ell} \neq 0$ and $d_{4n+2+\ell - i_0} \neq 0.$
From this we deduce that the condition \eqref{varthetah} holds for all $h$ unless
$w$ has the additional property

(5) There exists $i_0$ such that $w(\ell) > w(i_0)$ and $w(\ell) > w( 4n+2-i_0).$

However, if $\ell > n$ it is easy to check that no permutations with properties
(1),(2), (4) and (5) exist.

Thus $\mathcal{J}_{N_\ell, \vartheta}(I_w)=\{0\}$  for all $w$ and hence
$ \mathcal{J}_{N_\ell, \vartheta}(
^{un}Ind_{P(F)}^{G_{4n}(F)} \tau \otimes |\det |^{\frac12} \boxtimes
\omega )=\{0\}$ by exactness of the Jacquet functor.

If $\ell =n,$ there is exactly one permutation $w$ which satisfies (1)-(4).  For this
permutation, condition (4) is satisfied only with $i_0=4n+2-i_0=2n+1.$  The orbit
of $\vartheta$ contains characters such that $d_i=0$ for all $i\ne 2n+1$ iff $\Invt(\vartheta)=\square.$
\end{proof}

\begin{prop}\label{p:JacVanishes2}
Let $\tau = Ind_{B(GL_{2n})(F)}^{GL_{2n}(F)} \mu,$
with $\mu$ of the form \eqref{e:muFormTwo}, and let $P$ denote the Siegel parabolic subgroup.
Then for $\ell >n$ and $\vartheta$ in general postion,
the Jacquet module
$\mathcal{J}_{N_\ell, \vartheta}(
^{un}Ind_{P(F)}^{G_{4n+1}(F)} \tau \otimes |\det |^{\frac12} \boxtimes
\omega )$
is trivial.  The same is true if $\ell =n$ and $\Invt(\vartheta) = \square.$
\end{prop}
\begin{proof}
For $\ell >n,$ the proof is similar to that of Proposition \ref{p:JacVanishes1}.  Using Lemma
\ref{l:idOfUnIndsCaseTwo} in place of Lemma \ref{l:idOfUnIndsCaseOne}, we
consider a representation induced from a character of $P_2$ rather than $P_1.$  The
effect is that in place of condition (2) from the proof of Proposition  \ref{p:JacVanishes1}, we
have the condition

(2$'$) $w^{-1}(2i-1)< w^{-1}(2i), \;1\le i <n, \quad w^{-1}(2n)<w^{-1}(2n+1).$

The set of permutations satisfying (1),(2$'$),(3),(4) is again empty.

The proof of  vanishing when $\ell =n$ and $\Invt(\vartheta)=\square$ is more
nuanced.  In this case we use both Lemma \ref{l:idOfUnIndsCaseTwo}
and Corollary \ref{c:indInStages}, obtaining {\it two} filtrations
of $$\Ind_{P_3(F)}^{G_{4n+1}(F)}
\mu'''\otimes \pi_0
\subset \;{}
\Ind_{P_2(F)}^{G_{4n+1}(F)}
\mu'',$$
indexed by  $(W \cap P_3) \backslash W / (W\cap \descentparabolic{\ell})$
and  $(W \cap P_2) \backslash W / (W\cap \descentparabolic{\ell}).$
The latter is a refinement of the former, in a manner which is
compatible with the natural projection
$$(W \cap P_2) \backslash W / (W\cap \descentparabolic{\ell})\to
(W \cap P_3) \backslash W / (W\cap \descentparabolic{\ell}).$$
Let us denote the elements of the first filtration by
$I_w, \;w \in (W \cap P_3) \backslash W / (W\cap \descentparabolic{\ell}),$
and the elements of the second by
$I'_w, \; w \in (W \cap P_2) \backslash W /(W\cap \descentparabolic{\ell}).$

Now, when $\ell =n$ there is a unique permuation $w_0$ satisfying (1)(2$'$),(3), (4),(5).
It is the shortest element of the double coset containing the longest element of $W.$
It follows that $\mathcal{J}_{N_n,\vartheta}(I_w')$ vanishes for every $w\ne w_0,$
and hence that  $\mathcal{J}_{N_n,\vartheta}(I_w)$ vanishes for every $w$ other than
the shortest element of $(W \cap P_3) \cdot w_0 \cdot (W\cap \descentparabolic{n}),$
which we denote $w_0'.$

The permutation $w_0'$ can be described explicitly as follows:
$$w_0'(i) =
\begin{cases}
4n+2-2i & 1\le i \le n-1,\\
2n -1 & i=n,\\
2i-2n-1 &  n+1\le i \le 2n-1, \; 2n+3 \le i \le 3n+1,\\
i & 2n \le i \le 2n+2,\\
2n+3 & i=3n+2,\\
8n+4-2i & 3n+3 \le i \le 4n+1.
\end{cases}
$$

Furthermore, the space $I_{w_0'}$ is equal to the subspace of
$^{un}\Ind_{P_3(F)}^{G_{4n+1}(F)}
\mu'''\otimes \pi_0$
consisting of smooth functions having support in the open double coset
$P_3(F)\cdot  w_0' \cdot Q_n(F).$  Take such a function $f$ and take $N_n^0\subset N_n(F),$
compact.  Consider the integral
$$\int_{N_n^0}f(g n) \overline{\vartheta(n)} \; dn.$$
We may assume $g = w_0' q$ for some $q \in Q_n(F).$  Then we get
$$\int_{qN_n^0q^{-1}}f( w_0 n q) \overline{q\cdot  \vartheta(n)}\; dn,$$
where $q \cdot \vartheta(n) = \vartheta(q^{-1}nq).$
Hence, we consider
\begin{equation}\label{e:JacVanInt}
\int_{{N_n^0}'}f'(w_0n) \overline{\vartheta'(n)}\; dn,
\end{equation}
for $\vartheta'$ a character of $N_n$ such that $\Invt(\vartheta')=\square,$ $f'\in  I_{w_0'},$
and ${N_n^0}' \subset N_n(F)$ compact.
Observe that $w_0 N_n w_0^{-1}$ contains the unipotent radical $U_R$ of the parabolic $R$ of
$GSp_4$ used to define $\pi_0.$
Indeed, if  $\hat N_n = \{ u \in N_n: u_{n,2n}=u_{n,2n+1} = 0\},$ then $\hat N_n
\normal N_n$ and $N_n = w_0^{-1}U_R w_0 \cdot \hat N_n.$
  If $U\subset \maxunip,$ write $U(\mathfrak p^N)$ for
$\{ u \in U: u_{ij}\in \mathfrak p^N \forall i,j\}.$

For each $h \in G_{4n+1}(F),$ the function $g \mapsto f'(gh), \; g\in GSp_4(F)$ is an
element of $\pi_0.$  By Lemma \ref{l:kudlarallis}, for each $h$ there exists
$N$ such that
$$\int_{w_0^{-1}U_R(\mathfrak p^N ) w_0 } f'(w_0 u h) \vartheta'(u) \; du =0.$$

Clearly, $N$ depends on $f'$ and $\vartheta',$ and hence, if $f'(g)=f(g\cdot q)$ and
$\vartheta'=q\cdot \vartheta,$ on $q.$  However,  $f$ is smooth and has
support which is compact modulo $P_3(F),$ so $f'$ takes only finitely many values. Furthermore,
the $q\cdot \vartheta$ is a continuous function of $q$ in the sense discussed above.  Thus, $N$ may
be made uniform in $q.$
\end{proof}
Define a character $\Psi_n$ of $N_n(F)$ by the  same formula as in
Definition \ref{d:SpecificCharacters}.
In the proof of Lemma \ref{l:G^aBorel}, we fixed a specific isomorphism $\inc:G_{2n}
\to (L_n^{\Psi_n})^0.$  For the next proposition only, we let
$B$ denote the image under $\inc$ of the Borel $B(G_{2n})$
corresponding to our choices of maximal torus and simple roots for $G_{2n}.$  It is
equal to $(L_n^{\Psi_n})^0\cap B(G_{4n+1}).$  The corresponding maximal
torus $T$  is the subtorus $\langle e_i^*: i=0, \text{ or }n+1\le i \le 2n\rangle.$
Because of this $\sum_{i=0}^{2n} c_i e_i$ makes sense as a character
of $T(F).$  (But depends only on $c_i, \;i=0, \text{ or }n+1\le i \le 2n.$)
\begin{prop}\label{p:JFcompSplit}
Let $P_1,$ and $\mu'$ be defined as in Lemma \ref{l:idOfUnIndsCaseOne}.
Then we have isomorphisms
$$\mathcal{J}_{N_n, {\Psi_n}}(
\Ind_{P_1(F)}^{G_{4n+1}(F)}\mu')
\iso \Ind_{B(F)}^{(L_n^{\Psi_n})(F)}\mu^*\iso
 \Ind_{B(F)}^{(L_n^{\Psi_n})(F)}\mu^{**}
\qquad (\text{ of }L_n^{\Psi_n}-\text{modules}),$$
$$\mathcal{J}_{N_n, {\Psi_n}}(
\Ind_{P_1(F)}^{G_{4n+1}(F)}\mu')
\iso \Ind_{B(F)}^{(L_n^{\Psi_n})^0(F)}\mu^*\oplus
 \Ind_{B(F)}^{(L_n^{\Psi_n})^0(F)}\mu^{**}
\qquad (\text{ of }(L_n^{\Psi_n})^0-\text{modules}),$$
where
$$\mu^*=\sum_{i=1}^n \mu_ie_{n+i}+\Omega e_0, \qquad \mu^{**} =\sum_{i=1}^{n-1}
\mu_i e_{n+i} +(\Omega - \mu_n)e_{2n}+\Omega e_0.
$$
\end{prop}
\begin{proof}
As before,  we filter $\Ind_{P_{1}(F)}^{G_{4n+1}(F)}\mu' $ in terms of $Q_{n}(F)$-modules $I_w.$
This time, $\mathcal{J}_{N_{n}, {\Psi_n}}(I_w)=\{0\}$ for all $w$
except possibly for one.  This one Weyl element, which we denote $w_0,$
corresponds to
the unique permutation satisfying (1) and (2)
of Proposition \ref{p:JacVanishes1}, together with
$w_0(i) = 4n-2i+2$ for $i=1$ to $n.$
Exactness yields
$$\mathcal{J}_{N_{n},\Psi_n}\left(
^{un}Ind_{P(F)}^{G_{4n+1}(F)}\tau \otimes |\det |^{\frac12} \boxtimes
\omega \right) \iso \mathcal{J}_{N_n, {\Psi_n}}(I_{w_0}).$$
(This is an isomorphism of $Q_{n}^{{\Psi_n}}(F)$-modules, where $Q_{n}^{{\Psi_n}}=N_{n} \cdot L_{n}^{{\Psi_n}} \subset Q_{n},$
is the stabilizer of ${\Psi_n}$ in $Q_{n}$ (cf. $L^\vartheta$ above).)

Now, recall that for each $h \in Q_{n}(F)$ the character
$h\cdot \Psi_n( u) = \Psi_n(h^{-1} u h)$ is a character of $N_{n}$ in
general position, and as such determines coefficients
$^hc_1, \dots , ^hc_{n-1}$ and $^hd_1, \dots, ^hd_{2n+1}$ as in  \eqref{e:Action}.
Clearly,
$$Q_{n}^{o}:= \left\{ \left.h \in Q_{n}(F)\right|
d_i^h \neq 0 \text{ for some }i \neq n+1,
\right\}$$
is open.
Moreover, one may see from the description of $w_0$ that
for $h$ in this set the condition \eqref{varthetah}, which assures vanishing, is satisfied.

We have an  exact sequence of $Q_{n}^{\Psi_{n}}(F)$-modules
$$0\to   I_{w_0}^* \to I_{w_0} \to \bar I_{w_0} \to 0,$$
where $I_w^*$ consists of those functions in $I_w$ whose
compact support happens to be contained in $Q_{n}^{o},$
and the third arrow is restriction to the complement of $Q_{n}^{o}.$
This complement is slightly larger than
$Q_{n}^{\Psi_{n}}(F)$ in that it contains the full torus
of $Q_{n}(F),$ but restriction of functions is an isomorphism
of $Q_{n}^{\Psi_{n}}(F)$-modules,
$$\bar I_{w_0} \to c-ind_{Q_{n}^{\Psi_{n}}(F) \cap w_0^{-1} P_{1} (F)w_0 }^{Q_{n}^{\Psi_{n}}(F)}
\left(\mu' +\rho_{P_{1}}\right) \circ Ad( w_0).$$

Clearly $\mathcal{J}_{N_{n},\Psi_{n}}\left(I_{w_0}^*\right) =\{0\},$
and hence we have the isomorphism
$$\mathcal{J}_{N_{n},\Psi_{n}}\left(Ind_{P_{1}(F)}^{G_{4n+1}(F)}\mu'
\right)
\iso
\mathcal{J}_{N_{n},\Psi_{n}}\left(
c-ind_{Q_{n}^{\Psi_{n}}(F) \cap w_0^{-1} P_{1} (F)w_0 }^{Q_{n}^{\Psi_{n}}(F)}
\left(\mu' +\rho_{P_{1}}\right)\circ Ad( w_0)
\right)
$$ of $Q_n^{\Psi_n}$-modules.

 Let us denote
$$c-ind_{Q_{n}^{\Psi_{n}(F)} \cap w_0^{-1} P_{1}(F) w_0 }^{Q_{n}^{\Psi_{n}}(F)}
\left(\mu' +\rho_{P_{1}}\right) \circ Ad( w_0)$$ by $V.$
A straightforward computation shows that the functions in
 $V$
satisfy
$$f(bq) = b^{\mu^*+\rho_B-J} f(q) \qquad \forall b \in B(F), \; q\in Q_n^{\Psi_n}(F),$$
where
$$J = \sum_{i=1}^n (i-n-1)e_{n+i}.$$

For $f \in V,$
let
$$W(f)(q) = \int_{N_n(F) \cap w_0^{-1} \overline{\maxunip}(F) w_0 } f( u q) \bar \Psi_{n} ( u )
du.$$
Then the character $J$ computed above matches exactly  the Jacobian of $Ad(b),\; b \in B(F),$
 acting on $N_{n}(F) \cap w_0^{-1} \overline{\maxunip}(F) w_0.$
It follows that
$$W(f)(bq)=b^{\mu^*+\rho_B}f(g) \qquad \forall \; b \in B(F), \; q \in Q(F).$$

Now let $\mathcal{W}$ denote
$$\left\{ f: Q_{n}^{\Psi_{n}}(F) \to \C\left|
\begin{array}{l}
f( uq) = \Psi_{n}( u ) f( q)\; \forall\; u \in N_{n}(F), \;q \in Q_{n}^{\Psi_{n}}(F),\\
f( b m ) = b^{\mu^{*}+\rho_{B}}
 f( m)\; \forall\; b \in B(F), \;m \in L_n^{\Psi_n}(F)
 \end{array}\right.
 \right\}.$$
Then $W$ maps $V$ into $\mathcal{W}.$

Denote by $V(N_{n}, \Psi_{n})$
the kernel of the linear
map $V \to \mathcal{J}_{N_{n}, \Psi_{n}}(V).$
It is easy to show that $V(N_{n}, \Psi_{n})$ is contained in the kernel of $W.$
In the Lemma \ref{l:kerW=V} below, we show that in fact, they are equal.
Restriction from $Q_{n}^{\Psi_{n}}(F)$ to $L_{n}^{\Psi_{n}}(F)$ is clearly an isomorphism
$\mathcal{W} \to
\Ind_{B(F)}^{L_{n}^{\Psi_{n}}(F)}\mu^{*}.$

The proof that this is isomorphic to $\Ind_{B(F)}^{L_{n}^{\Psi_{n}}(F)}\mu^{** }$
and decomposes into $(L_n^{\Psi_n})^0$-modules in the manner described is straightforward.
\end{proof}
The next proposition is similar.   However, there is an interesting difference between
the two.  In the previous proposition, we let $B$ denote the Borel subgroup
$B(G_{4n+1})\cap (L_n^{\Psi_n})^0$ of $(L_n^{\Psi_n})^0\iso G_{2n}^{\bf 1}.$
For the next, we use it to denote $Q_{2n-1} \cap G_{2n}^a,$ which is a Borel subgroup
of $G_{2n}^a.$  The corresponding maximal torus, $G_{2n}^a \cap L_{2n-1},$
is given by
$$
\left\{ h_a \prod_{i=1}^{n-1}e_{n+i}^*(t_i) \cdot e_{2n}^*\left((x+y\sqrt{a})\cdot (x-y\sqrt{a})^{-1}\right)
e_0(x-y\sqrt{a})h_a^{-1} : t_i \in F, \; x, y\in F, x^2-ay^2\ne 0\right\},
$$
as in Lemma \ref{l:G^aBorel}\eqref{maxtor}.
Here $\sqrt{a}$ may be taken to be either of the solutions to $\zeta^2=a$ in the algebraic
closure of $F.$  We assume $\sqrt{a}\notin F.$
The lattice of $F$-rational
characters of this torus is $\langle e_{n+i}:\; 1\le i \le n-1 , \; e_{2n}+2e_0\rangle.$
The character $e_{2n}+2e_0$ is the restriction of a rational character of the
$L_{2n-1} \iso GL_1^{2n-1} \times GSpin_3.$
To be precise, it is the {\it inverse} of the character $\det_0$ introduced earlier.
(Cf. Lemma \ref{l:idOfUnIndsCaseTwo}.)  Thus, a general  rational
character of this torus may be expressed as
$$\sum_{i=1}^{n-1} c_i e_{n+i} + c_0 \det{}_0,$$
with $c_i \in \Z.$
In particular map, the  restriction map from $X(L_{2n-1})$ is surjective.
A general unramified character of this torus may be expressed in the same form
with $c_i \in \C.$  Then $c_i$

Observe that for any $t$ in this torus $\det{}_0(t)$ is a norm from $F(\sqrt{a}).$
When $a$ is in the square class which contains the non-square units (i.e.,
when $F(\sqrt{a})$ is the unique unramified quadratic extension of $F,$) the absolute
value of a norm is always an even power of $q_F,$ and so
$c_0$ is defined only up to $\frac{\pi i}{\log q_F}.$
(whereas the others are  defined up to $\frac{2 \pi i}{\log q_F}$ for $1\le i \le n-1.$)

We also
let $\tilde B$ denote $Q_{2n-1} \cap L_n^{\Psi_n^a}.$  (Recall that $G_{2n}^a := (L_n^{\Psi_n^a})^0.$)
 It is not difficult to see that $L_{2n-1} \cap L_n^{\Psi_n^a}$ is properly larger that
 $L_{2n-1} \cap (L_n^{\Psi_n^a})^0,$ i.e., contains elements of the non-identity component of
 $L_n^{\Psi_n^a}.$ A character of $B$ may be extended trivially to $\tilde B.$  And any
 character of $\tilde B$ which is obtained as the restriction of a character of $Q_{2n-1}$
 is such a trivial extension.

\begin{prop}\label{p:JFcompNonSplit}
Let $P_2,$ and $\mu''$ be defined as in Lemma \ref{l:idOfUnIndsCaseTwo}.
Then we have isomorphisms
$$\mathcal{J}_{N_n, {\Psi_n}}(
\Ind_{P_2(F)}^{G_{4n+1}(F)}\mu'')
\iso \Ind_{\tilde B(F)}^{(L_n^{\Psi_n^a})(F)}\mu^*
\qquad (\text{ of }L_n^{\Psi_n^a}-\text{modules}),$$
$$\mathcal{J}_{N_n, {\Psi_n}}(
\Ind_{P_2(F)}^{G_{4n+1}(F)}\mu'')
\iso \Ind_{B(F)}^{G_{2n}^a(F)}\mu^*
\qquad (\text{ of }(L_n^{\Psi_n})^0-\text{modules}),$$
where
$$\mu^*=\sum_{i=1}^{n-1} \mu_ie_{n+i}-\left( \frac \Omega2 +\frac{\pi i}{\log q_F}\right)\det{}_0.
$$
\end{prop}
\begin{proof}
We use Lemma \ref{l:idOfUnIndsCaseTwo}, and filter by $Q_n$-modules.  As in
Proposition \ref{p:JFcompSplit}, there is a unique permutation $w_1$ such that
the corresponding $Q_n$-module $I_{w_1}$ does not vanish.   This permutation is
given by
$$w_1(i) =\begin{cases}
4n+2-2i& 1 \le i \le n-1,\\
2n+3& i = n,\\
2i-2n-1& n+1 \le i \le 2n-1,\\
i & 2n \le i \le 2n+2, \\
2i-2n-1&2n+3\le i \le 3n+1,\\
2n-1&i=3n+2,\\
2(4n+2-i)&3n+3\le i \le 4n+1.
\end{cases}$$
The group $Q_n \cap w_1^{-1}P_2 w_1$ contains $L_{2n-1}.$  Since
$L_{2n-1} \cdot Q_n^{\Psi_n^a} = Q_n,$ restriction of functions is an
isomorphism of $Q_n^{\Psi_n^a}$-modules,
$$I_{w_1}\to c-ind_{Q_n^{\Psi_n^a}\cap w_1^{-1}P_2 w_1}^{Q_n^{\Psi_n^a}} (\mu'' +\rho_{P_2}) \circ \Ad(w_1).$$

 This time, let $V$ denote
$$ c-ind_{Q_n^{\Psi_n^a}\cap w_1^{-1}P_2 w_1}^{Q_n^{\Psi_n^a}} (\mu'' +\rho_{P_2}) \circ \Ad(w_1).$$
Once again the functions in
 $V$
satisfy
$$f(bq) = b^{\mu^*+\rho_B-J} f(q) \qquad \forall b \in B(F), \; q\in Q_n^{\Psi_n}(F),$$
with $J$ as before.
We define
$$W(f)(q) = \int_{N_n(F) \cap w_1^{-1} \overline{\maxunip}(F) w_1 } f( u q) \bar \Psi_{n} ( u )
du,$$
and find that $W$ maps $V$ to
$$\mathcal{W}:=\left\{ f: Q_{n}^{\Psi_{n}}(F) \to \C\left|
\begin{array}{l}
f( uq) = \Psi_{n}( u ) f( q)\; \forall\; u \in N_{n}(F), \;q \in Q_{n}^{\Psi_{n}}(F),\\
f( b m ) = b^{\mu^{*}+\rho_{B}}
 f( m)\; \forall\; b \in B(F), \;m \in L_n^{\Psi_n}(F)
 \end{array}\right.
 \right\},$$
 which is easily seen to be isomorphic to each of the induced representations specified.
As before,
the kernel of the linear
map $V \to \mathcal{J}_{N_{n}, \Psi_{n}}(V)$
 is contained in the kernel of $W.$
In  Lemma \ref{l:kerW=V}, we show that in fact, they are equal to complete the proof.
\end{proof}

\begin{lem}\label{l:kerW=V}
Let $\vartheta$ be a character of $N_n$ in general position, $H$
its stabilizer in $L_n,$ $U_1$ and $U_2$ two subgroups of
$N_n$ such that $U_1\cap U_2=1$ and $U_1U_2=U_2U_1=N_n.$
Let $B$ denote a Borel subgroup of the identity component of
$H$ and $\chi$ a character of $B.$  Assume
\begin{equation}\label{iw2}
B(F)H(\mathfrak o)=H(F).
\end{equation}
Let $V$ denote a space of functions on $N_n(F) \cdot H(F)$
which are compactly supported modulo $U_1(F)$ on the left  and satisfy
$$f(u_1bq) = \chi(b) f(q) \qquad \forall u_1\in U_1(F), \; b \in B(F),\; q \in H(F)N_n(F).$$
Let $V(N_{n}, \vartheta)$ denote the kernel of the usual projection from $V$ to
its twisted Jacquet module.

Let $$W(f)(q) = \int_{U_2(F)} f(u_2 q) \bar\vartheta (u_2) du_2.$$
Then
$Ker(W) \subset V(N_{n}, \vartheta).$
\end{lem}
\begin{proof}
We assume that
$$\int_{U_2(F)} f( uq) \bar\vartheta( u) du = 0,$$
for all $q\in H(F)N_n(F).$
What must be shown is that there is a compact subset $C$ of $N_{n}(F)$ such
that
$$\int_C f(gu) \bar\vartheta(u) du = 0,$$
for all $q\in H(F)N_n(F).$

Consider first $m \in H( \mathfrak o).$   Let $\mathfrak p$ denote the unique
maximal ideal in $\mathfrak o.$  If $U$ is a unipotent subgroup and $M$ an integer, we define
$$U(\primeideal^M)= \{u \in U(F): u_{ij} \in \primeideal^M \; \forall i \neq j \}.$$
Observe that for each $M\in \N,$ $N_{n}(\mathfrak p^M)$ is a subgroup of $N_{n}(F)$ which is
preserved by conjugation by  elements of $H( \mathfrak o).$
 We may choose $M$ sufficiently large that
$supp( f) \subset U_1(F) U_2( \mathfrak{p}^{-M} ) H(F).$  Then
we prove the desired assertion with $C= N_{n}( \mathfrak p^{-M}).$
Indeed, for  $m \in H( \mathfrak o),$  we have
$$\int_{N_{n}( \mathfrak p^{-M})} f(mu) \bar\vartheta(u) du
=\int_{N_{}( \mathfrak p^{-M})} f(um) \bar\vartheta(u) du,$$
because $Ad(m)$ preserves the subgroup $N_{n}(\mathfrak p^{-M}),$
and has Jacobian 1.  Let $c=\Vol( U_1(\mathfrak p^{-M})),$ which is
finite.  Then by $U_1$-invariance of $f,$ the above equals
$$ =c \int_{U_2( \mathfrak p^{-M})} f(um) \bar\vartheta(u) du.$$
This, in turn, is equal to
$$=c \int_{U_2(F)} f(um) \bar\vartheta(u) du,
$$
since none of the points we have added to the domain of integration
are in the support of $f,$ and this last integral is equal to zero by hypothesis.

Next, suppose $q = u_2m$ with $u_2 \in U_2(F)$ and $m\in H( \mathfrak o).$
If $u_2 \in  U_2(F)-U_2( \mathfrak p^{-M})$ then
$qu$ is not in the support of $f$ for any $u \in U_2( \mathfrak p^{-M}).$  On the other hand,
if $u_2 \in U_2( \mathfrak p^{-M}),$ then
$$\int_{N_{n}( \mathfrak p^{-M})} f(u_2mu) \bar\vartheta(u) du
=\int_{N_{n}( \mathfrak p^{-M})} f(u_2u m) \bar\vartheta(u) du$$
$$= \vartheta( u_2 ) \int_{N_{n}( \mathfrak p^{-M})} f(u m) \bar\vartheta(u) du,$$
and now we continue as in the case $u_1=1.$

The result for general $q$ now follows from the left-equivariance properties of $f$
and \eqref{iw2}.
\end{proof}

\section{Appendix VII:  Identities of Unipotent Periods}
\label{s: appV}

\begin{lem}\label{u1u2lemma}
Let $(U_1^a, \psi_1^a)$  and $(U_2, \psi_2^a)$ be defined as in Theorem \ref{t:maintheorem}.
Then $(U_1^a, \psi_1^a) \sim (U_2, \psi_2^a),$ for all $a\in F.$
\end{lem}
\begin{proof}
We regard $a$ as fixed and omit it from the notation.
We define some additional unipotent periods which appear at intermediate
stages in the argument.
Let $U_4$ be the subgroup defined by $u_{n, j}=0$ for $j = n+1$ to $2n-1$ and
$u_{2n,2n+1}=0.$
We define  a character $\psi_4$ of $U_4$ by the same formula as
$\psi_1.$  Then $(U_1,\psi_1)$ may be swapped for $(U_4, \psi_4).$
(See definition \ref{d:swap}.)

Now, for each $k$ from $1$ to $n,$ define $(U_5^{(k)},\psi_5^{(k)})$ as follows.
 First, for each $k,$ the group $U_5^{(k)}$ is contained in the subgroup of $\maxunip$
defined by, $u_{2n,2n+1}=0.$
In addition,
$u_{n+k-1,j}= 0$ for $j < 2n,$ and $u_{i,i+1}=0$ if $n-k+1 \leq i < n+k-1$ and $i \equiv n-k+1 \mod 2,$
and
$$\psi_5^{(k)}(u)=\baseaddchar\left(
\sum_{i=1}^{n-k} u_{i,i+1} + \sum_{i=n-k+1}^{n+k-2} u_{i,i+2}
+ u_{n+k-1, 2n}+ \frac a2u_{n+k-1, 2n+2} + \sum_{i=n+k}^{2n-2} u_{i,i+1}+u_{2n-1,2n+2}
\right).$$
(Note that one or more of the sums here may be empty.)

  Next, let $U_6^{(k)}$  be the subgroup of $\maxunip$ defined by the conditions
 $u_{2n,2n+1}=0,$
$u_{n+k-1,j}= 0$ for $j < 2n,$ and $u_{i,i+1}=0$ if $n-k+1 \leq i < n+k-1$ and $i \equiv n-k \mod 2.$
    The same formula which defines $\psi_5^{(k)}$
  also defines a character of $U_6^{(k)}.$  We denote this character
  by $\psi_6^{(k)}.$

  We make the following observations:
  \begin{itemize}
  \item{$(U_5^{(1)},\psi_5^{(1)})$ is precisely $(U_4, \psi_4).$}
  \item{For each $k,$ $(U_5^{(k)},\psi_5^{(k)})$ is conjugate to $(U_6^{(k+1)}, \psi_6^{(k+1)}).$
  The conjugation is accomplished by any preimage of the
  permutation matrix which transposes $i$ and $i+1$
  for
  $n-k \leq i < n+k$ and $i \equiv n-k+1 \mod 2.$
  }
  \item{
  $(U_6^{(k)} ,\psi_6^{(k)})$ may be swapped for  $(U_5^{(k)}, \psi_5^{(k)}).$  }
  \end{itemize}

 Thus $(U_4, \psi_4) \sim (U_5^{(n-1)}, \psi_5^{(n-1)}).$

 Now, let  $U_2'=U_5^{(n-1)},$ and let
$$\psi_2'(u)
 = \baseaddchar(u_{1,3} + \dots + u_{2n-2, 2n} +u_{2n-2, 2n+1} +\frac a2 u_{2n-1, 2n}+u_{2n-1,2n+2}).$$
  Then $(U_5^{(n-1)},\psi_5^{(n-1)})$ is conjugate to $(U_2', \psi_2'),$ which may be swapped
  for $(U_2,\psi_2).$
  \end{proof}

\begin{lem}\label{u2u3deep}
Let $(U_3, \psi_3)$ and $(U_2, \psi_2^0)$  be defined as in Theorem \ref{t:maintheorem}.
  Then
$$( U_3, \psi_3) \in \langle ( U_2, \psi_2^0) ,
\{ (N_\ell, \vartheta): n \leq \ell < 2n \text{ and } \vartheta \text{ in general position.} \}\rangle.$$
\end{lem}
\begin{proof}
To prove this assertion we introduce some additional unipotent periods.
For $k= n$ to $2n-1$ let $U_7^{(k)}$ denote the subgroup of $\maxunip$ defined
by $u_{2n,2n+1}=0,$ and $u_{i,2n}=0$ for $k+1\le i \le 2n-1,$ and let
 $$\psi_7^{(k)}(u)= \baseaddchar
\left( \sum_{i=1}^{k-1} u_{i,i+1}
+ u_{k,2n}
+ \sum_{i=k+1}^{2n-2} u_{i,i+1}+u_{2n-1,2n+2}\right).$$
Let $U_8^{(k)}$ denote the subgroup defined by
by $u_{2n-1,2n+1}=0,$  $u_{k,j}=0$ for $k+1\le j < 2n,$ and
let $U_9^{(k)}$ denote the subgroup defined by the additional condition
$u_{k,2n}=0.$
The same formula which defines $\psi_7^{(k)}$ may be used to specify a character of
$U_8^{(k)},$ which we denote $\psi_8^{(k)}.$
In addition, let
 $$\baseaddchar
\left( \sum_{i=1}^{k-1} u_{i,i+1}
+ u_{k,2n+2}
+ \sum_{i=k+1}^{2n-1} u_{i,i+1}\right),$$
be denoted by  $\tilde \psi_8^{(k)}$ for $u\in U_8^{(k)}$ or  $\psi_9^{(k)}$
for  $u\in U_9^{(k)}.$

Now, we need the following observations:
\begin{itemize}
\item $(U_7^{(n)},\psi_7^{(n)})$ is just the period $(U_1^0,\psi_1^0)$ from
theorem \ref{t:maintheorem}, and so is equivalent to $(U_2,\psi_2^0)$ by the previous
result.  \item For each $k,$ $(U_7^{(k)},\psi_7^{(k)})$ is conjugate
to $(U_9^{(k+1)},\psi_9^{(k+1)}).$  (One conjugates by a preimage of a permutation
matrix and then by a toral element to fix a minus sign which is introduced.)
\item $(U_8^{(k+1)},\tilde \psi_8^{(k+1)})$ is spanned by $(U_9^{(k+1)},\psi_9^{(k+1)})$
and $\{(N_k,\vartheta): \vartheta \text{ in general position}\}.$  More precisely, if $\vartheta$
is any extension of $\psi_9^{(k+1)}$ which is {\em not} in general position, then the restriction
of $\vartheta$ to $U_8$ is $\tilde \psi_8^{(k+1)}$. (Cf. Corollary \ref{unipercor}.)
\item $(U_8^{(k)},\tilde \psi_8^{(k)})$ is conjugate to $(U_8^{(k)}, \psi_8^{(k)}).$
\item $(U_8^{(k)}, \psi_8^{(k)})$ may be swapped for $(U_7^{(k)},\psi_7^{(k)}).$
\end{itemize}
We deduce that $(U_2,\psi_2^0)$ divides $(U_8^{(2n-1)},\psi_8^{(2n-1)}),$ a period which
differs from $(U_3,\psi_3)$ only in that integration over $u_{2n,2n+1}$ is omitted.  Thus
$(U_3,\psi_3)$ is the constant term in the Fourier expansion of  $(U_8^{(2n-1)},\psi_8^{(2n-1)}),$
in the variable $u_{2n,2n+1},$ while all of the nonconstant terms are Whittaker integrals
with respect to various generic characters of $\maxunip.$  As $\residuerep$ is non-generic,
they all vanish.  The result follows.
 \end{proof}

\begin{lem}  Take $a \in F^\times.$  We regard $a$ as fixed throughout and, for the most part
we suppress it
from the notation.
 As in Theorem \ref{t:maintheorem}, let
 $V_i$ denote the unipotent radical of the standard parabolic of $G_{4n+1}$ having Levi
 isomorphic to $GL_i \times G_{4n-2i+1}$ (for $1 \leq i \leq 2n$).
 For $1\le j < 2n,$
 let $V_i^{4n-2j}$ denote the unipotent radical of the standard maximal
 parabolic of $G^a_{4n-2j}$
 having Levi isomorphic to $GL_i \times G_{4n-2j-2i}^a$ (for $1 \leq i \leq 2n-j-2$ in
 the split case and $1 \leq i \leq 2n-j-2$ in the nonsplit cases).
 Let  $(N_\ell, \Psi_\ell^a)$ be the period used to define the descent, as usual,
 and let
 $(N_\ell, \Psi_\ell^a)^{(4n-2k+1)}$  denote the analogue for $G_{4n-2k+1},$ embedded
 into $G_{4n+1}$ inside the Levi of a maximal parabolic.

 Then,
\label{l:cuspidality-unip-id}
$(V_k^{2n }, {\bf 1} ) \circ ( N_{n} , \Psi_{n})$
is an element of $$ \langle
(N_{n+k} , \Psi_{n+k}) , \{
(N_{n+ j}, \Psi_{n+j})^{(4n-2k+2j)}
\circ ( V_{k-j} , {\bf 1} ): \; \; 1 \leq j < k
\}
\rangle.$$
\end{lem}
\begin{proof}
Let $m=(m_1,m_2,m_3)$ be a triple of integers satisfying:
$0\le m_1 < m_2\le m_3+1\le 2n.$  We associate to this data a unipotent group $U_m$
and two characters $\psi_m,\psi_m'$ as follows:
\begin{itemize}
\item $U_m$ is defined by the condition that $u_{i,j}=0$ whenever $m_1 < i < m_2-1$ and
$j< m_2,$ or $m_3 < i,$
\item $\psi_m(u) = \baseaddchar\left(\sum_{i=1}^{m_1} u_{i,i+1} + u_{m_1+1, m_2}+\sum_{i=m_2}^{m_3-1}
u_{i,i+1} + u_{m_3,2n} +\frac a2 u_{m_3, 2n+2}
\right),$
\item $\psi_m'(u) = \baseaddchar\left(\sum_{i=1}^{m_1-1} u_{i,i+1} + u_{m_1, m_2-1}+
\sum_{i=m_2-1}^{m_3-1}
u_{i,i+1} + u_{m_3,2n} +\frac a2 u_{m_3, 2n+2}
\right).$
\end{itemize}
Then $(U_m,\psi_m')$ is conjugate to $(U_m,\psi_m)$ and may be swapped for
$(U_{m'},\psi_{m'}),$ where
$(m_1,m_2,m_3)' = (m_1-1,m_2-1,m_3).$
Furthermore, for any $k<n,$
$(V_k^{2n }, {\bf 1} ) \circ ( N_{n} , \Psi_{n})$
is an integral over the subgroup of $U_{n,n+k+1,n+k}$ defined by the conditions,
$u_{i,2n}=-\frac a2u_{i,2n+2},$ for $n < i \le n+k.$  It may be swapped for the period
$(U_m,\psi'')$ corresponding to $m=(n-1,n+k+1,n+k),$ and
$$\psi''(u) =
\baseaddchar\left(\sum_{i=1}^{n-1} u_{i,i+1} +  u_{n,2n} +\frac a2 u_{n, 2n+2}
\right),$$
and this period is conjugate to $(U_m,\psi_m')$ for this value of $m.$
It follows that
$(V_k^{2n }, {\bf 1} ) \circ ( N_{n} , \Psi_{n})$
is equivalent to $(U_m,\psi_m')$ for the triple $m=(0,k+2, n+k).$

Now, it's easy to see that
$(U_{(0,1,m_3)},\psi_{(0,1,m_3)}') =(N_{m_3},\Psi_{m_3}^a),$ and that for $m_2>2$
there are two orbits of extensions of
$\psi_{(0,m_2,m_3)}$ to $U_{(0,m_2-1,m_3)},$ namely, the one containing
$\psi_{(0,m_2-1,m_3)}',$ and the trivial extension, which yields the period
$(N_{m_3-m_2+2},\psi_{m_3-m_2+2}^a)^{(4n-2m_2+5)}\circ (V_{m_2-2},{\bf 1}).$
This proves the assertions regarding all cases except for the two parabolics with Levi
isomorphic to $GL_1\times GL_n$ in the split case.

As noted previously, it is enough to consider one of them, because they are conjugate in $G_{4n+1}.$
Furthermore, we may conjugate by $h_a,$ and use the more convenient embedding
of $G_{2n}^\square$ into $G_{4n+1}$ as $(L_n^{\Psi_n})^0.$

For this case we take $m\in \Z$ with $0\le m \le n,$
and define $U_m$ to be the subgroup of $\maxunip$
defined by $u_{i,j}=0$ whenever $m<i<j\le m+n+1.$
Take
$$\psi_m'(u) = \baseaddchar\left( \sum_{i=1}^{m-1}u_{i,i+1} +u_{m,m+n+1}
+\sum_{i=m+n+2}^{2n}u_{i,i+1}
\right),$$
$$\psi_m''(u) = \baseaddchar\left( \sum_{i=1}^{m}u_{i,i+1} +u_{m+1,m+n+2}
+\sum_{i=m+n+3}^{2n}u_{i,i+1}
\right).$$
Then $(V_n^{2n},{\bf 1}) \circ (N_n,\Psi_n) = (U_n,\psi'_n).$  Furthermore
$(U_m,\psi'_m)$ is conjugate to $(U_m,\psi_m'')$ and may be swapped for
$(U_{m-1},\psi_{m-1}'').$
Furthermore,
$(U_0, \psi'_0)$ is easily seen to be in the span of the periods
$$(\maxunip^{4n-2k+1},\vartheta) \circ (V_k, {\bf 1})$$ for $0\le k < n$
and $\vartheta$ a generic character of the maximal unipotent subgroup
of $G_{4n-2k+1}$ (embedded into $G_{4n+1}$) as a component of a standard
Levi as usual.
This completes the proof.
\end{proof}

So far, we have proved relations of two forms
\begin{itemize}
\item Equivalencies, in which the unipotent subgroup $U$ is replaced by another of the
same dimension, and the character $\psi$ by another in the same orbit.
\item Relations where one replaces $U$ by a group of properly larger dimension, and
considers all orbits of extensions of $\psi.$
\end{itemize}

The statement that $(U_2,\psi_2^0)$ is spanned by $\{(U_2,\psi_2^a): \; a \in F^\times\}$
is of a different nature,
and requires the use of theta functions, as in section
\ref{section with theta functions}.
\begin{lem} \label{thetaLemma}
Let the group $U_2,$ and the character $\psi_2^a$ for each $a\in F$
be defined as in the main theorem.  Then
$(U_2,\psi_2^0)\in \langle \{(U_2,\psi_2^a):a\in F^\times \}\rangle.$
\end{lem}
\begin{proof}
The $R=LN$ be the unique standard parabolic subgroup
of $GSpin_{4n+1}$
such that the Levi, $L$ is isomorphic to
$GL_2^{n-1} \times GSpin_5.$ Define a
character $\psi_N$ of the unipotent radical $N$
by
$$
\psi_N(u) = \baseaddchar\left(
\sum_{i=1}^{2n-2} u_{i, i+2}
\right).
$$
Let $\on{Stab}_L(\psi_N)$ denote
the stabilizer of $\psi_N$ in $L.$
Then  $\on{Stab}_L(\psi_N)$
is equal to the product of
a reductive group isomorphic to $GL_2 \times GL_1$
and three dimensional unipotent group.
The image in $SO_{4n+1}$ consists
of matrices of the form
$$
\diag( g, g, \dots, g, g', \,_tg^{-1}
\dots , \,_tg^{-1}), \;
g \in GL_2, \;
g' = \bpm g&*&*\\&1&*\\&&_tg^{-1} \epm \in SO_5.
$$
In particular, $\on{Stab}_L(\psi_N)$ maps isomorphically
onto the Siegel parabolic of $GSpin_5,$ which is to say
the Klingen parabolic of $GSp_4.$  This group
has a subgroup which was identified with $G^J$ above.
Now
$$
\varphi^{(U_2, \psi_2^a)}
= \left( \varphi^{(N, \psi_N)}\right)\FC\Sieg a,
$$
so this result follows from corollary \ref{theta lemma Sp4 second cor}.
\end{proof}

\end{document}